\documentclass[american]{IEEEtran}
\usepackage[T1]{fontenc}
\usepackage[latin9]{inputenc}
\usepackage{float}
\usepackage{calc}
\usepackage{amsmath}
\usepackage{amssymb}
\usepackage{graphicx}

\makeatletter

\providecommand{\tabularnewline}{\\}
\floatstyle{ruled}
\newfloat{algorithm}{tbp}{loa}
\providecommand{\algorithmname}{Algorithm}
\floatname{algorithm}{\protect\algorithmname}

\usepackage{url}
\usepackage{subfigure}
\usepackage{epsfig}
\usepackage{subfigmat}
\usepackage{amsthm}\usepackage{tensor}\usepackage{amsfonts}
\usepackage{mathrsfs}

\theoremstyle{plain}
\newtheorem{theorem}{\bf Theorem}
\newtheorem{corollary}{\bf Corollary}

\theoremstyle{remark}
\newtheorem{assumption}{\bf Assumption}
\newtheorem{definition}{\bf Definition}
\newtheorem{remark}{\bf Remark}

\renewenvironment{proof}[1][Proof:]{\begin{trivlist}
\item[\hskip \labelsep {\bfseries #1}]}{\end{trivlist}}

\makeatother

\usepackage{babel}
\begin{document}

\title{Probabilistic and Distributed Control of a Large-Scale Swarm of Autonomous
Agents}

\author{Saptarshi Bandyopadhyay \IEEEmembership{Member,~IEEE}, Soon-Jo Chung
\IEEEmembership{Senior~Member,~IEEE}, and \\*Fred Y. Hadaegh \IEEEmembership{Fellow,~IEEE}
\thanks{S. Bandyopadhyay and F. Y. Hadaegh is with the Jet Propulsion Laboratory,
California Institute of Technology, Pasadena, CA 91109, USA (e-mail:
saptarshi.bandyopadhyay@jpl.nasa.gov and fred.y.hadaegh@jpl.nasa.gov).
S.-J. Chung is with the University of Illinois at Urbana-Champaign,
Urbana, IL 61801, USA. (email: sjchung@alum.mit.edu). This research
was supported in part by AFOSR grant FA95501210193 and NSF IIS 1253758.}}
\maketitle
\begin{abstract}
We present a novel method for guiding a large-scale swarm of autonomous
agents into a desired formation shape in a distributed and scalable
manner. Our Probabilistic Swarm Guidance using Inhomogeneous Markov
Chains (PSG\textendash IMC) algorithm adopts an Eulerian framework,
where the physical space is partitioned into bins and the swarm's
density distribution over each bin is controlled. Each agent determines
its bin transition probabilities using a time-inhomogeneous Markov
chain. These time-varying Markov matrices are constructed by each
agent in real-time using the feedback from the current swarm distribution,
which is estimated in a distributed manner. The PSG\textendash IMC
algorithm minimizes the expected cost of the transitions per time
instant, required to achieve and maintain the desired formation shape,
even when agents are added to or removed from the swarm. The algorithm
scales well with a large number of agents and complex formation shapes,
and can also be adapted for area exploration applications. We demonstrate
the effectiveness of this proposed swarm guidance algorithm by using
results of numerical simulations and hardware experiments with multiple
quadrotors.
\end{abstract}


\section{Introduction \label{sec:Introduction}}

\begin{figure}
\begin{centering}
\includegraphics[bb=75bp 65bp 890bp 470bp,clip,width=3in]{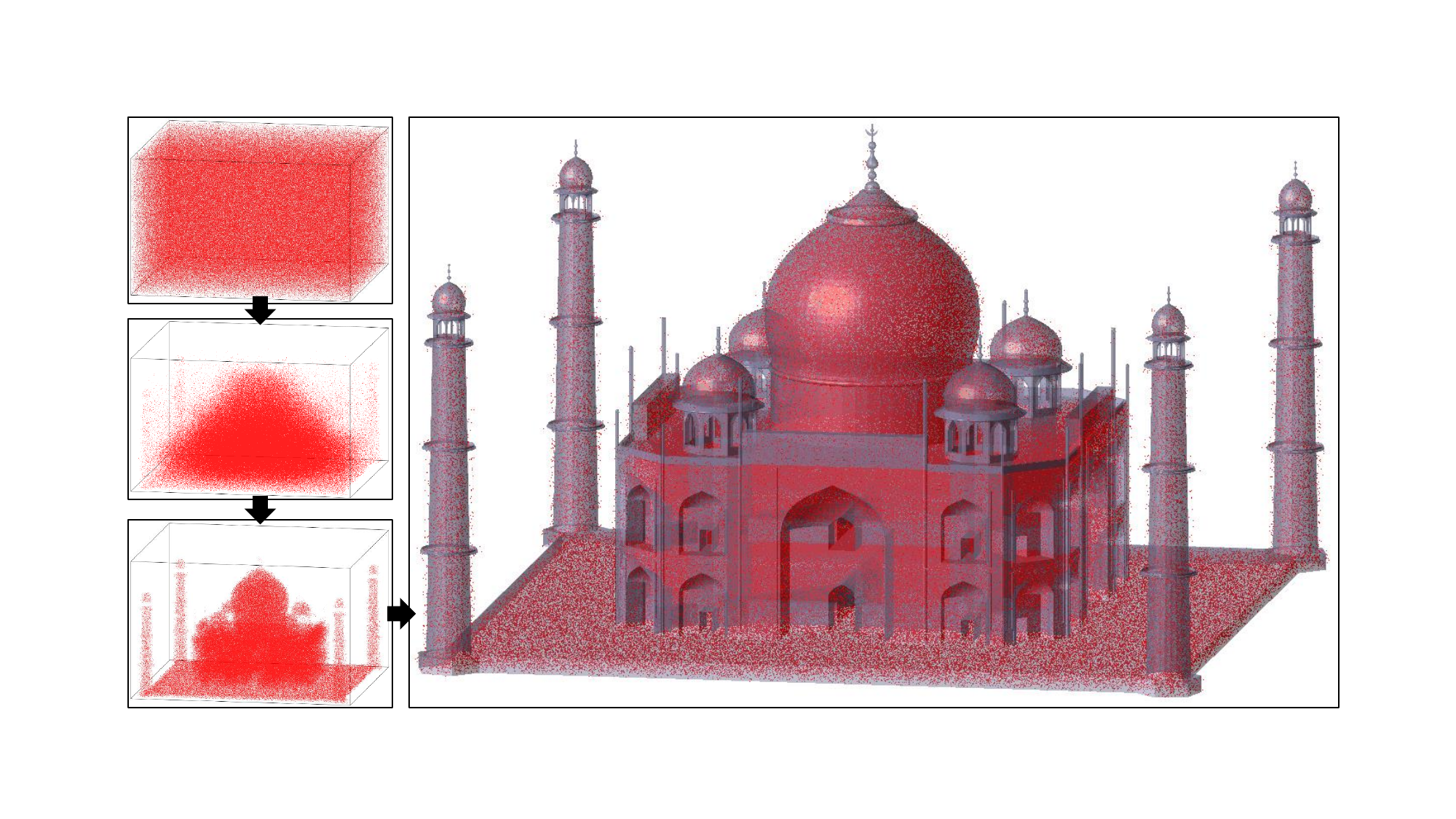}
\par\end{centering}
\caption{Using the PSG\textendash IMC algorithm for shape formation, a million
swarm agents (shown in red) attain the complex 3-D shape of the Taj
Mahal (translucent silhouette shown in gray). The physical space is
partitioned into $100\times100\times70$ bins. See the supplementary
video (SV1). \label{fig:swarm-Taj-Mahal}}
\end{figure}

A large-scale swarm of space robotic systems could collaboratively
complete tasks that are very difficult for a single agent, with significantly
enhanced flexibility, adaptability, and robustness \cite{Ref:Hadaegh13}.
Moreover, a large-scale swarm (having $10^{3}$\textendash $10^{6}$
or more agents) may be deployed for challenging missions like constructing
a complex formation shape (see Fig.~\ref{fig:swarm-Taj-Mahal}) \nocite{Ref:Nagpal14,Ref:Kushleyev13,Ref:Dorigo06,Ref:Lynch08}\nocite{Ref:Chung12,Ref:Morgan15_SATO,Ref:Morgan14}\cite{Ref:Nagpal14}\textendash \cite{Ref:Morgan14}
or exploring an unknown environment \nocite{Ref:Bullo04,Ref:Correll07,Ref:Dorigo12}\nocite{Ref:Bruemmer02,Ref:Hurtado2004,Ref:Liu07}\cite{Ref:Bullo04}\textendash \cite{Ref:Liu07}.
The control algorithm for such a large-scale swarm of autonomous agents
should be: 
\begin{itemize}
\item Distributed: The algorithm should not need a centralized supervisor
or controller.
\item Versatile: The algorithm can be easily tailored for multiple applications
such as maintaining the formation shape or exploring the target area. 
\item Robust: Since a small fraction of agents in the swarm might get lost
during the course of an operation or new agents might get added to
the swarm, the algorithm should seamlessly adapt to loss or addition
of agents. Moreover, the algorithm should accommodate measurement
errors, actuation errors, and other uncertainties. 
\item Scalable: The algorithm should scale well with the number of agents
and the size of the area.
\end{itemize}
In this paper, we lay the theoretical foundations of a distributed,
versatile, robust, and scalable algorithm for controlling the shape
of large-scale swarms. 

In fluid mechanics, the motion of fluids is described using either
a \textit{Lagrangian} framework, where each fluid particle's motion
is tracked; or an \textit{Eulerian} framework, where the fluid's density
distribution at specific locations is tracked. Analogously, distributed
control algorithms for swarms can be classified into two categories
\cite{Ref:Okubo94,Ref:Passino04}: the individual-agent-based Lagrangian
framework, where each agent's trajectory is generated separately \cite{Ref:Nagpal14}\textendash \cite{Ref:Liu07};
and the continuum-based Eulerian framework, where the collective properties
of the swarm (e.g., its density distribution) is controlled. In the
Lagrangian framework, the computation cost for generating each agent's
trajectory increases rapidly with a large number ($10^{3}$\textendash $10^{6}$
or more) of agents. For example, the computational complexity of each
agent's target assignment problem increases at least quadratically
with the number of agents \cite{Ref:Lynch08}\textendash \cite{Ref:Morgan14}.
The Eulerian framework decouples these two tasks by first solving
the target assignment problem at a coarser spatial resolution. Moreover,
the Lagrangian framework does not efficiently handle loss or addition
of agents, nor does it scale well with the size of the area and arbitrary
formation shapes \cite{Ref:Nagpal14}\textendash \cite{Ref:Dorigo06}.
Consequently, we adopt the Eulerian framework in this paper.

\subsection{Literature Review}

Numerous path planning algorithms within the Lagrangian framework
are discussed in the survey papers on swarm robotics \cite{Ref:Sahin05,Ref:Brambilla13}.
In this subsection, we focus on guidance algorithms that use an Eulerian
approach \nocite{Ref:Canuto08,Ref:Menon04,Ref:Egerstedt10}\cite{Ref:Canuto08}\textendash \cite{Ref:Egerstedt10}.
For shape formation and reconfiguration applications, the physical
space over which the swarm is distributed is first tessellated (i.e.,
partitioned) into discrete bins \cite{Ref:Brooks85,Ref:Tsiotras07}.
The bin size is determined by the spatial resolution of the desired
formation shape. Assuming that the number of agents is much larger
than the number of non-empty bins, the density distribution of the
swarm over these bins is controlled to achieve the desired formation
shape. 

Within the Eulerian framework, Homogeneous Markov Chain (HMC) based
algorithms are a popular choice for shape formation \nocite{Ref:Chattopadhyay09,Ref:Acikmese12,Ref:Acikmese14_AJC,Ref:Acikmese15_AR}\cite{Ref:Chattopadhyay09}\textendash \cite{Ref:Acikmese15_AR},
area exploration \cite{Ref:Hereford10,Ref:Hespanha08}, task allocation
\cite{Ref:Kumar09,Ref:Mather12}, and surveillance applications \cite{Ref:Grace05,Ref:Bullo14}.
In such algorithms, the agent's transition probability between bins
is encoded in a constant Markov matrix that has the desired formation
shape as its stationary distribution. Such an approach is probabilistic,
as opposed to deterministic, because each agent determines its next
bin location by inverse-transform sampling of the Markov matrix \cite{Ref:Devroye86}.
Since the number of agents in the swarm can vary with time and the
agents do not keep track of the time-varying number of agents in the
swarm, a probabilistic approach is preferred because a deterministic
target assignment algorithm needs to keep track of the changes in
the number of agents and targets \cite{Ref:Morgan15_SATO}. Moreover,
as shall be seen in the paper, our probabilistic approach can also
handle measurement uncertainties and actuation errors in a robust
manner. The HMC-based algorithms possess the aforementioned benefits
of robustness and scalability, because addition or removal of agents
from the swarm does not affect the convergence of the HMC to the stationary
distribution. 

However, the major drawback of these HMC-based algorithms is that
they are inherently open-loop strategies which cannot incorporate
any feedback. Clearly, the effectiveness of these algorithms can be
greatly improved by refining the Markov matrix at each time step using
some feedback. Such refinement results in an Inhomogeneous Markov
Chain (IMC), which is at the core of our algorithm. 

In this paper, we derive the Probabilistic Swarm Guidance using Inhomogeneous
Markov Chains (PSG\textendash IMC) algorithm, which incorporates the
feedback from the current swarm density distribution at each time
step. The PSG\textendash IMC algorithm is a closed-loop distributed
guidance strategy that retains the original robustness and scalability
properties associated with a Markovian approach. Another disadvantage
of HMC-based algorithms is that they suffer undesirable transitions,
i.e., transitions from bins that are deficient in agents to bins with
surplus agents. Such undesirable transitions prevent the swarm from
converging to the desired formation. The PSG\textendash IMC algorithm
suppresses such undesirable transitions between bins, thereby reducing
the control effort needed for achieving and maintaining the formation.
This benefit also results in smaller convergence errors than HMC-based
algorithms. 

The motion of swarm agents can be formulated as an optimal transport
problem with respect to a given cost function \cite{Ref:Bertsekas91}.
If perfect feedback of the current swarm distribution is available,
then our prior work on the Probabilistic Swarm Guidance using Optimal
Transport (PSG\textendash OT) algorithm \cite{Ref:Bandyopadhyay14MSC}
gives good performance. However, there are two major disadvantages
of such an approach. First, the performance of an optimal transport-based
algorithm drops precipitously with estimation errors in the feedback
loop. Measurement and estimation errors are routinely encountered
in practice and it is often impossible or impractical to generate
perfect feedback of the current swarm distribution. Second, the computation
time of the optimization problem increases very fast with an increasing
number of bins. This is a notable drawback because a large number
of bins are necessary for capturing fine spatial details in the desired
formation shape. The PSG\textendash IMC algorithm proposed in the
present paper can overcome both challenges, since it works effectively
in the presence of error-prone feedback and scales well with a large
number of bins.

A different approach to swarm formation within the Eulerian framework
is to model the continuum dynamics using partial differential equations
(PDEs) \nocite{Ref:Milutinovic06,Ref:Ferrari14,Ref:Krstic15}\cite{Ref:Milutinovic06}\textendash \cite{Ref:Krstic15}.
Since our goal is to achieve arbitrary formation shapes that are not
limited to to the equilibrium states of the PDE, we do not consider
a PDE-based approach. Our approach is also different from the multiagent
Markov decision process approach in \cite{Ref:Durfee05,Ref:How14}
because our agents do not keep track of the states and actions of
other agents.

Distributed Eulerian approaches for area exploration applications
using region-based shape controllers and attraction-repulsion forcing
functions are discussed in \nocite{Ref:Mogilner99,Ref:Slotine09,Ref:MKumar11}\cite{Ref:Mogilner99}\textendash \cite{Ref:MKumar11}.
We show that a slight modification of our PSG\textendash IMC algorithm
also results in an efficient area-exploration algorithm. 

\subsection{Main Contributions}

The first contribution of this paper is a novel technique for constructing
feedback-based Markov matrices for a given stationary distribution
which represents the desired formation, where the expected cost of
transitions at each time instant is minimized. Each Markov matrix
satisfies the motion constraints that might arise due to the dynamics
or other physical constraints. The Markov matrix converges to the
identity matrix when the swarm converges to the desired formation,
thereby ensuring that the swarm settles down and reducing unnecessary
transitions. 

Second, we describe the PSG\textendash IMC algorithm for shape formation.
We rigorously derive the convergence proofs for the PSG\textendash IMC
algorithm based on the analysis of IMC, which is more involved than
the convergence proof for HMC. We show that each agent's IMC strongly
ergodically converges to the desired formation shape. Further, we
provide a time-varying probabilistic-bound on the convergence error
as well as a lower bound on the number of agents for ensuring that
the final convergence error is below the desired threshold. Furthermore,
we present an extension of the PSG\textendash IMC algorithm for area
exploration applications. 

Along with a lower-level collision-free motion planner, we demonstrate
using multiple quadrotors that the PSG\textendash IMC algorithm can
be executed in real-time to achieve multiple desired formation shapes.
Using results of numerical simulations, we show that the PSG\textendash IMC
algorithm yields a smaller convergence error and more robust convergence
than the HMC-based and PSG\textendash OT algorithms, while significantly
reducing the number of transitions in the presence of estimation errors.
Thus, the PSG\textendash IMC algorithm is best-suited for large-scale
swarms with error-prone feedback and complex desired formations with
a large number of bins.

Compared to our conference paper \cite{Ref:Bandyopadhyay_IROS16},
we have added detailed proofs in the convergence analysis section
and extensions of our algorithm for shape formation and area exploration
applications. We have also added extensive numerical and experimental
results to this paper. 

This paper is organized as follows. The problem statement for shape
formation is discussed in Section \ref{sec:Problem-Statement}. In
Section \ref{sec:Construction-of-Markov-Matrix}, we present our novel
techniques for constructing Markov matrices. The PSG\textendash IMC
algorithm for shape formation is presented in Section \ref{sec:PSG-IMC-pattern-formation}.
The PSG\textendash IMC algorithm for area exploration is presented
in Section \ref{sec:Goal-Searching}. Results of numerical simulations
and experimentation are presented in Section \ref{sec:Numerical-Simulation-pattern-formation}.
This paper is concluded in Section \ref{sec:Conclusion}.

\begin{table}[h]
\caption{List of frequently used symbols \label{tab:List-of-symbols}}

\centering{}%
\begin{tabular}{ll}
\hline 
Symbol & Explanation\tabularnewline
\hline 
\hline 
$\boldsymbol{A}_{k}^{j}$ & Motion constraints matrix of the $j^{\textrm{th}}$ agent (Definition~\ref{def:Motion-Constraints})\tabularnewline
$B[i]$ & Bins, where $1\leq i\leq n_{\textrm{bin}}$ (Definition~\ref{def:bins})\tabularnewline
$\boldsymbol{M}_{k}^{j}$ & Markov matrix of the $j^{\textrm{th}}$ agent (\ref{eq:Markov_matrix})\tabularnewline
$\boldsymbol{S}_{k}^{j}$ & Condition for escaping transient bins (\ref{eq:S_matrix})\tabularnewline
$\boldsymbol{C}_{k}$ & Cost matrix at the $k^{\textrm{th}}$ time instant (Definition~\ref{def:cost-function})\tabularnewline
$m_{k}$ & Number of agents in the swarm (Assumption~\ref{assump:m_k})\tabularnewline
$\boldsymbol{r}_{k}^{j}$ & Actual bin position of the $j^{\textrm{th}}$ agent (Assumption~\ref{assump:r_known})\tabularnewline
$\boldsymbol{x}_{k}^{j}$ & Probability vector of the $j^{\textrm{th}}$ agent predicted position
(\ref{eq:def_prob_vector})\tabularnewline
$\boldsymbol{\Theta}$ & Probability vector of the desired formation (Definition~\ref{def:Desired-formation})\tabularnewline
$\epsilon_{\mathrm{est}}$ & Estimation error between $\boldsymbol{\mu}_{k}^{\star}$ and $\boldsymbol{\mu}_{k}^{j}$
(\ref{eq:estimation_error})\tabularnewline
$\boldsymbol{\mu}_{k}^{\star}$ & Current swarm distribution (Definition~\ref{def:swarm-distribution}) \tabularnewline
$\boldsymbol{\mu}_{k}^{j}$ & Estimate of $\boldsymbol{\mu}_{k}^{\star}$ by the $j^{\textrm{th}}$
agent\tabularnewline
$\xi_{k}^{j}$ & Hellinger-distance-based feedback gain (Definition~\ref{def:Feedback-gain}) \tabularnewline
\hline 
\end{tabular}
\end{table}

\subsection{Notations}

The \textit{time index} is denoted by a right subscript and the \textit{agent
index} is denoted by a lower-case right superscript. Frequently used
symbols are listed in Table \ref{tab:List-of-symbols}.

The symbol $\mathbb{P}(\cdot)$ refers to the probability of an event.
Let $\mathbb{N}$ and $\mathbb{R}$ be the set of natural numbers
(positive integers) and real numbers respectively. The matrix $\mathrm{diag}(\boldsymbol{\alpha})$
denotes the diagonal matrix of appropriate size with the vector $\boldsymbol{\alpha}$
as its diagonal elements. Let $\mathbf{1}=[1,1,\ldots,1]^{T}$, $\mathbf{I}$,
$\mathbf{0}$, and $\emptyset$ denote the ones (column) vector, the
identity matrix, the zero matrix of appropriate sizes, and the empty
set respectively. Let $\|\cdot\|_{p}$ denote the $\ell_{p}$ vector
norm. Let $\min{}^{+}$ denote the minimum of the positive elements. 

A probability vector $\boldsymbol{a}\in\mathbb{R}^{n}$ is a row vector
with non-negative elements, where the sum of all the elements is 1
(i.e., $\boldsymbol{a}\geq0$, $\boldsymbol{a}\mathbf{1}=1$) \cite[pp. 92]{Ref:Seneta06}.
The metric $D_{(\cdot)}(\boldsymbol{a},\boldsymbol{b})$ represents
the distance between probability vectors $\boldsymbol{a}$ and $\boldsymbol{b}$,
where the subscript represents the type of metric.

\section{Preliminaries and Problem Statement \label{sec:Problem-Statement}}

In this section, we state the problem statement for shape formation
in Section~\ref{subsec:Problem-Statement-pattern-formation} and
introduce the PSG\textendash IMC algorithm for shape formation in
Section~\ref{subsec:PSG-IMC-flowchart}.

\begin{figure}[h]
\begin{centering}
\includegraphics[width=1.5in]{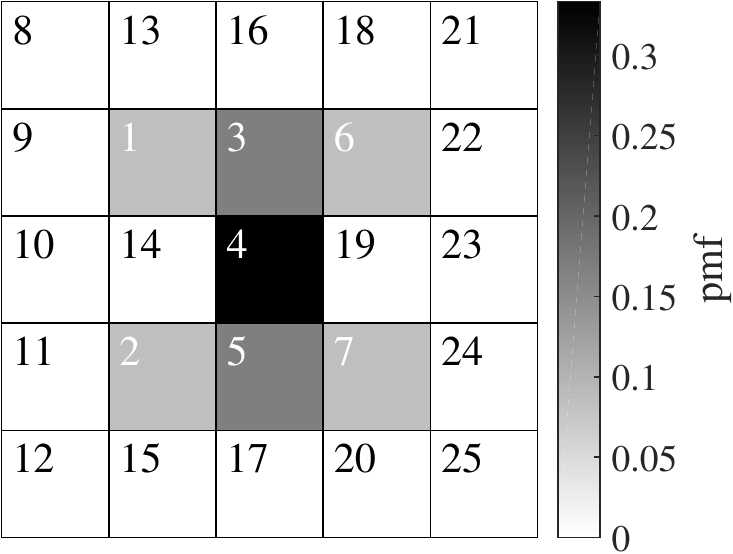}
\par\end{centering}
\caption{In this example, the desired formation $\boldsymbol{\Theta}=[\frac{1}{12},\frac{1}{12},\frac{1}{6},\frac{1}{3},\frac{1}{6},\frac{1}{12},\frac{1}{12},\mathbf{0}^{1\times18}]$.
The bins 1 to 7 are recurrent bins. \label{fig:I-example}}
\end{figure}

\begin{definition} \label{def:bins} \textit{(Bins $B[i]$)} The
compact physical space over which the swarm is distributed is partitioned
into $n_{\textrm{bin}}$ disjoint bins. These bins are represented
by $B[i]$ for all $i\in\{1,\ldots,n_{\textrm{bin}}\}$. The size
of the bins is determined by the spatial resolution of the desired
formation shape. For example,  $n_{\textrm{bin}}=25$ in Fig.~\ref{fig:I-example}.
\hfill$\Box$ \end{definition}

\begin{definition} \label{def:Desired-formation} \textit{(Desired
Formation $\boldsymbol{\Theta}$ and Recurrent Bins)} The desired
formation shape $\boldsymbol{\Theta}$ is a probability vector in
$\mathbb{R}^{n_{\textrm{bin}}}$. Each element $\boldsymbol{\Theta}[i]$
represents the desired swarm density distribution in the corresponding
bin $B[i]$. The bins that have nonzero elements in $\boldsymbol{\Theta}$
are called recurrent bins. Let $n_{\mathrm{rec}}$ denote the number
of recurrent bins. The remaining bins, with zero elements in $\boldsymbol{\Theta}$,
are called transient bins. Without loss of generality, we re-label
the bins such that the first $n_{\mathrm{rec}}$ bins are recurrent
bins (i.e., $\boldsymbol{\Theta}[i]>0$ for all $i\in\{1,\ldots,n_{\textrm{rec}}\}$)
and the remaining bins are transient bins (i.e., $\boldsymbol{\Theta}[i]=0$
for all $i\in\{n_{\textrm{rec}}+1,\ldots,n_{\textrm{bin}}\}$). For
example, see Fig.~\ref{fig:I-example}. \hfill$\Box$ \end{definition}

Note that representing the desired formation as a distribution over
bins is analogous to representing a 2-D image using pixels or a 3-D
shape using voxels. In this paper, the complex desired formation shapes
of the Taj Mahal in Fig.~\ref{fig:swarm-Taj-Mahal} and the Eiffel
Tower in Fig.~\ref{fig:Multi-resolution-Images-of-Eiffel-Tower}
are generated in this manner. 

\begin{assumption} \label{assump:m_k} Let the scalar $m_{k}\in\mathbb{N}$
denote the number of agents in the swarm at the $k^{\textrm{th}}$
time instant. We assume that $m_{k}\gg n_{\mathrm{rec}}$ and the
agents do not keep track of $m_{k}$. In Section~\ref{sec:Convergence-Analysis},
we provide a lower-bound on $m_{k}$ for achieving the desired convergence
error. Moreover, we can only achieve the best quantized representation
of $\boldsymbol{\Theta}$ using $m_{k}$ agents, due to the quantization
error of $\frac{1}{m_{k}}$. For example, if $\boldsymbol{\Theta}=[\frac{1}{3},\thinspace\frac{2}{3}]$
and $m_{k}=10$, then the best-quantized representation of $\boldsymbol{\Theta}$
that can be achieved is $[0.3,\thinspace0.7]$. \hfill$\Box$ \end{assumption}

\begin{assumption} \label{assump:identical_agents} We assume that
the agents are \textit{anonymous} and \textit{identical}, i.e., the
agents do not have any global identifiers and all agents execute the
same algorithm \cite{Ref:Olshevsky10}. If the agents are indexed
(non-anonymous), then a spanning-tree-based algorithm can be executed
\cite{Ref:McLurkin14}, but this is not possible in our case. \hfill$\Box$
\end{assumption}

\begin{assumption} \label{assump:r_known} We assume that each agent
can sense the bin it belongs to. Note that this requirement is less
stringent than sensing the precise location. For example, outdoor
robots and spacecraft in low Earth orbit can sense their position
in a distributed manner using the Global Positioning System (GPS)
or other relative navigation technologies (e.g., cellphone towers,
radio beacons). 

The indicator (row) vector $\boldsymbol{r}_{k}^{j}\in\mathbb{R}^{n_{\textrm{bin}}}$
represents the actual bin position of the $j^{\textrm{th}}$ agent
at the $k^{\textrm{th}}$ time instant. If the element $\boldsymbol{r}_{k}^{j}[i]=1$,
then the $j^{\textrm{th}}$ agent is present inside the bin $B[i]$
at the $k^{\textrm{th}}$ time instant; otherwise $\boldsymbol{r}_{k}^{j}[i]=0$.
\hfill$\Box$ \end{assumption}

\begin{definition} \label{def:swarm-distribution} \textit{(Current
Swarm Distribution $\boldsymbol{\mu}_{k}^{\star}$)} The current swarm
distribution $\boldsymbol{\mathcal{\mu}}_{k}^{\star}$ is a probability
vector in $\mathbb{R}^{n_{\textrm{bin}}}$. It is given by the ensemble
mean of actual bin positions of the agents:
\begin{equation}
\boldsymbol{\mathcal{\mu}}_{k}^{\star}:=\frac{1}{m_{k}}\sum_{j=1}^{m_{k}}\boldsymbol{r}_{k}^{j}\thinspace.\label{eq:current_swarm_distribution}
\end{equation}
Each element $\boldsymbol{\mathcal{\mu}}_{k}^{\star}[i]$ gives the
swarm density distribution in the corresponding bin $B[i]$ at the
$k^{\textrm{th}}$ time instant. \hfill$\Box$ \end{definition}

\begin{figure}[h]
\begin{centering}
\includegraphics[width=2.5in]{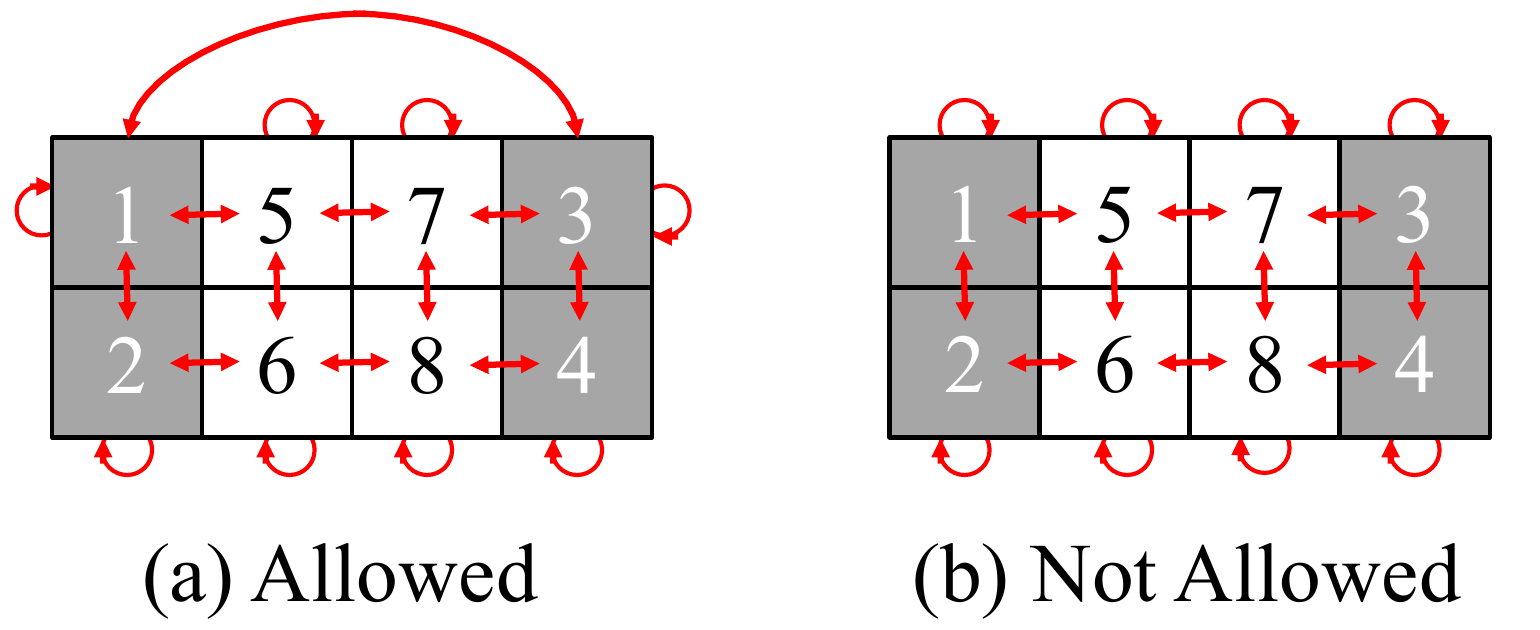}
\par\end{centering}
\caption{In this example, the bins 1 to 4 are recurrent bins. The allowed transitions
(motion constraints) are shown as red arrows. Case (a) satisfies all
the properties in Definition~\ref{def:Motion-Constraints}. Case
(b) does not satisfy property (3) because the recurrent bins are not
strongly connected. \label{fig:motion-constraint}}
\end{figure}

\begin{definition} \label{def:Motion-Constraints} \textit{(Matrix
$\boldsymbol{A}_{k}^{j}$ of Motion Constraints)} An agent in a particular
bin can only transition to some bins but cannot transition to other
bins, because of the dynamics or physical constraints. These (possibly
time-varying) motion constraints are specified by the matrix $\boldsymbol{A}_{k}^{j}\in\mathbb{R}^{n_{\textrm{bin}}\times n_{\textrm{bin}}}$,
where each element is given by:
\begin{align}
\boldsymbol{A}_{k}^{j}[i,\ell] & =\begin{cases}
1 & \textrm{if the transition from bin \ensuremath{B[i]}}\textrm{ to bin }\\
 & \ensuremath{B[\ell]}\textrm{ is allowed at the }k^{\textrm{th}}\textrm{ time instant},\\
0 & \textrm{if this transition is not allowed}.
\end{cases}\label{eq:motion_constraints}
\end{align}
We assume that $\boldsymbol{A}_{k}^{j}$ satisfies the following properties: 
\begin{enumerate}
\item The matrix $\boldsymbol{A}_{k}^{j}$ is symmetric and irreducible\footnote{A matrix is irreducible if and only if the graph conforming to that
matrix is strongly connected.},
\item $\boldsymbol{A}_{k}^{j}[i,i]=1$ for all agents, bins, and time instants, 
\item Since the first $n_{\mathrm{rec}}$ bins are recurrent bins, the sub-matrix
$\boldsymbol{A}_{k,\textrm{sub}}^{j}:=\boldsymbol{A}_{k}^{j}[1\!:\!n_{\textrm{rec}},1\!:\!n_{\textrm{rec}}]$
encapsulate the motion constraints between the recurrent bins. The
matrix $\boldsymbol{A}_{k,\textrm{sub}}^{j}$ is irreducible.
\end{enumerate}
These properties are visualized in Fig.~\ref{fig:motion-constraint}.
\hfill$\Box$ \end{definition}

As shown in Fig.~\ref{fig:motion-constraint}, the recurrent bins
need not be contiguous. Therefore, the desired distribution can have
multiple disconnected components. Note that the matrix $\boldsymbol{A}_{k}^{j}$
is different from the Markov matrix introduced in Section~\ref{sec:Construction-of-Markov-Matrix}.

\begin{definition} \label{def:cost-function} \textit{(Cost Matrix
$\boldsymbol{C}_{k}$)} Consider a matrix $\boldsymbol{C}_{k}\in\mathbb{R}^{n_{\textrm{bin}}\times n_{\textrm{bin}}}$
that encapsulates the cost of transitions between bins. Each element
$\boldsymbol{C}_{k}[i,\ell]$ denotes the cost incurred by an agent
while transitioning from bin $B[i]$ to bin $B[\ell]$ at the $k^{\textrm{th}}$
time instant. This cost represents the control effort or the fuel
consumed or any other metric that the agents seek to minimize. We
assume that the agents do not incur any cost if they remain in their
present bin and the agents incur some positive cost if they transition
out of their present bin (i.e., $\boldsymbol{C}_{k}[i,i]=0$ and $\boldsymbol{C}_{k}[i,\ell]>0$
for all $i,\ell\in\{1,\ldots,n_{\textrm{bin}}\}$ and $i\not=\ell$).
\hfill$\Box$ \end{definition}

\begin{assumption} \label{assump:prior-knowledge} \textit{(A Priori
Information Available with Each Agent)} We assume that the physical
space over which the swarm is distributed, the position of the bins
$B[i]$, and the desired formation shape $\boldsymbol{\Theta}$ are
known before the algorithm starts. Subsequently, the time-varying
cost matrices $\boldsymbol{C}_{k}$ and the motion constraints matrices
$\boldsymbol{A}_{k}^{j}$, which depend on the position of the bins,
are computed and stored beforehand. Finally, the four design variables
(namely $\varepsilon_{M}$, $\varepsilon_{C}$, $\beta^{j}$, and
$\tau^{j}$ ), which are introduced later, are also stored on board
each agent.\hfill$\Box$ \end{assumption}

\subsection{Distributed Estimation of the Current Swarm Distribution \label{subsec:Generating-Feedback-swarm-dist}}

The algorithms in this paper use feedback of the current swarm distribution
$\boldsymbol{\mathcal{\mu}}_{k}^{\star}$. In order to generate this
estimate in a distributed manner, we need the following assumption. 

\begin{assumption} \label{assump:Strongly-connected-communication}
The time-varying communication network topology of the swarm is assumed
to be strongly-connected. \hfill$\Box$ \end{assumption}

Under Assumption \ref{assump:Strongly-connected-communication}, several
algorithms exist in the literature for estimating $\boldsymbol{\mathcal{\mu}}_{k}^{\star}$
in a distributed manner \cite{Ref:Chen02,Ref:Saber09,Ref:Bandyopadhyay14_ACC_BCF}.
For example, the distributed consensus algorithm \cite{Ref:Bandyopadhyay14_ACC_BCF}
is used to estimate $\boldsymbol{\mathcal{\mu}}_{k}^{\star}$ in Remark
\ref{rem:consensus-distributed-estimation} in Appendix.

Any distributed estimation algorithm will have some residual estimation
error between the current swarm distribution $\boldsymbol{\mathcal{\mu}}_{k}^{\star}$
and the $j^{\textrm{th}}$ agent's estimate of the current swarm distribution
at the $k^{\textrm{th}}$ time instant, which is represented by the
probability vector $\boldsymbol{\mathcal{\mu}}_{k}^{j}\in\mathbb{R}^{n_{\textrm{bin}}}$.
Let the positive parameter $\epsilon_{\mathrm{est}}$ represent the
maximum estimation error between $\boldsymbol{\mathcal{\mu}}_{k}^{\star}$
and $\boldsymbol{\mathcal{\mu}}_{k}^{j}$:
\begin{equation}
D_{\mathcal{L}_{1}}(\boldsymbol{\mathcal{\mu}}_{k}^{\star},\boldsymbol{\mathcal{\mu}}_{k}^{j})=\sum_{i=1}^{n_{\mathrm{bin}}}\left|\boldsymbol{\mathcal{\mu}}_{k}^{\star}[i]-\boldsymbol{\mathcal{\mu}}_{k}^{j}[i]\right|\leq\epsilon_{\mathrm{est}}\thinspace,\thinspace\forall k\in\mathbb{N}\thinspace.\label{eq:estimation_error}
\end{equation}
We later show that our algorithm works remarkably well in the presence
of this estimation error $\epsilon_{\mathrm{est}}$ (\ref{eq:estimation_error}). 

\subsection{Problem Statement for Shape Formation \label{subsec:Problem-Statement-pattern-formation}}

Under Assumptions~\ref{assump:m_k}\textendash \ref{assump:Strongly-connected-communication},
the objectives of the PSG\textendash IMC algorithm for shape formation
are as follows:
\begin{itemize}
\item Each agent independently determines its bin-to-bin trajectory using
a Markov chain, which obeys motion constraints $\boldsymbol{A}_{k}^{j}$,
so that the overall swarm converges to a desired formation shape $\boldsymbol{\Theta}$.
\item The algorithm automatically detects and repairs damages to the formation. 
\item The algorithm minimizes the expected cost of transitions at every
time instant (see Definition~\ref{def:expected-cost-transitions})
for all the agents, where the cost matrix $\boldsymbol{C}_{k}$ is
defined in Definition~\ref{def:cost-function}.
\end{itemize}
\begin{figure}[h]
\begin{centering}
\includegraphics[bb=270bp 450bp 1060bp 1580bp,clip,width=3in]{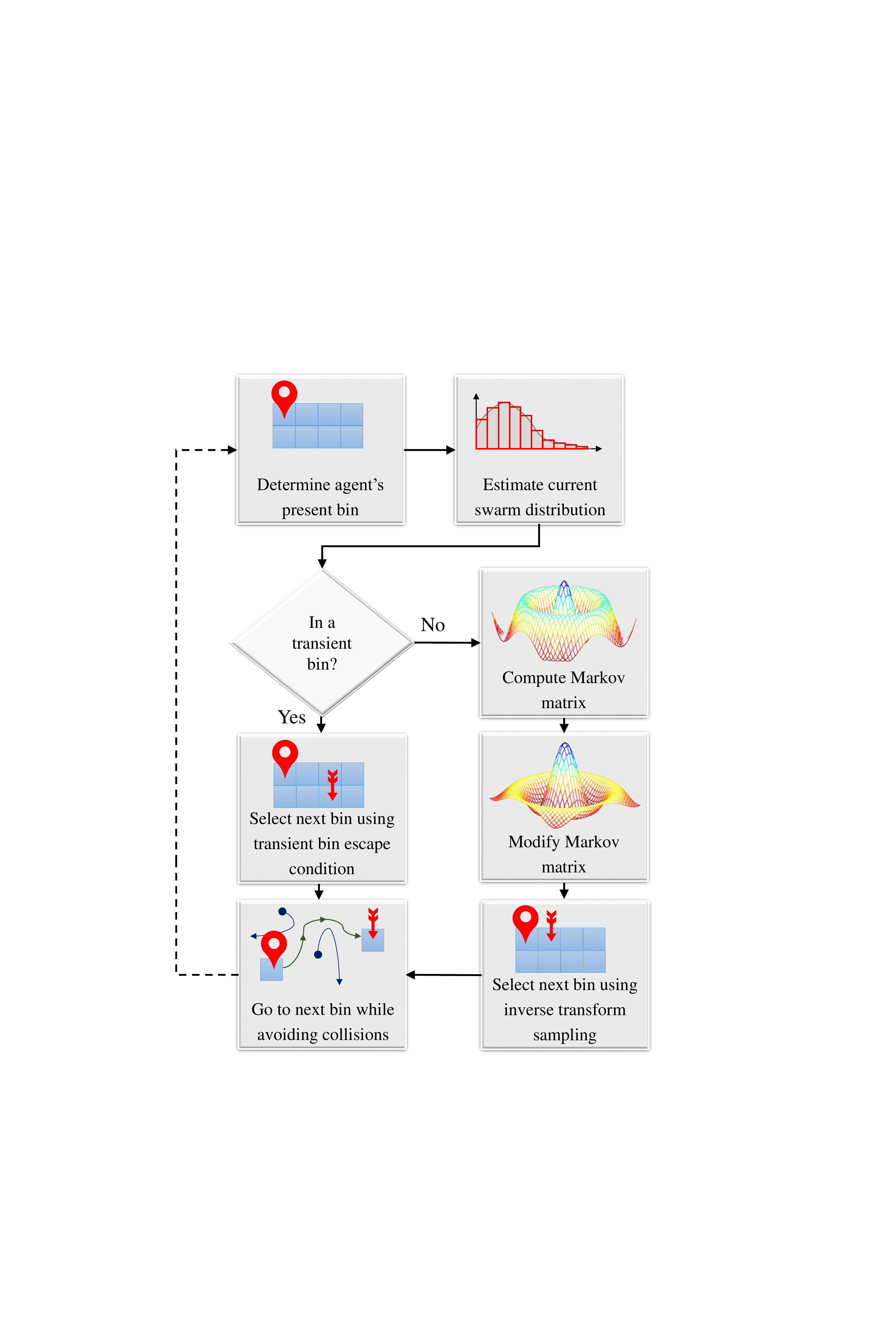}
\par\end{centering}
\caption{Flowchart of the PSG\textendash IMC algorithm for shape formation.
\label{fig:Flowchart-PSGIMC}}
\end{figure}

\subsection{Outline of the PSG\textendash IMC Algorithm \label{subsec:PSG-IMC-flowchart}}

The key steps in the proposed PSG\textendash IMC algorithm for shape
formation are shown in Fig.~\ref{fig:Flowchart-PSGIMC}. The agent
first determines its present bin and estimates the current swarm distribution
(Section~\ref{subsec:Generating-Feedback-swarm-dist}). If the agent
is in a transient bin, then it selects another bin using the condition
for escaping transient bins (Section~\ref{subsec:Condition-for-Escaping-transient-bins}).
Otherwise, the agent computes the Markov matrix (Section~\ref{subsec:Construction-of-Minimum-cost-Markov-matrix})
and then modifies it to suppress undesirable transitions (Section~\ref{sec:PSG-IMC-pattern-formation}).
Finally, the agent uses inverse transform sampling to select the next
bin (Remark~\ref{rem:Inverse-Transform-Sampling} in Appendix). The
agent uses a lower-level guidance and control algorithm, which depends
on the agent's dynamics, to go from its present bin to the selected
bin in a collision-free manner. Such lower-level algorithms based
on real-time optimal control and Voronoi partitions are presented
in \cite{Ref:Bandyopadhyay14MSC,Ref:Giri14,Ref:Morgan15_SATO} and
also discussed briefly in Remark~\ref{rem:collision-free-motion}
in Appendix. The pseudo-code for the PSG\textendash IMC algorithm
for shape formation is given in \textbf{Method~2} in Section \ref{sec:PSG-IMC-pattern-formation}.

\section{Construction of Feedback-based Markov Matrix \label{sec:Construction-of-Markov-Matrix}}

In this section, we present our novel techniques for constructing
Markov matrices.

\subsection{Construction of Minimum Cost Markov Matrix \label{subsec:Construction-of-Minimum-cost-Markov-matrix}}

In this subsection, we construct Markov matrices that minimize the
expected cost of transitions at each time instant. 

\begin{definition} \label{def:Feedback-gain} \textit{(Feedback Gain
$\xi_{k}^{j}$ and Desired Convergence Error $\xi_{\textrm{des}}$)}
The feedback gain $\xi_{k}^{j}$ is given by the Hellinger distance
(HD) between the current swarm distribution $\boldsymbol{\mathcal{\mu}}_{k}^{j}$
and the desired formation $\boldsymbol{\Theta}$:
\begin{align}
\xi_{k}^{j} & =D_{H}(\boldsymbol{\Theta},\boldsymbol{\mathcal{\mu}}_{k}^{j}):=\frac{1}{\sqrt{2}}\sqrt{\sum_{i=1}^{n_{\textrm{bin}}}\left(\sqrt{\boldsymbol{\Theta}[i]}-\sqrt{\boldsymbol{\mathcal{\mu}}_{k}^{j}[i]}\right)^{2}}.\label{eq:tuning_eqn}
\end{align}
The HD is a symmetric measure of the difference between two probability
distributions and bounded by $1$ \cite{Ref:Torgerson91,Ref:Cha07}.

Let $\xi_{\textrm{des}}$ represent the desired convergence error
threshold between the final swarm distribution and $\boldsymbol{\Theta}$.\hfill$\Box$
\end{definition}

\begin{figure}[h]
\begin{centering}
\includegraphics[bb=0bp 270bp 900bp 540bp,clip,width=3.4in]{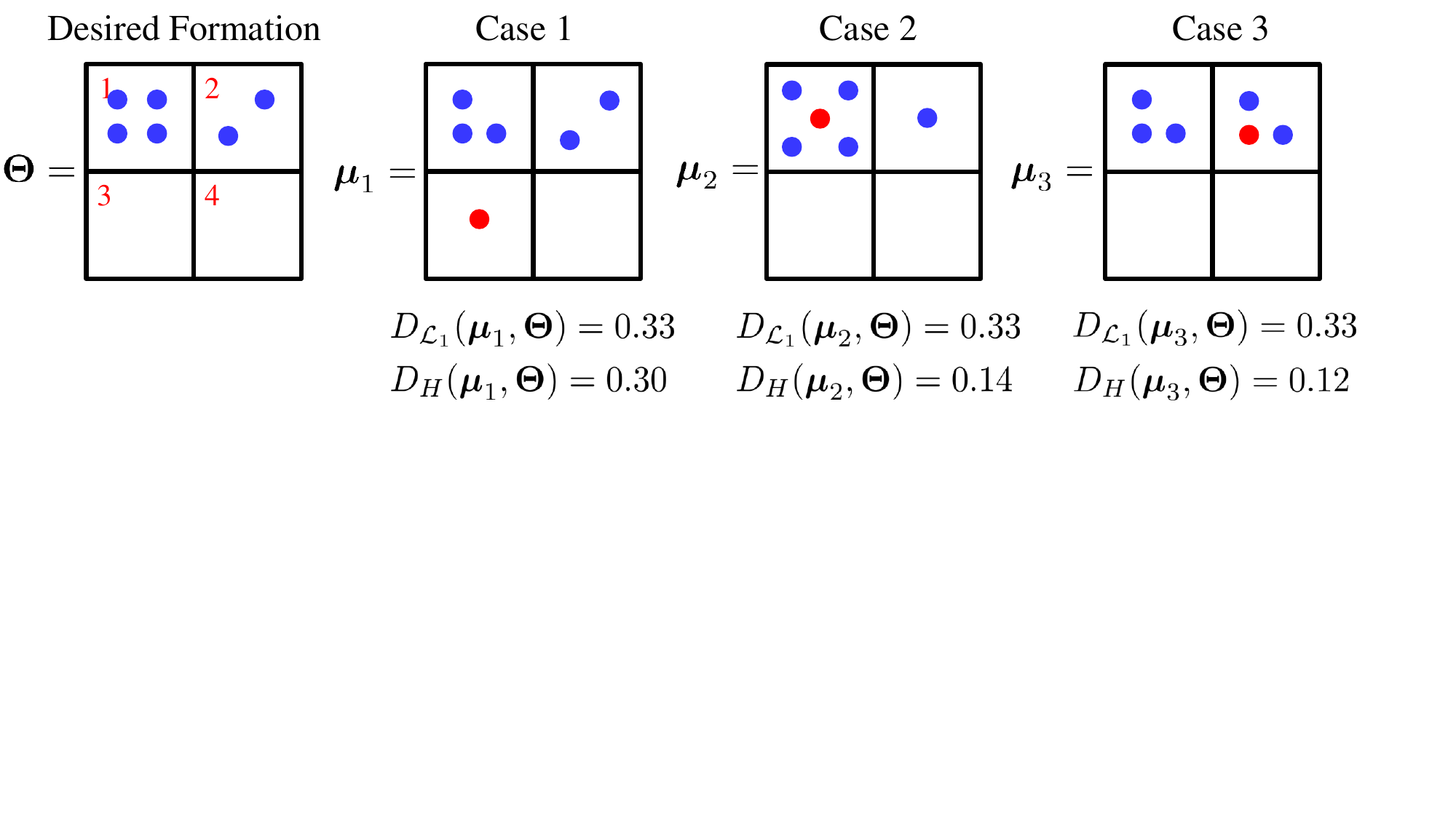}
\par\end{centering}
\caption{In this example, the desired distribution $\boldsymbol{\Theta}$ has
$4$ and $2$ agents in bins 1 and 2 respectively. In the three cases,
one agent (marked in red) is not in its correct bin. The $\mathcal{L}_{1}$
distances are equal, but the HD are different. \label{fig:HD_comparision}}
\end{figure}

\begin{remark} \label{rem:Hellinger-Distance} \textit{(Advantages
of Hellinger Distance)} The HD between $\boldsymbol{\mathcal{\mu}}_{k}^{j}$
and $\boldsymbol{\Theta}$ in (\ref{eq:tuning_eqn}) is bounded as
follows \cite{Ref:Pollard00asymptopia}:
\begin{equation}
\frac{1}{2\sqrt{2}}D_{\mathcal{L}_{1}}(\boldsymbol{\Theta},\boldsymbol{\mathcal{\mu}}_{k}^{j})\leq D_{H}(\boldsymbol{\Theta},\boldsymbol{\mathcal{\mu}}_{k}^{j})\leq\frac{1}{\sqrt{2}}D_{\mathcal{L}_{1}}(\boldsymbol{\Theta},\boldsymbol{\mathcal{\mu}}_{k}^{j})^{\frac{1}{2}}\thinspace.\label{eq:HD_L1_compare}
\end{equation}
We choose HD, over other popular metrics like $\mathcal{L}_{1}$ and
$\mathcal{L}_{2}$ distances, because of its properties illustrated
in Fig.~\ref{fig:HD_comparision}. The $\mathcal{L}_{1}$ distances
for the cases ($\boldsymbol{\mathcal{\mu}}_{1}$, $\boldsymbol{\mu}_{2}$,
$\boldsymbol{\mathcal{\mu}}_{3}$) from $\boldsymbol{\Theta}$ are
equal. But in Case 1, the wrong agent is in a bin where there should
be no agent, hence HD heavily penalizes this case. If all the agents
are only in those bins which have positive weights in $\boldsymbol{\Theta}$,
then HD is significantly smaller. Finally, if an agent is missing
from a bin that has fewer agents in $\boldsymbol{\Theta}$ (Case 2)
compared to a bin that has more agents in $\boldsymbol{\Theta}$ (Case
3), then HD penalizes Case 2 slightly more than Case 3. These properties
are useful for swarm guidance. \hfill$\Box$ \end{remark}

Consider the Markov matrix $\boldsymbol{M}_{k}^{j}$ in $\mathbb{R}^{n_{\textrm{bin}}\times n_{\textrm{bin}}}$
that encapsulates the transition probabilities between bins. Each
element $\boldsymbol{M}_{k}^{j}[i,\ell]$ represents the probability
that the $j^{\textrm{th}}$ agent in bin $B[i]$ at the $k^{\textrm{th}}$
time instant will transition to bin $B[\ell]$ at the $(k+1)^{\textrm{th}}$
time instant:
\begin{equation}
\boldsymbol{M}_{k}^{j}[i,\ell]:=\mathbb{P}\left(\boldsymbol{r}_{k+1}^{j}[\ell]=1|\boldsymbol{r}_{k}^{j}[i]=1\right)\thinspace.\label{eq:Markov_matrix}
\end{equation}
Therefore, the Markov matrix $\boldsymbol{M}_{k}^{j}$ is row stochastic
(i.e., $\boldsymbol{M}_{k}^{j}\boldsymbol{1}=\boldsymbol{1}$). Its
stationary distribution is defined as follows. 

\begin{definition} \label{def:stat_dist} \textit{(Stationary Distribution)}
The stationary distribution $\boldsymbol{e}_{k}^{j}$ of the Markov
matrix $\boldsymbol{M}_{k}^{j}$ is given by the solution of $\boldsymbol{e}_{k}^{j}\boldsymbol{M}_{k}^{j}=\boldsymbol{e}_{k}^{j}$,
where $\boldsymbol{e}_{k}^{j}$ is a probability (row) vector in $\mathbb{R}^{n_{\textrm{bin}}}$
(i.e., $\boldsymbol{e}_{k}^{j}\geq0$, $\boldsymbol{e}_{k}^{j}\boldsymbol{1}=1$).
The stationary distribution is unique if the Markov matrix is irreducible
\cite[pp. 119]{Ref:Seneta06}. \hfill$\Box$ \end{definition}

\begin{definition} \label{def:expected-cost-transitions} \textit{(Expected
Cost of Transitions at Each Time Instant)} The expected cost of transitions
for the $j^{\textrm{th}}$ agent at the $k^{\textrm{th}}$ time instant
is given by $\sum_{i=1}^{n_{\mathrm{bin}}}\sum_{\ell=1}^{n_{\mathrm{bin}}}\boldsymbol{C}_{k}[i,\ell]\boldsymbol{M}_{k}^{j}[i,\ell]$,
where the cost matrix $\boldsymbol{C}_{k}$ is defined in Definition
\ref{def:cost-function}. \hfill$\Box$ \end{definition}

The following \textbf{Method~1,} Theorem~\ref{thm:Markov-optimization},
and Corollary~\ref{cor:Markov-optimization-scalable} present our
construction of the optimal Markov matrix $\boldsymbol{M}_{k}^{j}$
that minimizes this expected cost of transitions at the each time
instant. Our construction technique has no relation with the well-known
Metropolis-Hastings (MH) algorithm, which is commonly used for constructing
Markov matrices with a given stationary distribution \cite{Ref:Chib95,Ref:Diaconis01}.
In the MH algorithm, the proposal distribution is used to iteratively
generate the next sample, which is accepted or rejected based on the
desired stationary distribution. There is no direct method for incorporating
feedback into the MH algorithm. In contrast, the feedback of the current
swarm distribution is directly incorporated within our construction
process using the feedback gain.

\vspace{5pt}

\noindent\fbox{\begin{minipage}[t]{1\columnwidth - 2\fboxsep - 2\fboxrule}%
\textbf{Method~1: Construction of Optimal Markov Matrix}

Under Assumptions~\ref{assump:m_k}\textendash \ref{assump:Strongly-connected-communication},
the optimal Markov matrix $\boldsymbol{M}_{k}^{j}$ that minimizes
the expected cost of transitions at each time instant is constructed
as follows:

(A) If $\xi_{k}^{j}<\xi_{\mathrm{des}}$, then set $\boldsymbol{M}_{k}^{j}=\mathbf{I}$.

(B) Otherwise, $\boldsymbol{M}_{k}^{j}$ is computed as follows:\\
$\mathrm{(CS1)}$ If $\boldsymbol{A}_{k}^{j}[i,\ell]=0$, then set
$\boldsymbol{M}_{k}^{j}[i,\ell]=0$ for all bins $i,\ell\in\{1,\ldots,n_{\mathrm{bin}}\}$.\\
$\mathrm{(CS2)}$ If $\boldsymbol{\Theta}[\ell]=0$, then set $\boldsymbol{M}_{k}^{j}[i,\ell]=0$
for all bins $i,\ell\in\{1,\ldots,n_{\mathrm{bin}}\}$ with $i\not=\ell$.\\
The remaining elements in $\boldsymbol{M}_{k}^{j}$ are computed using
the following linear program (LP):
\begin{align}
 & \mathrm{minimize}\thinspace\sum_{i=1}^{n_{\mathrm{bin}}}\sum_{\ell=1}^{n_{\mathrm{bin}}}\boldsymbol{C}_{k}[i,\ell]\boldsymbol{M}_{k}^{j}[i,\ell]\thinspace,\label{eq:Markov_optimization_problem}\\
\textrm{subject to } & \sum_{\ell=1}^{n_{\mathrm{bin}}}\boldsymbol{M}_{k}^{j}[i,\ell]=1,\thinspace\forall i\thinspace, & \mathrm{(LP1)}\nonumber \\
 & \sum_{i=1}^{n_{\mathrm{bin}}}\boldsymbol{\Theta}[i]\boldsymbol{M}_{k}^{j}[i,\ell]=\boldsymbol{\Theta}[\ell],\thinspace\forall\ell\thinspace, & \mathrm{(LP2)}\nonumber \\
 & (1-\xi_{k}^{j})\leq\boldsymbol{M}_{k}^{j}[i,i]\leq1,\thinspace\forall i\thinspace, & \mathrm{(LP3)}\nonumber \\
 & \varepsilon_{M}\xi_{k}^{j}\boldsymbol{\Theta}[\ell]\left(1-\frac{\boldsymbol{C}_{k}[i,\ell]}{C_{k,\max}+\varepsilon_{C}}\right)\nonumber \\
 & \qquad\qquad\leq\boldsymbol{M}_{k}^{j}[i,\ell]\leq\frac{\xi_{k}^{j}}{\varepsilon_{M}},\thinspace\forall i\not=\ell\thinspace, & \mathrm{(LP4)}\nonumber 
\end{align}
where $\varepsilon_{M}$ is a positive scalar constant in $(0,1]$,
$C_{k,\max}$ is the maximum transition cost (i.e., $C_{k,\max}=\max_{i,\ell}\boldsymbol{C}_{k}[i,\ell]$),
and $\varepsilon_{C}$ is a positive scalar constant. %
\end{minipage}}

\begin{remark} \label{rem:explain-LP-constraints} The Markov matrix
$\boldsymbol{M}_{k}^{j}$ designed in \textbf{Method~1} has the following
desirable properties:
\begin{itemize}
\item If $\xi_{k}^{j}<\xi_{\mathrm{des}}$ (see Definition~\ref{def:Feedback-gain}),
the swarm is deemed to have converged to the desired formation. Then,
$\boldsymbol{M}_{k}^{j}$ is set to an identity matrix so that the
agents do not transition anymore. Hence, Step (A) ensures that the
swarm remains converged without additional movement. 
\item If the swarm has not converged to the desired formation (i.e., $\xi_{k}^{j}\geq\xi_{\mathrm{des}}$),
then Step (B) is initiated. 
\item (CS1) prevents those transitions that are not allowed by motion constraints.
\item (CS2) prevents transitions into transient bins. 
\item The objective function in LP (\ref{eq:Markov_optimization_problem})
minimizes the expected cost of transitions at the current time instant
(see Definition~\ref{def:expected-cost-transitions}). 
\item (LP1) ensures that $\boldsymbol{M}_{k}^{j}$ is row stochastic.
\item (LP2) ensures that $\boldsymbol{M}_{k}^{j}$ has $\boldsymbol{\Theta}$
as its stationary distribution (i.e., $\boldsymbol{\Theta}\boldsymbol{M}_{k}^{j}=\boldsymbol{\Theta}$).
\item The lower bound in (LP3) ensures that there is non-zero probability
for agents to remain in their present bin when $\xi_{k}^{j}<1$. The
upper bound in (LP3) is derived from (LP1).
\item The lower bound in (LP4) ensures that the minimum transition probability
to a target bin is non-zero and directly proportional to both the
feedback gain $\xi_{k}^{j}$ and the target bin's desired distribution
$\boldsymbol{\Theta}[\ell]$. But the minimum transition probability
decreases with increasing cost of transition to the target bin. 
\item The upper bound in (LP4) ensures that the maximum transition probability
is also directly proportional to the feedback gain $\xi_{k}^{j}$. 
\end{itemize}
A salient feature of the constraints (LP3,4) is that they depend on
the feedback gain $\xi_{k}^{j}$. Therefore, if the swarm distribution
$\boldsymbol{\mathcal{\mu}}_{k}^{j}$ converges to $\boldsymbol{\Theta}$
(i.e., $\boldsymbol{\mathcal{\mu}}_{k}^{j}\rightarrow\boldsymbol{\Theta}$),
then $\xi_{k}^{j}\rightarrow0$ (because $\xi_{k}^{j}=D_{H}(\boldsymbol{\Theta},\boldsymbol{\mathcal{\mu}}_{k}^{j})$)
and $\boldsymbol{M}_{k}^{j}\rightarrow\mathbf{I}$ based on these
constraints. The identity matrix ensures that agents settle down after
the desired formation is achieved, thereby reducing unnecessary transitions.
In Section~\ref{sec:Convergence-Analysis}, we show that these constraints
also help prove the convergence of the algorithm. \hfill$\Box$ \end{remark}

\begin{theorem} \label{thm:Markov-optimization} The feasible set
of Markov matrices that satisfy the constraints (CS1,2) and the linear
constraints in LP (\ref{eq:Markov_optimization_problem}) in \textbf{Method~1}
is non-empty. The optimal Markov matrix $\boldsymbol{M}_{k}^{j}$
is row-stochastic, has $\boldsymbol{\Theta}$ as its stationary distribution,
and only allows transitions into recurrent bins. \end{theorem}

\begin{proof} The optimization problem in (\ref{eq:Markov_optimization_problem})
is an LP because the constraints are all linear inequalities or equalities
and the objective function is linear. An optimal solution for the
LP exists if the feasible set of Markov matrices is non-empty. We
now show that the following family of Markov matrices $\boldsymbol{Q}_{k}^{j}$
are within the set of feasible solutions:
\begin{align}
\boldsymbol{Q}_{k}^{j}[i,\ell] & =\begin{cases}
0 & \textrm{ if }\boldsymbol{A}_{k}^{j}[i,\ell]=0\\
\frac{\xi_{k}^{j}}{\boldsymbol{\Theta}\boldsymbol{\alpha}_{k}^{j}}\left(\boldsymbol{\alpha}_{k}^{j}[i]\boldsymbol{\alpha}_{k}^{j}[\ell]\boldsymbol{\Theta}[\ell]\right) & \textrm{ otherwise }
\end{cases}\thinspace,\nonumber \\
 & \qquad\forall i,\ell\in\{1,\ldots,n_{\textrm{bin}}\}\textrm{ and }i\not=\ell\thinspace,\label{eq:Q_off_diag}\\
\boldsymbol{Q}_{k}^{j}[i,i] & =\frac{\xi_{k}^{j}}{\boldsymbol{\Theta}\boldsymbol{\alpha}_{k}^{j}}\left(\boldsymbol{\alpha}_{k}^{j}[i]\boldsymbol{\alpha}_{k}^{j}[i]\boldsymbol{\Theta}[\ell]\right)+\left(1-\xi_{k}^{j}\boldsymbol{\alpha}_{k}^{j}[i]\right)\nonumber \\
 & \qquad+\!\!\!\!\!\!\!\sum_{\ell\in\left\{ \boldsymbol{A}_{k}^{j}[i,\ell]=0\right\} }\frac{\xi_{k}^{j}}{\boldsymbol{\Theta}\boldsymbol{\alpha}_{k}^{j}}\left(\boldsymbol{\alpha}_{k}^{j}[i]\boldsymbol{\alpha}_{k}^{j}[\ell]\boldsymbol{\Theta}[\ell]\right)\thinspace.\label{eq:Q_diag}
\end{align}
where $\varepsilon_{\alpha}=\sqrt{\varepsilon_{M}}$ and $\boldsymbol{\alpha}_{k}^{j}$
is a positive column vector in $\mathbb{R}^{n_{\mathrm{bin}}}$, with
$\varepsilon_{\alpha}\leq\boldsymbol{\alpha}_{k}^{j}[i]\leq1$ for
all bins. 

The matrix $\boldsymbol{Q}_{k}^{j}$ satisfies (CS1) due to (\ref{eq:Q_off_diag}).
If $\boldsymbol{\Theta}[\ell]=0$, then the off-diagonal element satisfies
$\boldsymbol{Q}_{k}^{j}[i,\ell]=0$ and the matrix $\boldsymbol{Q}_{k}^{j}$
satisfies (CS2). 

We now show that the matrix $\boldsymbol{Q}_{k}^{j}$ satisfies (LP1):
\[
\sum_{\ell=1}^{n_{\mathrm{bin}}}\boldsymbol{Q}_{k}^{j}[i,\ell]=\tfrac{\xi_{k}^{j}\boldsymbol{\alpha}_{k}^{j}[i]}{\boldsymbol{\Theta}\boldsymbol{\alpha}_{k}^{j}}\sum_{\ell=1}^{n_{\mathrm{bin}}}\boldsymbol{\alpha}_{k}^{j}[\ell]\boldsymbol{\Theta}[\ell]+1-\xi_{k}^{j}\boldsymbol{\alpha}_{k}^{j}[i]=1\thinspace,
\]
where $\sum_{\ell=1}^{n_{\mathrm{bin}}}\left(\boldsymbol{\alpha}_{k}^{j}[\ell]\boldsymbol{\Theta}[\ell]\right)=\boldsymbol{\Theta}\boldsymbol{\alpha}_{k}^{j}$.
We now prove that the matrix $\boldsymbol{Q}_{k}^{j}$ satisfies (LP2):
\begin{align*}
 & \sum_{i=1}^{n_{\mathrm{bin}}}\boldsymbol{\Theta}[i]\boldsymbol{Q}_{k}^{j}[i,\ell]=\frac{\xi_{k}^{j}\boldsymbol{\alpha}_{k}^{j}[\ell]\boldsymbol{\Theta}[\ell]}{\boldsymbol{\Theta}\boldsymbol{\alpha}_{k}^{j}}\sum_{i=1}^{n_{\mathrm{bin}}}\left(\boldsymbol{\alpha}_{k}^{j}[i]\boldsymbol{\Theta}[i]\right)\\
 & \qquad\qquad\qquad\qquad+\left(\boldsymbol{\Theta}[\ell]-\xi_{k}^{j}\boldsymbol{\Theta}[\ell]\boldsymbol{\alpha}_{k}^{j}[\ell]\right)=\boldsymbol{\Theta}[\ell]\thinspace.
\end{align*}

The matrix $\boldsymbol{Q}_{k}^{j}$ satisfies (LP3) because each
diagonal element $\boldsymbol{Q}_{k}^{j}[i,i]$ is lower bounded by
$\left(1-\xi_{k}^{j}\boldsymbol{\alpha}_{k}^{j}[i]\right)$, which
is one of the positive terms in (\ref{eq:Q_diag}). The term $\left(1-\xi_{k}^{j}\boldsymbol{\alpha}_{k}^{j}[i]\right)$
is lower bounded by $(1-\xi_{k}^{j})$ because $\boldsymbol{\alpha}_{k}^{j}[i]\leq1$. 

We now prove that the matrix $\boldsymbol{Q}_{k}^{j}$ satisfies (LP4).
Since the element $\boldsymbol{\Theta}[\ell]>0$, the off-diagonal
element $\boldsymbol{Q}_{k}^{j}[i,\ell]$ is upper bounded by $\frac{\xi_{k}^{j}}{\varepsilon_{M}}$
because $\boldsymbol{\alpha}_{k}^{j}[i]\leq1$, $\boldsymbol{\Theta}[\ell]\leq1$,
and $\left(\frac{\xi_{k}^{j}}{\boldsymbol{\Theta}\boldsymbol{\alpha}_{k}^{j}}\right)\leq\frac{\xi_{k}^{j}}{\varepsilon_{\alpha}}\leq\frac{\xi_{k}^{j}}{\varepsilon_{M}}$.
On the other hand, the off-diagonal element $\boldsymbol{Q}_{k}^{j}[i,\ell]$
is lower bounded by $\varepsilon_{M}\xi_{k}^{j}\boldsymbol{\Theta}[\ell]$
because $\boldsymbol{\alpha}_{k}^{j}[i]\geq\varepsilon_{\alpha}$
and $\left(\frac{\xi_{k}^{j}}{\boldsymbol{\Theta}\boldsymbol{\alpha}_{k}^{j}}\right)\geq\xi_{k}^{j}$
as $\boldsymbol{\Theta}\boldsymbol{1}=1$. Therefore, $\varepsilon_{M}\xi_{k}^{j}\boldsymbol{\Theta}[\ell]\left(1-\frac{\boldsymbol{C}_{k}[i,\ell]}{C_{k,\max}+\varepsilon_{C}}\right)\leq\boldsymbol{Q}_{k}^{j}[i,\ell]$
because $\left(1-\frac{\boldsymbol{C}_{k}[i,\ell]}{C_{k,\max}+\varepsilon_{C}}\right)<1$.
Thus the matrix $\boldsymbol{Q}_{k}^{j}$ satisfies (LP4). 

Therefore, the feasible set is non-empty and the optimal Markov matrix
$\boldsymbol{M}_{k}^{j}$ posses the desirable properties discussed
in Remark~\ref{rem:explain-LP-constraints}. \hfill$\blacksquare$
\end{proof} 

\begin{figure}[h]
\begin{centering}
\includegraphics[bb=0bp 275bp 890bp 540bp,clip,width=3.4in]{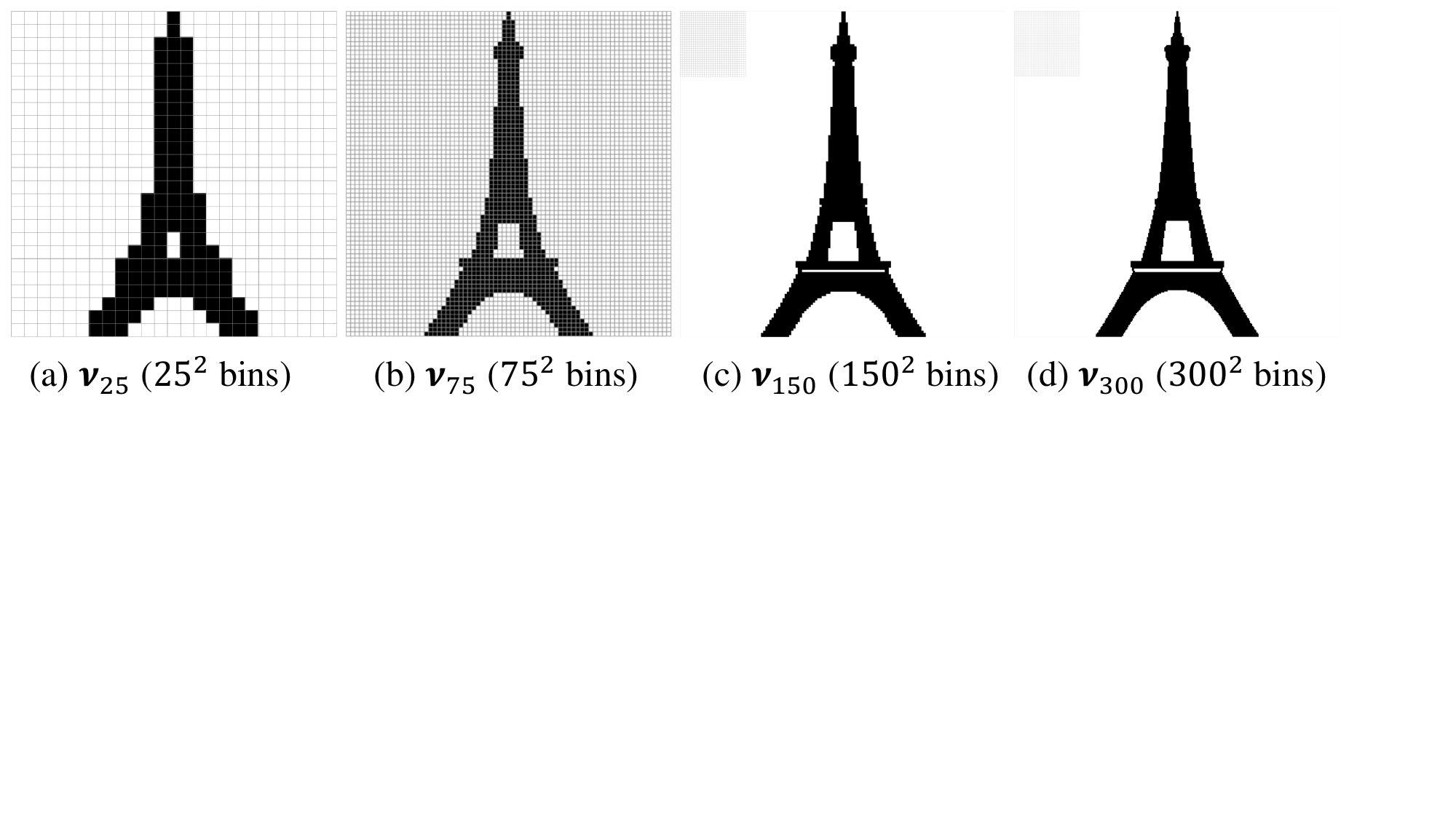}
\par\end{centering}
\caption{Multiresolution images of the Eiffel Tower are shown, where the spatial
resolution increases from (a) to (d). All the bins are shown in (a)
and (b), whereas only a few bins are shown in the left-top corner
in (c) and (d). \label{fig:Multi-resolution-Images-of-Eiffel-Tower}}
\end{figure}

\begin{remark} \label{rem:Computation-time} \textit{(Computation
Time)} Although each agent only needs the row of the Markov matrix
$\boldsymbol{M}_{k}^{j}$ corresponding to its present bin, it has
to solve the entire LP (\ref{eq:Markov_optimization_problem}). The
computation time for an LP increases with an increasing number of
bins because the number of variables in $\boldsymbol{M}_{k}^{j}$
is approximately equal to $n_{\textrm{bin}}^{2}$. For example, if
the desired formation is given by $\boldsymbol{\nu}_{25}$ or $\boldsymbol{\nu}_{75}$
in Fig.~\ref{fig:Multi-resolution-Images-of-Eiffel-Tower}, then
the computation time is a few minutes on a standard desktop computer.
If the desired formation is given by $\boldsymbol{\nu}_{150}$ (with
$5\times10^{8}$ variables) or $\boldsymbol{\nu}_{300}$ (with $8\times10^{9}$
variables), then the LP is impractical for real-time computation.
This escalating computation time with increasing number of bins is
an issue for all LP-based algorithms \cite{Ref:Acikmese15_AR,Ref:Bandyopadhyay14MSC}.
\hfill$\Box$\end{remark}

Therefore, we need a faster method for computing the Markov matrices.
The following corollary gives the analytical formula of the optimal
Markov matrix when the cost matrix is symmetric.

\begin{corollary} \label{cor:Markov-optimization-scalable} The optimal
Markov matrix of the LP (\ref{eq:Markov_optimization_problem}) in
\textbf{Method~1} is given by:
\begin{align}
\boldsymbol{M}_{k}^{j}[i,\ell] & =\begin{cases}
0 & \textrm{ if }\boldsymbol{A}_{k}^{j}[i,\ell]=0\\
\varepsilon_{M}\xi_{k}^{j}\boldsymbol{\Theta}[\ell]\left(1-\frac{\boldsymbol{C}_{k}[i,\ell]}{C_{k,\max}+\varepsilon_{C}}\right) & \textrm{ otherwise }
\end{cases}\thinspace,\nonumber \\
 & \qquad\forall i,\ell\in\{1,\ldots,n_{\mathrm{bin}}\}\textrm{ and }i\not=\ell\thinspace,\label{eq:scalable_Markov1}
\end{align}
\begin{equation}
\boldsymbol{M}_{k}^{j}[i,i]=1-\sum_{\ell\in\{1,\ldots,n_{\mathrm{bin}}\}\backslash\{i\}}\boldsymbol{M}_{k}^{j}[i,\ell]\thinspace.\label{eq:scalable_Markov2}
\end{equation}
if the cost matrix $\boldsymbol{C}_{k}$ is symmetric (i.e., $\boldsymbol{C}_{k}=\boldsymbol{C}_{k}^{T}$).\end{corollary}

\begin{proof} We first state a simpler LP by neglecting the constraints
(LP1,2) from the original LP (\ref{eq:Markov_optimization_problem}).
We state this simpler LP (\ref{eq:Markov_optimization_problem_simpler})
using the following substitutions for all positive elements $\boldsymbol{R}_{k}^{j}[i,i]=\boldsymbol{M}_{k}^{j}[i,i]-(1-\xi_{k}^{j})$
and $\boldsymbol{R}_{k}^{j}[i,\ell]=\boldsymbol{M}_{k}^{j}[i,\ell]-\varepsilon_{M}\xi_{k}^{j}\boldsymbol{\Theta}[\ell]\left(1-\frac{\boldsymbol{C}_{k}[i,\ell]}{C_{k,\max}+\varepsilon_{C}}\right)$:
\begin{align}
 & \mathrm{minimize}\thinspace\sum_{i=1}^{n_{\mathrm{bin}}}\sum_{\ell=1}^{n_{\mathrm{bin}}}\boldsymbol{C}_{k}[i,\ell]\boldsymbol{R}_{k}^{j}[i,\ell]\label{eq:Markov_optimization_problem_simpler}\\
 & +\sum_{i=1}^{n_{\mathrm{bin}}}\sum_{\ell\in\{\boldsymbol{A}_{k}^{j}[i,\ell]=1,i\not=\ell\}}^{n_{\mathrm{bin}}}\!\!\!\!\!\!\!\!\!\!\!\!\boldsymbol{C}_{k}[i,\ell]\varepsilon_{M}\xi_{k}^{j}\boldsymbol{\Theta}[\ell]\left(1-\tfrac{\boldsymbol{C}_{k}[i,\ell]}{C_{k,\max}+\varepsilon_{C}}\right)\thinspace,\nonumber 
\end{align}
subject to:
\begin{align*}
 & 0\leq\boldsymbol{R}_{k}^{j}[i,i]\leq\xi_{k}^{j},\thinspace\forall i\thinspace,\qquad\qquad\qquad\qquad\qquad\qquad\thinspace\mathrm{(\widetilde{LP3})}\\
 & 0\leq\boldsymbol{R}_{k}^{j}[i,\ell]\leq\frac{\xi_{k}^{j}}{\varepsilon_{M}}-\varepsilon_{M}\xi_{k}^{j}\boldsymbol{\Theta}[\ell]\left(1-\tfrac{\boldsymbol{C}_{k}[i,\ell]}{C_{k,\max}+\varepsilon_{C}}\right).\thinspace\mathrm{(\widetilde{LP4})}
\end{align*}
According to Definition~\ref{def:cost-function}, $\boldsymbol{C}_{k}[i,i]=0$
and $\boldsymbol{C}_{k}[i,\ell]>0$ for all $i\not=\ell$. The minimum
cost of this simpler LP (\ref{eq:Markov_optimization_problem_simpler})
is obtained when $\sum_{i=1}^{n_{\mathrm{bin}}}\sum_{\ell=1}^{n_{\mathrm{bin}}}\boldsymbol{C}_{k}[i,\ell]\boldsymbol{R}_{k}^{j}[i,\ell]=0$,
because the second term in the cost function is a constant. Therefore,
all the positive off-diagonal elements $\boldsymbol{M}_{k}^{j}[i,\ell]$
are equal to their respective lower bounds $\varepsilon_{M}\xi_{k}^{j}\boldsymbol{\Theta}[\ell]\left(1-\frac{\boldsymbol{C}_{k}[i,\ell]}{C_{k,\max}+\varepsilon_{C}}\right)$
in the optimal solution of the simpler LP (\ref{eq:Markov_optimization_problem_simpler}).
This optimal solution of the simpler LP (\ref{eq:Markov_optimization_problem_simpler})
is given by the Markov matrix $\boldsymbol{M}_{k}^{j}$ (\ref{eq:scalable_Markov1})-(\ref{eq:scalable_Markov2}). 

If the optimal solution of the simpler LP (\ref{eq:Markov_optimization_problem_simpler})
also satisfies the constraints (LP1,2) that we neglected previously,
then it is the optimal solution of the original LP (\ref{eq:Markov_optimization_problem}).
It follows from the construction of the diagonal elements in $\boldsymbol{M}_{k}^{j}$
(\ref{eq:scalable_Markov2}) that it satisfies (LP1). The diagonal
elements of $\boldsymbol{M}_{k}^{j}$ are given by:
\begin{equation}
\boldsymbol{M}_{k}^{j}[i,i]=1-\!\!\!\!\!\!\!\!\!\!\!\!\sum_{\ell\in\{\boldsymbol{A}_{k}^{j}[i,\ell]=1,i\not=\ell\}}\!\!\!\!\!\!\!\!\!\!\!\!\varepsilon_{M}\xi_{k}^{j}\boldsymbol{\Theta}[\ell]\left(1-\tfrac{\boldsymbol{C}_{k}[i,\ell]}{C_{k,\max}+\varepsilon_{C}}\right)\thinspace.\label{eq:scalable_optimal_diagonal}
\end{equation}
Note that the matrix $\boldsymbol{M}_{k}^{j}$ is a reversible Markov
matrix because of the symmetric cost matrix, i.e., $\boldsymbol{\Theta}[i]\boldsymbol{M}_{k}^{j}[i,\ell]=\boldsymbol{\Theta}[\ell]\boldsymbol{M}_{k}^{j}[\ell,i]=\varepsilon_{M}\xi_{k}^{j}\boldsymbol{\Theta}[i]\boldsymbol{\Theta}[\ell]\left(1-\frac{\boldsymbol{C}_{k}[i,\ell]}{C_{k,\max}+\varepsilon_{C}}\right)$
for all $i,\ell$. Hence, the matrix $\boldsymbol{M}_{k}^{j}$ also
satisfies (LP2) because:
\begin{align*}
\sum_{i=1}^{n_{\mathrm{bin}}}\boldsymbol{\Theta}[i]\boldsymbol{M}_{k}^{j}[i,\ell] & =\sum_{i=1}^{n_{\mathrm{bin}}}\boldsymbol{\Theta}[\ell]\boldsymbol{M}_{k}^{j}[\ell,i]\\
 & =\boldsymbol{\Theta}[\ell]\left(\sum_{i=1}^{n_{\mathrm{bin}}}\boldsymbol{M}_{k}^{j}[\ell,i]\right)=\boldsymbol{\Theta}[\ell]\thinspace.
\end{align*}
Therefore, the matrix $\boldsymbol{M}_{k}^{j}$ is the optimal solution
of the original LP (\ref{eq:Markov_optimization_problem}). \hfill$\blacksquare$\end{proof}

If the cost matrix $\boldsymbol{C}_{k}$ is symmetric, then the Markov
matrix $\boldsymbol{M}_{k}^{j}$ (\ref{eq:scalable_Markov1})-(\ref{eq:scalable_Markov2})
gives significant savings in computation time because each agent can
directly compute the row of the optimal Markov matrix. For example,
the computation times for all four cases in Fig.~\ref{fig:Multi-resolution-Images-of-Eiffel-Tower}
are less than $2$ minutes on a standard desktop computer. 

\begin{remark} \label{rem:Alternative_func} \textit{(Alternative
Functions for Constraints)} Note that our construction technique holds
even if the term $\left(1-\frac{\boldsymbol{C}_{k}[i,\ell]}{C_{k,\max}+\varepsilon_{C}}\right)$
in the constraint (LP4) in \textbf{Method~1} and (\ref{eq:scalable_Markov1})
in Corollary~\ref{cor:Markov-optimization-scalable} is replaced
by any monotonic function in $(0,1]$ that decreases with an increasing
$\boldsymbol{C}_{k}[i,\ell]$. Similarly, the term $\xi_{k}^{j}$
in the constraints (LP3,4) can be replaced by any monotonic function
in $(0,1]$ that decreases with a decreasing $\xi_{k}^{j}$. For example,
see Fig.~\ref{fig:estimation-error-plot-compare}(b) in Section~\ref{subsec:Numerical-Simulations-Estimation-error}.
\hfill$\Box$\end{remark}

\subsection{Construction of Fastest Mixing IMC}

In this subsection, we construct the fastest mixing IMC, where the
IMC's convergence rate to the rank one matrix $\boldsymbol{1}\boldsymbol{\Theta}$
is optimized. The convergence rate of HMC, with time-invariant Markov
matrix $\boldsymbol{M}$, is determined by the second largest eigenvalue
modulus (i.e., $\max_{r\in\{2,\ldots,n_{\textrm{bin}}\}}|\lambda_{r}(\boldsymbol{M})|$)
\cite{Ref:Boyd04Markov,Ref:Akar15}. On the other hand, the convergence
rate of IMC is determined by the coefficient of ergodicity \cite[pp. 137]{Ref:Seneta06}.
Since the first $n_{\mathrm{rec}}$ bins are recurrent bins, the Markov
matrix $\boldsymbol{M}_{k}^{j}$ can be decomposed as follows:

{\scriptsize{}
\begin{equation}
\boldsymbol{M}_{k}^{j}\!=\!\left[\!\begin{array}{cc}
\boldsymbol{M}_{k,\textrm{sub}}^{j} & \boldsymbol{0}^{n_{\textrm{rec}}\times(n_{\textrm{bin}}-n_{\textrm{rec}})}\\
\boldsymbol{M}_{k}^{j}[n_{\textrm{rec}}\!+\!1\!:\!n_{\textrm{bin}},1\!:\!n_{\textrm{rec}}] & \!\!\!\boldsymbol{M}_{k}^{j}[n_{\textrm{rec}}\!+\!1\!:\!n_{\textrm{bin}},n_{\textrm{rec}}\!+\!1\!:\!n_{\textrm{bin}}]
\end{array}\!\right],\label{eq:Markov_decomposition}
\end{equation}
}where the sub-matrix $\boldsymbol{M}_{k,\textrm{sub}}^{j}:=\boldsymbol{M}_{k}^{j}[1\!:\!n_{\textrm{rec}},1\!:\!n_{\textrm{rec}}]$
encapsulates the bin transition probabilities between the recurrent
bins. 

\begin{definition}\textit{\label{def:Coefficient-of-egodicity} }\cite[pp. 137--139]{Ref:Seneta06}\textit{
(Coefficient of Ergodicity)} For the stochastic matrix $\boldsymbol{M}_{k,\textrm{sub}}^{j}$,
the coefficient of ergodicity $\tau_{1}(\boldsymbol{M}_{k,\textrm{sub}}^{j})$
is defined as:
\begin{align}
 & \tau_{1}(\boldsymbol{M}_{k,\textrm{sub}}^{j})=\sup_{\boldsymbol{v}_{1},\boldsymbol{v}_{2},\boldsymbol{v}_{1}\not=\boldsymbol{v}_{2}}\frac{D_{\mathcal{L}_{1}}(\boldsymbol{v}_{1}\boldsymbol{M}_{k,\textrm{sub}}^{j},\boldsymbol{v}_{2}\boldsymbol{M}_{k,\textrm{sub}}^{j})}{D_{\mathcal{L}_{1}}(\boldsymbol{v}_{1},\boldsymbol{v}_{2})}\thinspace,\nonumber \\
 & \quad=1-\min_{i,\ell}\sum_{s=1}^{n_{\textrm{rec}}}\min\left(\boldsymbol{M}_{k,\textrm{sub}}^{j}[i,s],\boldsymbol{M}_{k,\textrm{sub}}^{j}[\ell,s]\right)\thinspace,\label{eq:coeff-ergodicity}
\end{align}
where $\boldsymbol{v}_{1}$, $\boldsymbol{v}_{2}$ are probability
row vectors in $\mathbb{R}^{n_{\textrm{rec}}}$ and $i,\ell,s\in\{1,\ldots,n_{\textrm{rec}}\}$.
\hfill$\Box$\end{definition}

We define $n_{k,\mathrm{dia}}^{j}$ as the graph diameter in the graph
conforming to the matrix $\boldsymbol{A}_{k,\textrm{sub}}^{j}:=\boldsymbol{A}_{k}^{j}[1\!:\!n_{\textrm{rec}},1\!:\!n_{\textrm{rec}}]$;
i.e., it is the greatest number of edges in a shortest path between
any pair of recurrent bins \cite{Ref:West01}. If $n_{k,\mathrm{dia}}^{j}>2$,
then there exist recurrent bins $B[i]$ and $B[\ell]$ such that either
$\boldsymbol{M}_{k,\textrm{sub}}^{j}[i,s]=0$ or $\boldsymbol{M}_{k,\textrm{sub}}^{j}[\ell,s]=0$
for all $s\in\{1,\ldots,n_{\textrm{rec}}\}$. Substituting these bins
into (\ref{eq:coeff-ergodicity}) show that $\tau_{1}(\boldsymbol{M}_{k,\textrm{sub}}^{j})=1$
when $n_{k,\mathrm{dia}}^{j}>2$. In order to avoid this trivial case,
we minimize the coefficient of ergodicity of the positive matrix $(\boldsymbol{M}_{k,\textrm{sub}}^{j})^{n_{k,\mathrm{dia}}^{j}}$
\cite[Theorem 8.5.2, pp. 516]{Ref:Horn85}. 

\begin{corollary}\label{cor:Markov-coeff-ergo} \textit{(Construction
of Fastest-mixing Markov Matrix)} The following convex optimization
problem is used instead of the LP (\ref{eq:Markov_optimization_problem})
in \textbf{Method~1}:
\begin{align}
\min\tau_{1}\left((\boldsymbol{M}_{k,\mathrm{sub}}^{j})^{n_{k,\mathrm{dia}}^{j}}\right)\thinspace,\label{eq:Markov_coeff_ergo}
\end{align}
subject to $\mathrm{(LP1-4)}$ in (\ref{eq:Markov_optimization_problem}),
where $\tau_{1}$ is defined in Definition~\ref{def:Coefficient-of-egodicity}.
\end{corollary}

\begin{proof} The cost function $\tau_{1}\left((\boldsymbol{M}_{k,\textrm{sub}}^{j})^{n_{k,\mathrm{dia}}^{j}}\right)$
is a convex function of the stochastic matrix $\boldsymbol{M}_{k,\textrm{sub}}^{j}$
because it can be expressed as \cite[Lemma 4.3, pp. 139]{Ref:Seneta06}:
\[
\tau_{1}\left((\boldsymbol{M}_{k,\textrm{sub}}^{j})^{n_{k,\mathrm{dia}}^{j}}\right)=\sup_{\|\boldsymbol{\delta}\|_{2}=1,\thinspace\boldsymbol{\delta}\boldsymbol{1}=0}\left\Vert \boldsymbol{\delta}\cdot(\boldsymbol{M}_{k,\textrm{sub}}^{j})^{n_{k,\mathrm{dia}}^{j}}\right\Vert _{1}\thinspace,
\]
where $\boldsymbol{\delta}=\textrm{const}(\boldsymbol{v}_{1}-\boldsymbol{v}_{2})$
is a row vector in $\mathbb{R}^{n_{\textrm{rec}}}$. Hence, (\ref{eq:Markov_coeff_ergo})
is a convex optimization problem and the family of Markov matrices
$\boldsymbol{Q}_{k}^{j}$ (\ref{eq:Q_off_diag})-(\ref{eq:Q_diag})
is a subset of its feasible set. \hfill$\blacksquare$ \end{proof} 

\subsection{Condition for Escaping Transient Bins \label{subsec:Condition-for-Escaping-transient-bins}}

In this subsection, we discuss the condition for escaping transient
bins.

\begin{figure}[h]
\begin{centering}
\includegraphics[width=2.5in]{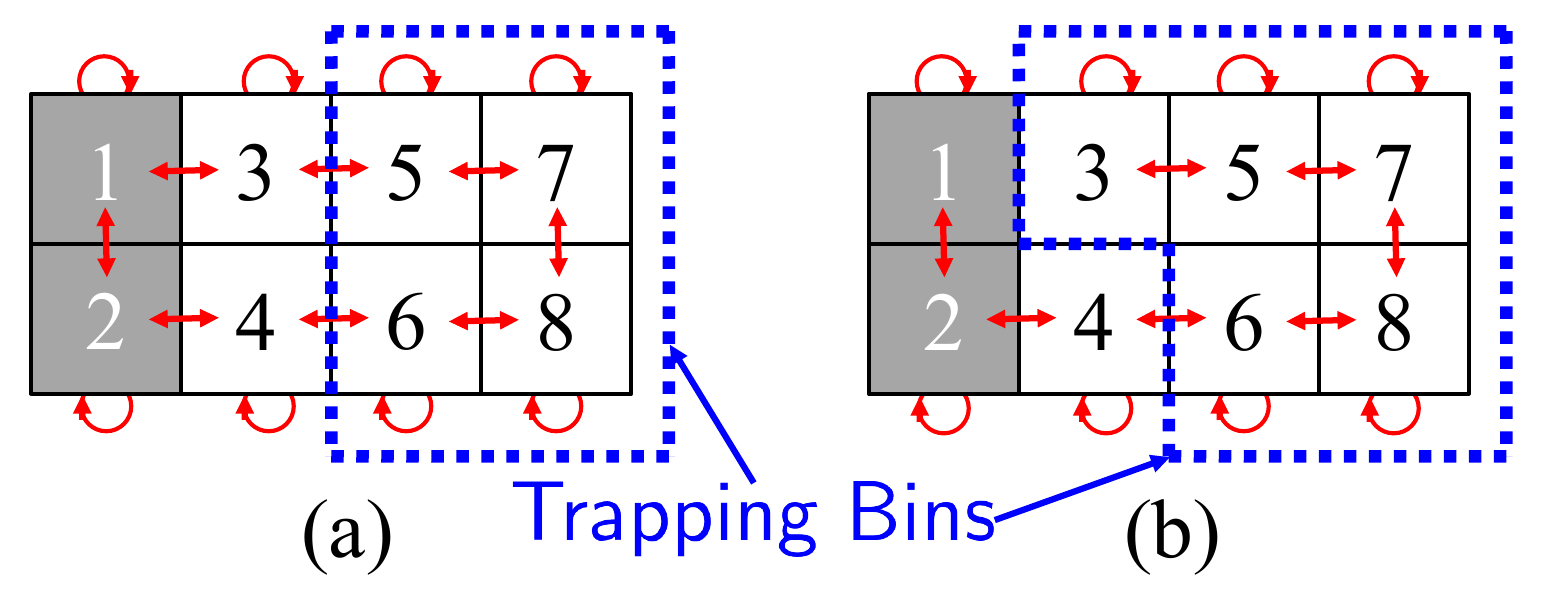}
\par\end{centering}
\caption{In this example, the bins 1 and 2 are recurrent bins. The allowed
transitions (motion constraints) are shown in red. The trapping bins
for the two cases are enclosed in blue. \label{fig:trapping_bins}}
\end{figure}

\begin{definition} \label{def:tapping_bin} \textit{(Trapping Bins)}
If an agent is inside a transient bin ($\boldsymbol{\Theta}[i]=0$)
and its motion constraints matrix $\boldsymbol{A}_{k}^{j}$ only allows
transitions to other transient bins, then that transient bin is called
a trapping bin. This agent is trapped in this bin because the Markov
matrix does not allow transitions out of this bin. Let $\boldsymbol{\mathcal{T}}_{k}^{j}$
represent the set of trapping bins for the $j^{\textrm{th}}$ agent
at the $k^{\textrm{th}}$ time instant. For example, see Fig.~\ref{fig:trapping_bins}.
\hfill$\Box$ \end{definition}

Since the irreducible motion constraints matrices $\boldsymbol{A}_{k}^{j}$
are known a priori (Assumption~\ref{assump:prior-knowledge} and
Definition~\ref{def:Motion-Constraints}), the deterministic path
for exiting the set of trapping bins is stored on board each agent.
For each trapping bin $B[i]\in\boldsymbol{\mathcal{T}}_{k}^{j}$,
the $j^{\textrm{th}}$ agent transitions to a transient bin $\Psi_{k}^{j}[i]$,
chosen a priori, such that the transition from bin $B[i]$ to bin
$\Psi_{k}^{j}[i]$ is allowed by motion constraints. This deterministic
path, which is chosen a priori, ensures that the agent exits the set
of trapping bins as soon as possible. This bin $\Psi_{k}^{j}[i]$
has to be chosen on a case-by-case basis depending on the motion constraints
matrix $\boldsymbol{A}_{k}^{j}$. For example, in Fig.~\ref{fig:trapping_bins},
for the trapping bin 5, the best option is bin 3 in case (a) and bin
7 in case (b). Therefore, the agent can follow this path to deterministically
exit the set of trapping bins in finite time instants. 

If an agent is in a transient bin, but not in a trapping bin, then
its motion constraint matrix allows transitions to some recurrent
bins. We can speed up the process of exiting this transient bin by
forcing the agent to transition to any reachable recurrent bin, with
equal probability, during the current time instant. Thus the agent
transitions from its current transient bin to a recurrent bin in one
time instant. 

The matrix $\boldsymbol{S}_{k}^{j}\in\mathbb{R}^{n_{\textrm{bin}}\times n_{\textrm{bin}}}$
encapsulates the condition for escaping transient bins. If $B[i]$
is a transient bin (i.e., $\boldsymbol{\Theta}[i]=0$), then each
element in the corresponding row $\boldsymbol{S}_{k}^{j}[i,1\!:\!n_{\textrm{bin}}]$
is given by:
\begin{align}
\boldsymbol{S}_{k}^{j}[i,\ell] & =\begin{cases}
1 & \textrm{ \textbf{if} }B[i]\in\boldsymbol{\mathcal{T}}_{k}^{j}\textrm{ and }B[\ell]=\Psi_{k}^{j}[i]\\
\frac{1}{n_{k,i}^{j}} & \begin{array}{l}
\textrm{\textbf{if} }B[i]\not\in\boldsymbol{\mathcal{T}}_{k}^{j}\textrm{ and }\boldsymbol{A}_{k}^{j}[i,\ell]=1\\
\textrm{ and }\boldsymbol{\Theta}[\ell]>0
\end{array}\\
0 & \textrm{ \textbf{otherwise}}
\end{cases},\label{eq:S_matrix}
\end{align}
where $n_{k,i}^{j}$ is the number of recurrent bins that the $j^{\textrm{th}}$
agent can transition to, from bin $B[i]$ at the $k^{\textrm{th}}$
time instant. This condition is used only when the agent is in a transient
bin, as shown in \textbf{Method~1}. In Section \ref{sec:Convergence-Analysis},
we show that the agent exits the set of transient bins within finite
time instants due to this condition. 

\section{PSG\textendash IMC Algorithm for Shape Formation \label{sec:PSG-IMC-pattern-formation}}

In this section, we first state the PSG\textendash IMC algorithm for
shape formation and then present its convergence analysis and its
property of robustness. 

\begin{algorithm}[h]
\begin{tabular}{cl}
\hline 
\multicolumn{2}{l}{\textbf{Method 2}: PSG\textendash IMC Algorithm for Shape Formation}\tabularnewline
\hline 
\hline 
{\small{}1:} & \multicolumn{1}{l}{{\small{}One iteration of $j^{\textrm{th}}$ agent during $k^{\textrm{th}}$
time instant, }}\tabularnewline
{\small{} } & {\small{}where the $j^{\textrm{th}}$ agent is in bin $B[i]$}\tabularnewline
{\small{}2:} & {\small{}Given $\boldsymbol{\Theta}$, $\boldsymbol{C}_{k}$, $\boldsymbol{A}_{k}^{j}$,
$\boldsymbol{\mu}_{k}^{j}$, $\varepsilon_{M}$, $\varepsilon_{C}$,
$\tau^{j}$, and $\beta^{j}$ }\tabularnewline
{\small{}3:} & \textbf{\small{}if}{\small{} }\textbf{\small{}$\boldsymbol{\Theta}[i]=0$}{\small{},
}\textbf{\small{}then}{\small{} }\tabularnewline
{\small{}4:} & \textbf{\small{}$\qquad$}{\small{}Compute $\boldsymbol{S}_{k}^{j}[i,1\!:\!n_{\textrm{bin}}]$
using (\ref{eq:S_matrix})}\tabularnewline
{\small{}5:} & \textbf{\small{}$\qquad$}{\small{}Generate a random number $z\in\textrm{unif}[0;1]$ }\tabularnewline
{\small{}6:} & \textbf{\small{}$\qquad$}{\small{}Select bin $B[q]$ such that }\tabularnewline
 & \textbf{\small{}$\qquad$$\qquad$}{\small{}$\sum_{\ell=1}^{q-1}\boldsymbol{S}_{k}^{j}[i,\ell]\leq z<\sum_{\ell=1}^{q}\boldsymbol{S}_{k}^{j}[i,\ell]$}\tabularnewline
{\small{}7:} & \textbf{\small{}else}{\small{} }\tabularnewline
{\small{}8:} & \textbf{\small{}$\qquad$}{\small{}Compute the feedback gain $\xi_{k}^{j}$
using (\ref{eq:tuning_eqn}) }\tabularnewline
{\small{}9:} & {\small{}$\qquad$Compute $\boldsymbol{M}_{k}^{j}[i,1\!:\!n_{\textrm{bin}}]$
using Corollary~\ref{cor:Markov-optimization-scalable}}\tabularnewline
 & {\small{}$\qquad$$\qquad$or compute $\boldsymbol{M}_{k}^{j}$ using
}\textbf{\small{}Method~1}\tabularnewline
{\small{}10:} & \textbf{\small{}$\qquad$}{\small{}Compute the term $\eta_{k,i}^{j}$
using (\ref{eq:eta_ki}) }\tabularnewline
{\small{}11:} & \textbf{\small{}$\qquad$}{\small{}Compute $\boldsymbol{P}_{k}^{j}[i,1\!:\!n_{\textrm{bin}}]$
using (\ref{eq:Pil}) and (\ref{eq:Pii}) }\tabularnewline
{\small{}12:} & \textbf{\small{}$\qquad$}{\small{}Generate a random number $z\in\textrm{unif}[0;1]$ }\tabularnewline
{\small{}13:} & \textbf{\small{}$\qquad$}{\small{}Select bin $B[q]$ such that }\tabularnewline
 & \textbf{\small{}$\qquad$$\qquad$}{\small{}$\sum_{\ell=1}^{q-1}\boldsymbol{P}_{k}^{j}[i,\ell]\leq z<\sum_{\ell=1}^{q}\boldsymbol{P}_{k}^{j}[i,\ell]$}\tabularnewline
{\small{}14:} & \textbf{\small{}end if}{\small{} }\tabularnewline
{\small{}15:} & {\small{}Go to bin $B[q]$ while avoiding collisions }\tabularnewline
\hline 
\end{tabular}
\end{algorithm}

The pseudo-code for the PSG\textendash IMC algorithm for shape formation
is given in \textbf{Method~2}, whose key steps are shown in Fig.~\ref{fig:Flowchart-PSGIMC}.
At the start, the $j^{\textrm{th}}$ agent knows the desired formation
shape $\boldsymbol{\Theta}$, the time-varying cost matrix $\boldsymbol{C}_{k}$,
and its time-varying motion constraints matrix $\boldsymbol{A}_{k}^{j}$
(Assumption~\ref{assump:prior-knowledge}). During each iteration,
the agent determines the bin it belongs to (Assumption~\ref{assump:r_known},
here we assume that the $j^{\textrm{th}}$ agent is in bin $B[i]$)
and the current swarm distribution $\boldsymbol{\mathcal{\mu}}_{k}^{j}$
from Section~\ref{subsec:Generating-Feedback-swarm-dist} (lines
1\textendash 2). 

If the agent is in a transient bin (line 3), then it uses inverse
transform sampling (Remark~\ref{rem:Inverse-Transform-Sampling}
in Appendix) to select the next bin $B[q]$ from the corresponding
row of the matrix $\boldsymbol{S}_{k}^{j}[i,1\!:\!n_{\textrm{bin}}]$
(lines 4\textendash 6).

Otherwise, the agent first computes the HD-based feedback gain $\xi_{k}^{j}$
(line 8). If the cost matrix $\boldsymbol{C}_{k}$ is symmetric, then
the agent can directly compute the row $\boldsymbol{M}_{k}^{j}[i,1\!:\!n_{\textrm{bin}}]$
using Corollary~\ref{cor:Markov-optimization-scalable} (line 9).
Otherwise, the agent computes the entire Markov matrix $\boldsymbol{M}_{k}^{j}$
using \textbf{Method~1} (line 9). 

In order to avoid undesirable transitions from bins that are deficient
in agents (i.e., where $\boldsymbol{\Theta}[i]>\boldsymbol{\mathcal{\mu}}_{k}^{j}[i]$),
the agent modifies its Markov matrix row $\boldsymbol{M}_{k}^{j}[i,1\!:\!n_{\textrm{bin}}]$
as follows:
\begin{align}
\boldsymbol{P}_{k}^{j}[i,\ell] & =\left(1-\eta_{k,i}^{j}\right)\boldsymbol{M}_{k}^{j}[i,\ell]\thinspace,\qquad\forall i\not=\ell\label{eq:Pil}\\
\boldsymbol{P}_{k}^{j}[i,i] & =\left(1-\eta_{k,i}^{j}\right)\boldsymbol{M}_{k}^{j}[i,i]+\eta_{k,i}^{j}\thinspace,\label{eq:Pii}\\
\textrm{where }\eta_{k,i}^{j} & =\exp(-\tau^{j}k)\frac{\exp\left(\beta^{j}(\boldsymbol{\Theta}[i]-\boldsymbol{\mathcal{\mu}}_{k}^{j}[i])\right)}{\exp\left(\beta^{j}|\boldsymbol{\Theta}[i]-\boldsymbol{\mathcal{\mu}}_{k}^{j}[i]|\right)}\thinspace,\label{eq:eta_ki}
\end{align}
and $\tau^{j}$ and $\beta^{j}$ are time-invariant positive constants
(lines 10\textendash 11). Then, the agent uses inverse transform sampling
(Remark~\ref{rem:Inverse-Transform-Sampling} in Appendix) to select
the next bin $B[q]$ from the bin transition probabilities $\boldsymbol{P}_{k}^{j}[i,1\!:\!n_{\textrm{bin}}]$
(lines 12\textendash 13). Finally, the agent goes to the selected
bin $B[q]$ from the current bin $B[i]$ in a collision-free manner
using the lower-level guidance and control algorithm discussed in
Remark~\ref{rem:collision-free-motion} in Appendix (line 15).

\begin{remark} \label{rem:explain_eta_term} \textit{(Explanation
of the term $\eta_{k,i}^{j}$)} The term $\eta_{k,i}^{j}$ (\ref{eq:eta_ki})
is designed to suppress the transition probabilities from a bin that
is deficient in agents; i.e., if $\boldsymbol{\Theta}[i]>\boldsymbol{\mathcal{\mu}}_{k}^{j}[i]$,
then $\eta_{k,i}^{j}=\exp(-\tau^{j}k)$. Its effect decreases with
increasing time instants. The design variables $\beta^{j}$ and $\tau^{j}$
dictate the amplitude and time constant of this suppression. It is
shown later in Section~\ref{subsec:Numerical-Simulations-fine-resolution}
that the undesirable transitions, suppressed using this term, greatly
reduce the total number of transitions, which in turn significantly
improves the convergence error.\hfill$\Box$\end{remark}

\begin{remark} \label{rem:PSG-IMC-distributed} \textit{(PSG\textendash IMC
Algorithm is a Distributed Algorithm)} Under Assumptions~\ref{assump:r_known}
and \ref{assump:Strongly-connected-communication}, the agent perceives
the bin it belongs to and estimates the current swarm distribution
in a distributed manner. The remaining terms in lines 1\textendash 2
are known a priori. Lines 3\textendash 14 are executed individually
by each agent. Finally, in line 15, the agent only needs to communicate
with its neighboring agents as shown in Remark~\ref{rem:collision-free-motion}
in Appendix. Thus, the steps in \textbf{Method~2} can be accomplished
in a distributed manner.\hfill$\Box$\end{remark}

\subsection{Main Result: Convergence Analysis \label{sec:Convergence-Analysis}}

In this subsection, we prove that the swarm distribution $\boldsymbol{\mu}_{k}^{\star}$
converges to the desired formation shape $\boldsymbol{\Theta}$ with
prescribed convergence errors using the PSG\textendash IMC algorithm
for shape formation given in \textbf{Method~2}. Unlike the convergence
proof for HMC, which is a direct application of the Perron\textendash Frobenius
theorem, the convergence proof for IMC is rather involved (e.g., see
\cite{Ref:Seneta06,Ref:Touri11}). We first state an assumption on
$\xi_{k}^{j}$ for this subsection.

\begin{assumption} \label{assump:xi_greater_xi_des} \textit{(Minimum
Value of $\xi_{k}^{j}$)} If $\xi_{k}^{j}<\xi_{\textrm{des}}$, then
the current swarm distribution is sufficiently close to the desired
formation (see Definition~\ref{def:Feedback-gain}). Moreover, the
Markov matrices in \textbf{Method~1} become an identity matrix, hence
the agents do not transition any more. Therefore, the swarm has converged
to the desired formation and no further convergence is necessary.

Therefore, in this subsection, we assume that the swarm has not converged
to the desired formation (i.e., $\xi_{k}^{j}\geq\xi_{\textrm{des}}$).
\hfill$\Box$ \end{assumption}

We first show that agents in recurrent bins transition using to the
modified Markov matrix $\boldsymbol{P}_{k}^{j}$ derived from (\ref{eq:Pil})-(\ref{eq:Pii}). 

\begin{theorem} \label{thm:Markov_Pkj} According to \textbf{Method~2},
if an agent is in a recurrent bin, then it transitions using the following
modified Markov matrix $\boldsymbol{P}_{k}^{j}$ from (\ref{eq:Pil})-(\ref{eq:Pii}):
\begin{align}
\boldsymbol{P}_{k}^{j} & =\left(\mathbf{I}-\boldsymbol{D}_{k}^{j}\right)\boldsymbol{M}_{k}^{j}+\boldsymbol{D}_{k}^{j}\thinspace,\label{eq:Pk}
\end{align}
where $\boldsymbol{D}_{k}^{j}=\textrm{diag}\left(\eta_{k,1}^{j},\ldots,\eta_{k,n_{\mathrm{bin}}}^{j}\right)$.
The Markov matrix $\boldsymbol{P}_{k}^{j}$ is row stochastic (i.e.,
$\boldsymbol{P}_{k}^{j}\boldsymbol{1}=\boldsymbol{1}$), asymptotically
homogeneous with respect to $\boldsymbol{\Theta}$ (i.e., $\lim_{k\rightarrow\infty}\boldsymbol{\Theta}\boldsymbol{P}_{k}^{j}=\boldsymbol{\Theta}$),
and only allows transitions into recurrent bins. \end{theorem}

\begin{proof} The modified Markov matrix $\boldsymbol{P}_{k}^{j}$
(\ref{eq:Pk}) is derived from (\ref{eq:Pil})\textendash (\ref{eq:eta_ki}).
It follows from lines 3 and 6 of \textbf{Method~2} that the agent
uses the Markov matrix $\boldsymbol{P}_{k}^{j}$ to transition if
and only if it is in a recurrent bin (i.e., $\boldsymbol{\Theta}[i]>0$). 

The matrix $\boldsymbol{P}_{k}^{j}$ is row stochastic because $\boldsymbol{M}_{k}^{j}\boldsymbol{1}=\boldsymbol{1}$.
The matrix $\boldsymbol{M}_{k}^{j}$ has $\boldsymbol{\Theta}$ as
its stationary distribution for all $k\in\mathbb{N}$. It follows
from the definition of the term $\eta_{k,i}^{j}$ (\ref{eq:eta_ki})
that $\lim_{k\rightarrow\infty}\boldsymbol{D}_{k}^{j}=\boldsymbol{0}^{n_{\mathrm{bin}}\times n_{\mathrm{bin}}}$,
because $\lim_{k\rightarrow\infty}\exp(-\tau^{j}k)=0$. Therefore
$\lim_{k\rightarrow\infty}\boldsymbol{P}_{k}^{j}=\lim_{k\rightarrow\infty}\boldsymbol{M}_{k}^{j}$.
Hence, the sequence of matrices $\boldsymbol{P}_{k}^{j}$ is asymptotically
homogeneous with respect to $\boldsymbol{\Theta}$ because $\lim_{k\rightarrow\infty}\boldsymbol{\Theta}\boldsymbol{P}_{k}^{j}=\lim_{k\rightarrow\infty}\boldsymbol{\Theta}\boldsymbol{M}_{k}^{j}=\boldsymbol{\Theta}$
(see Definition~\ref{def:Asymp-homo} in Appendix). 

Note that the element $\boldsymbol{P}_{k}^{j}[i,\ell]>0$ if and only
if the corresponding element $\boldsymbol{M}_{k}^{j}[i,\ell]>0$ for
all $i,\ell\in\{1,\ldots,n_{\textrm{bin}}\}$ and $k\in\mathbb{N}$.
Therefore, like the matrix $\boldsymbol{M}_{k}^{j}$, the matrix $\boldsymbol{P}_{k}^{j}$
only allows transitions into recurrent bins. \hfill$\blacksquare$
\end{proof} 

We now show that all the agents leave the transient bins and enter
the recurrent bins in finite time instants. 

\begin{theorem} \label{thm:finite-time-conv} According to \textbf{Method~2},
each agent is in a recurrent bin by the $T^{\textrm{th}}$ time instant,
where $T\leq(n_{\mathrm{bin}}-n_{\mathrm{rec}}+1)$. Once an agent
is inside a recurrent bin, it always remains within the set of recurrent
bins. \end{theorem} 

\begin{proof} If an agent is in a recurrent bin, then it follows
from Theorem~\ref{thm:Markov_Pkj} that it cannot transition to any
transient bin. 

If the agent is in a trapping bin, then the matrix $\boldsymbol{S}_{k}^{j}$
(\ref{eq:S_matrix}) ensures that the agent exits the set of trapping
bins as soon as possible in a deterministic manner. Therefore, the
maximum number of steps inside the set of trapping bins is upper bounded
by the number of transient bins $(n_{\mathrm{bin}}-n_{\mathrm{rec}})$. 

If an agent is in a transient bin, but not in a trapping bin, then
the matrix $\boldsymbol{S}_{k}^{j}$ (\ref{eq:S_matrix}) ensures
that the agent transitions to a recurrent bin in one time instant.
Hence each agent enters a recurrent bin in at most $(n_{\mathrm{bin}}-n_{\mathrm{rec}}+1)$
time instants. \hfill$\blacksquare$ \end{proof} 

Consider a probability (row) vector $\boldsymbol{x}_{k}^{j}\in\mathbb{R}^{n_{\textrm{bin}}}$,
which denotes the probability mass function (PMF) of the predicted
bin position of the $j^{\textrm{th}}$ agent at the $k^{\textrm{th}}$
time instant. Each element $\boldsymbol{x}_{k}^{j}[i]$ gives the
probability that the $j^{\textrm{th}}$ agent is in bin $B[i]$ at
the $k^{\textrm{th}}$ time instant: 
\begin{equation}
\boldsymbol{x}_{k}^{j}[i]=\mathbb{P}(\boldsymbol{r}_{k}^{j}[i]=1),\quad\forall i\in\{1,\ldots,n_{\mathrm{bin}}\}\thinspace.\label{eq:def_prob_vector}
\end{equation}
We now discuss convergence of each agent's predicted bin position
$\boldsymbol{x}_{k}^{j}$ to the desired formation $\boldsymbol{\Theta}$. 

\begin{theorem} \label{thm:Convergence-of-imhomo-MC} \textit{(Convergence
of IMC)} According to \textbf{Method~2}, each agent's time evolution
of the vector $\boldsymbol{x}_{k}^{j}$ converges pointwise to the
desired stationary distribution $\boldsymbol{\Theta}$ irrespective
of the initial condition, i.e., $\lim_{k\rightarrow\infty}\boldsymbol{x}_{k}^{j}=\boldsymbol{\Theta}$
pointwise for all agents. \end{theorem} 

\begin{proof} It follows from Theorem~\ref{thm:finite-time-conv}
that all agents are always in the set of recurrent bins from the $T^{\textrm{th}}$
time instant onwards. Since the first $n_{\mathrm{rec}}$ bins are
recurrent bins, we decompose the vector $\boldsymbol{x}_{k}^{j}=[\bar{\boldsymbol{x}}_{k}^{j},\thinspace0,\ldots,0]$
for all $k\geq T$, where the probability row vector $\bar{\boldsymbol{x}}_{k}^{j}:=[\boldsymbol{x}_{k}^{j}[1],\ldots,\boldsymbol{x}_{k}^{j}[n_{\mathrm{rec}}]]\in\mathbb{R}^{n_{\mathrm{rec}}}$
denotes the agent's PMF over the set of recurrent bins. Similarly,
we decompose $\boldsymbol{\Theta}=[\bar{\boldsymbol{\Theta}},\thinspace0,\ldots,0]$,
where $\bar{\boldsymbol{\Theta}}:=\left[\boldsymbol{\Theta}[1],\ldots,\boldsymbol{\Theta}[n_{\mathrm{rec}}]\right]$.
Note that convergence of $\bar{\boldsymbol{x}}_{k}^{j}$ to $\bar{\boldsymbol{\Theta}}$,
implies the convergences of $\boldsymbol{x}_{k}^{j}$ to $\boldsymbol{\Theta}$. 

According to Theorem~\ref{thm:Markov_Pkj}, the time evolution of
the PMF vector $\bar{\boldsymbol{x}}_{k}^{j}$ is given by:
\begin{align}
\bar{\boldsymbol{x}}_{k+1}^{j} & =\bar{\boldsymbol{x}}_{k}^{j}\boldsymbol{P}_{k,\textrm{sub}}^{j}\thinspace,\qquad\forall k\geq T\thinspace,\label{eq:proof_time_evolution}
\end{align}
where the row-stochastic sub-matrix $\boldsymbol{P}_{k,\textrm{sub}}^{j}:=\boldsymbol{P}_{k}^{j}[1\!:\!n_{\textrm{rec}},1\!:\!n_{\textrm{rec}}]$
encapsulates the bin transition probabilities between the recurrent
bins. The matrix $\boldsymbol{P}_{k,\textrm{sub}}^{j}$, like the
matrix $\boldsymbol{M}_{k,\textrm{sub}}^{j}$ in (\ref{eq:Markov_decomposition}),
is irreducible because the matrix $\boldsymbol{A}_{k,\textrm{sub}}^{j}$
is irreducible (Definition~\ref{def:Motion-Constraints}).

It follows from (\ref{eq:tuning_eqn}) that $D_{H}(\boldsymbol{\Theta},\boldsymbol{\mathcal{\mu}}_{k}^{j})=1$
if and only if $\boldsymbol{\mathcal{\mu}}_{k}^{j}[i]=0$ for all
recurrent bins $i\in\{1,\ldots,n_{\textrm{rec}}\}$, because $\boldsymbol{\Theta}[i]>0$
only in recurrent bins (Definition~\ref{def:Desired-formation}).
Therefore, the feedback gain $\xi_{k}^{j}<1$ for all time instant
$k\geq T$ because all agents are in recurrent bins. Hence, the diagonal
elements $\boldsymbol{M}_{k,\textrm{sub}}^{j}[i,i]$ and $\boldsymbol{P}_{k,\textrm{sub}}^{j}[i,i]$
are positive for all $i\in\{1,\ldots,n_{\textrm{rec}}\}$ and $k\geq T$.

The overall time evolution of the agent's PMF vector is given by the
IMC for all $r>T$:
\begin{align}
\bar{\boldsymbol{x}}_{r}^{j} & =\bar{\boldsymbol{x}}_{T}^{j}\boldsymbol{P}_{T,\textrm{sub}}^{j}\boldsymbol{P}_{T+1,\textrm{sub}}^{j}\ldots\boldsymbol{P}_{r-1,\textrm{sub}}^{j}=\bar{\boldsymbol{x}}_{T}^{j}\boldsymbol{U}_{T,r}^{j}.\label{eq:proof_overall_evolution}
\end{align}
We now show that this forward matrix product $\boldsymbol{U}_{T,r}^{j}$
is strongly ergodic (see Definition~\ref{def:Strong-ergodic} in
Appendix) and $\bar{\boldsymbol{\Theta}}$ is its unique limit vector
(i.e., $\lim_{r\rightarrow\infty}\boldsymbol{U}_{T,r}=\boldsymbol{1\bar{\Theta}}$).

The matrix $\boldsymbol{U}_{T,r}^{j}$ is a product of nonnegative
matrices, hence it is a nonnegative matrix. If $\boldsymbol{P}_{k,\textrm{sub}}^{j}[i,\ell]>0$
for some $k\in\{T,\ldots,r-1\}$ and $i,\ell\in\{1,\ldots,n_{\textrm{rec}}\}$,
then the corresponding element $\boldsymbol{U}_{T,r}^{j}[i,\ell]>0$
because, as shown below, the value of $\boldsymbol{U}_{T,r}^{j}[i,\ell]$
is lower bounded by the product of positive diagonal elements and
$\boldsymbol{P}_{k,\textrm{sub}}^{j}[i,\ell]$:
\begin{align}
 & \boldsymbol{U}_{T,r}^{j}[i,\ell]\geq\boldsymbol{P}_{T,\textrm{sub}}^{j}[i,\ell]\Bigl(\prod_{q=T+1}^{r-1}\boldsymbol{P}_{q,\textrm{sub}}^{j}[\ell,\ell]\Bigr)\nonumber \\
 & +\sum_{s=T+1}^{r-2}\left(\Bigl(\prod_{q=T}^{s-1}\boldsymbol{P}_{q,\textrm{sub}}^{j}[i,i]\Bigr)\boldsymbol{P}_{s,\textrm{sub}}^{j}[i,\ell]\Bigl(\prod_{q=s+1}^{r-1}\boldsymbol{P}_{q,\textrm{sub}}^{j}[\ell,\ell]\Bigr)\right)\nonumber \\
 & +\Bigl(\prod_{q=T}^{r-2}\boldsymbol{P}_{q,\textrm{sub}}^{j}[i,i]\Bigr)\boldsymbol{P}_{r-1,\textrm{sub}}^{j}[i,\ell]\thinspace,\qquad\qquad\textrm{if }i\not=\ell\thinspace,\label{eq:UTr_il}\\
 & \boldsymbol{U}_{T,r}^{j}[i,i]\geq\Bigl(\prod_{q=T}^{r-1}\boldsymbol{P}_{q,\textrm{sub}}^{j}[i,i]\Bigr)\thinspace,\qquad\qquad\quad\textrm{if }i=\ell\thinspace.\label{eq:UTr_il_diag}
\end{align}
Therefore, the matrix $\boldsymbol{U}_{T,r}^{j}$, like the matrix
$\boldsymbol{P}_{k,\textrm{sub}}^{j}$, is irreducible because $\boldsymbol{U}_{T,r}^{j}[i,\ell]>0$
if $\boldsymbol{P}_{k,\textrm{sub}}^{j}[i,\ell]>0$ for all $i,\ell\in\{1,\ldots,n_{\textrm{rec}}\}$.
Since the irreducible matrix $\boldsymbol{U}_{T,r}^{j}$ has positive
diagonal elements (\ref{eq:UTr_il_diag}), it is a primitive matrix
\cite[Lemma 8.5.4, pp. 516]{Ref:Horn85}.

Some of off-diagonal elements in $\boldsymbol{M}_{k,\textrm{sub}}^{j}$
and $\boldsymbol{P}_{k,\textrm{sub}}^{j}$ are zero due to the constraints
(CS1,2) in Theorem~\ref{thm:Markov-optimization}. The lower bound
$\gamma^{j}$, which is independent of $k$, for the remaining positive
elements in $\boldsymbol{P}_{k,\textrm{sub}}^{j}$ is given by the
constraint (LP4) in Theorem~\ref{thm:Markov-optimization}, the lower
bound of $\xi_{k}^{j}$ in Assumption~\ref{assump:xi_greater_xi_des},
and the upper bound of the term $\eta_{k,i}^{j}$ (\ref{eq:eta_ki}):
\begin{align}
 & \gamma^{j}=\left(1-\exp(-\tau^{j}T)\right)\xi_{\textrm{des}}\varepsilon_{M}(\min{}^{+}\boldsymbol{\Theta})\left(1-\tfrac{C_{\max}}{C_{\max}+\varepsilon_{C}}\right)\nonumber \\
 & \qquad\qquad\leq\min_{i,\thinspace\ell}{}^{+}\boldsymbol{P}_{k,\textrm{sub}}^{j}[i,\ell]\thinspace,\qquad k\geq T\thinspace,\label{eq:gamma_value}
\end{align}
where $\min{}^{+}$ refers to the minimum of the positive elements
and $C_{\max}=\max_{k\in\mathbb{N}}C_{k,\max}$. It follows from Theorem~\ref{thm:Markov_Pkj}
that the sequence of matrices $\boldsymbol{P}_{k,\textrm{sub}}^{j}$,
$k\geq T$ is asymptotically homogeneous with respect to $\bar{\boldsymbol{\Theta}}$.
Since (i) the forward matrix product $\boldsymbol{U}_{T,r}$ is primitive
and (ii) there exists $\gamma^{j}$ (independent of $k$), it follows
from Theorem~\ref{thm:asymp-homo-ns-strong-ergo} in Appendix that
the forward matrix product $\boldsymbol{U}_{T,r}^{j}$ is strongly
ergodic. Since (i) the matrices $\boldsymbol{P}_{k,\textrm{sub}}^{j},\thinspace k\geq T$
are irreducible and (ii) there exists $\gamma^{j}$ (independent of
$k$), it follows from Theorem~\ref{thm:limit-of-strong-ergo} in
Appendix that the limit vector $\boldsymbol{e}=\bar{\boldsymbol{\Theta}}$.
Since (iii) $\boldsymbol{U}_{T,r}^{j}$ is strongly ergodic, it follows
from Corollary~\ref{cor:limit-of-strong-ergo} in Appendix that the
unique limit vector is given by $\bar{\boldsymbol{\Theta}}$ (i.e.,
$\lim_{r\rightarrow\infty}\boldsymbol{U}_{T,r}^{j}=\boldsymbol{1}\bar{\boldsymbol{\Theta}}$).
Hence, each agent's PMF vector converges to: 
\[
\lim_{r\rightarrow\infty}\bar{\boldsymbol{x}}_{r}^{j}=\lim_{r\rightarrow\infty}\bar{\boldsymbol{x}}_{T}^{j}\boldsymbol{U}_{T,r}^{j}=\bar{\boldsymbol{x}}_{T}^{j}\boldsymbol{1}\bar{\boldsymbol{\Theta}}=\bar{\boldsymbol{\Theta}}\thinspace.
\]
Therefore, $\lim_{k\rightarrow\infty}\boldsymbol{x}_{k}^{j}=\boldsymbol{\Theta}$
pointwise for all agents. \hfill$\blacksquare$ \end{proof}

\begin{theorem} \label{thm:DL1_x_Theta} Since $\lim_{r\rightarrow\infty}\boldsymbol{U}_{T,r}^{j}=\boldsymbol{1}\bar{\boldsymbol{\Theta}}$,
for all $\varepsilon_{\mathrm{lim}}>0$, there exists a $k_{\epsilon,\mathrm{lim}}^{j}\in\mathbb{N}$
such that $D_{\mathcal{L}_{1}}(\bar{\boldsymbol{\Theta}}\boldsymbol{U}_{T,r}^{j},\bar{\boldsymbol{\Theta}})\leq\varepsilon_{\mathrm{lim}}$
for all $r\geq k_{\epsilon,\mathrm{lim}}^{j}$. The convergence error
between the $j^{\textrm{th}}$ agent's PMF vector $\boldsymbol{x}_{r}^{j}$
and the desired formation $\boldsymbol{\Theta}$ is bounded for all
$r\geq k_{\epsilon,\mathrm{lim}}^{j}$ by:
\begin{align*}
 & D_{\mathcal{L}_{1}}(\boldsymbol{x}_{r}^{j},\boldsymbol{\Theta})\leq\varepsilon_{\mathrm{lim}}\\
 & +D_{\mathcal{L}_{1}}(\boldsymbol{x}_{T}^{j},\boldsymbol{\Theta})\!\!\!\!\!\!\prod_{s=0}^{\left\lfloor \frac{r-T}{n_{\mathrm{rec}}-1}\right\rfloor -1}\!\!\!\!\left(1-n_{\mathrm{rec}}\left(\prod_{q=T+s(n_{\mathrm{rec}}-1)}^{T+(s+1)(n_{\mathrm{rec}}-1)}\delta_{q}^{j}\right)\right),
\end{align*}
where $\left\lfloor \cdot\right\rfloor $ is the floor function and
$\delta_{q}^{j}=\min_{i,\thinspace\ell}{}^{+}\boldsymbol{P}_{q,\mathrm{sub}}^{j}[i,\ell]$.\end{theorem}

\begin{proof}It follows from the definition of $\tau_{1}(\boldsymbol{U}_{T,r}^{j})$
in (\ref{eq:coeff-ergodicity}) that:
\[
D_{\mathcal{L}_{1}}(\bar{\boldsymbol{x}}_{T}^{j}\boldsymbol{U}_{T,r}^{j},\bar{\boldsymbol{\Theta}}\boldsymbol{U}_{T,r}^{j})\leq\tau_{1}(\boldsymbol{U}_{T,r}^{j})\thinspace D_{\mathcal{L}_{1}}(\bar{\boldsymbol{x}}_{T}^{j},\bar{\boldsymbol{\Theta}})\thinspace.
\]
Since $\bar{\boldsymbol{x}}_{r}^{j}=\bar{\boldsymbol{x}}_{T}^{j}\boldsymbol{U}_{T,r}^{j}$
(\ref{eq:proof_overall_evolution}), we get from the triangle inequality:
\begin{align*}
D_{\mathcal{L}_{1}}(\bar{\boldsymbol{x}}_{r}^{j},\bar{\boldsymbol{\Theta}}) & \leq D_{\mathcal{L}_{1}}(\bar{\boldsymbol{x}}_{T}^{j}\boldsymbol{U}_{T,r}^{j},\bar{\boldsymbol{\Theta}}\boldsymbol{U}_{T,r}^{j})+D_{\mathcal{L}_{1}}(\bar{\boldsymbol{\Theta}}\boldsymbol{U}_{T,r}^{j},\bar{\boldsymbol{\Theta}})\\
 & \leq\tau_{1}(\boldsymbol{U}_{T,r}^{j})D_{\mathcal{L}_{1}}(\bar{\boldsymbol{x}}_{T}^{j},\bar{\boldsymbol{\Theta}})+\varepsilon_{\mathrm{lim}}\thinspace.
\end{align*}
The sub-multiplicative property of $\tau_{1}(\boldsymbol{U}_{T,r}^{j})$
\cite[Lemma 4.3, pp. 139]{Ref:Seneta06} gives:
\[
\tau_{1}(\boldsymbol{U}_{T,r}^{j})\leq\prod_{s=0}^{\left\lfloor \frac{r-T}{n_{\mathrm{rec}}-1}\right\rfloor -1}\tau_{1}\left(\boldsymbol{U}_{T+s(n_{\mathrm{rec}}-1),T+(s+1)(n_{\mathrm{rec}}-1)}^{j}\right)\thinspace.
\]
Here, if $r>T+\left\lfloor \frac{r-T}{n_{\mathrm{rec}}-1}\right\rfloor (n_{\mathrm{rec}}-1)$,
then we neglect the contribution of the residual term by assuming
$\tau_{1}\left(\boldsymbol{U}_{T+\left\lfloor \frac{r-T}{n_{\mathrm{rec}}-1}\right\rfloor (n_{\mathrm{rec}}-1),r}^{j}\right)=1$. 

The matrix $\boldsymbol{U}_{k,k+n_{\textrm{rec}}-1}^{j}$, for any
$k\geq T$, is a positive matrix because there exists a path of length
smaller than $(n_{\mathrm{rec}}-1)$ between every two recurrent bins
(see Theorem~\ref{thm:positive-matrix} in Appendix). A conservative
lower bound on the elements in the positive matrix $\boldsymbol{U}_{T+s(n_{\mathrm{rec}}-1),T+(s+1)(n_{\mathrm{rec}}-1)}^{j}$
is given by the product of the smallest positive elements in all the
matrices, i.e., $\boldsymbol{U}_{T+s(n_{\mathrm{rec}}-1),T+(s+1)(n_{\mathrm{rec}}-1)}^{j}[i,\ell]\geq\left(\prod_{q=T+s(n_{\mathrm{rec}}-1)}^{T+(s+1)(n_{\mathrm{rec}}-1)}\delta_{q}^{j}\right)$
for all $i,\ell\in\{1,\ldots,n_{\textrm{rec}}\}$. Therefore, it follows
from (\ref{eq:coeff-ergodicity}) that $\tau_{1}\left(\boldsymbol{U}_{T+s(n_{\mathrm{rec}}-1),T+(s+1)(n_{\mathrm{rec}}-1)}^{j}\right)\leq1-n_{\textrm{rec}}\left(\prod_{q=T+s(n_{\mathrm{rec}}-1)}^{T+(s+1)(n_{\mathrm{rec}}-1)}\delta_{q}^{j}\right)<1$.
\hfill$\blacksquare$ \end{proof}

We now focus on the convergence of the swarm distribution to the desired
formation. In practical scenarios, the number of agents is finite,
hence the following theorem gives a lower bound on the number of agents.

\begin{theorem} \label{thm:Convergence-error} Let $\varepsilon_{\mathrm{lim}}>0$,
$\varepsilon_{\mathrm{bin}}>0$, and $\varepsilon_{\mathrm{conv}}>0$
represent convergence error thresholds. Let $\kappa$ represent the
latest time instant when an agent is added to or removed from the
swarm, i.e., the number of agents $m_{k}=m_{\kappa}$ for all $k\geq\kappa$.
Since $\lim_{k\rightarrow\infty}\boldsymbol{U}_{\kappa+T,k}^{j}=\boldsymbol{1}\bar{\boldsymbol{\Theta}}$
for all agents, there exists $k_{\epsilon,\mathrm{lim}}\in\mathbb{N}$
such that $D_{\mathcal{L}_{1}}(\bar{\boldsymbol{\Theta}}\boldsymbol{U}_{\kappa+T,k}^{j},\bar{\boldsymbol{\Theta}})\leq\varepsilon_{\mathrm{lim}}$
for all $k\geq k_{\epsilon,\mathrm{lim}}$ and $j\in\{1,\ldots,m_{\kappa}\}$. 

The convergence error between the swarm distribution $\boldsymbol{\mathcal{\mu}}_{k}^{\star}$
and the desired formation $\boldsymbol{\Theta}$ is probabilistically
bounded for all $k\geq k_{\epsilon,\mathrm{lim}}$ by:
\begin{align}
\mathbb{P}\left(D_{\mathcal{L}_{1}}(\boldsymbol{\mathcal{\mu}}_{k}^{\star},\boldsymbol{\Theta})\geq\varepsilon_{\mathrm{bin}}+\sigma_{k}\right) & \leq\frac{n_{\mathrm{rec}}}{4m_{\kappa}\varepsilon_{\mathrm{bin}}^{2}}\thinspace,\label{eq:L1_distance_swarm}\\
\mathbb{P}\left(D_{H}(\boldsymbol{\mathcal{\mu}}_{k}^{\star},\boldsymbol{\Theta})\geq\frac{1}{\sqrt{2}}\sqrt{\varepsilon_{\mathrm{bin}}+\sigma_{k}}\right) & \leq\frac{n_{\mathrm{rec}}}{4m_{\kappa}\varepsilon_{\mathrm{bin}}^{2}}\thinspace,\label{eq:H_distance_swarm}
\end{align}
where $\delta_{q}=\min_{j\in\{1,\ldots,m_{\kappa}\}}\delta_{q}^{j}$
and $\sigma_{k}=\varepsilon_{\mathrm{lim}}+2\prod_{s=0}^{\left\lfloor \frac{k-\kappa-T}{n_{\mathrm{rec}}-1}\right\rfloor -1}\left(1-n_{\mathrm{rec}}\left(\prod_{q=T+s(n_{\mathrm{rec}}-1)}^{T+(s+1)(n_{\mathrm{rec}}-1)}\delta_{q}\right)\right)$. 

If the number of agents satisfies the inequality:
\begin{equation}
m_{\kappa}\geq\frac{n_{\mathrm{rec}}}{16\xi_{\mathrm{des}}^{4}\varepsilon_{\mathrm{conv}}}\thinspace,
\end{equation}
then the HD between the final swarm distribution and the desired formation
is probabilistically bounded by $\varepsilon_{\mathrm{conv}}$, i.e.,

\begin{equation}
\mathbb{P}\left(D_{H}\Bigl(\lim_{k\rightarrow\infty}\boldsymbol{\mathcal{\mu}}_{k}^{\star},\boldsymbol{\Theta}\Bigr)\geq\xi_{\mathrm{des}}\right)\leq\varepsilon_{\mathrm{conv}}\thinspace,
\end{equation}
where $\xi_{\mathrm{des}}$ is the desired convergence error defined
in Definition~\ref{def:Feedback-gain}. \end{theorem} 

\begin{proof} Let $X_{k,i}^{j}$ denote the Bernoulli random variable,
where $X_{k,i}^{j}=1$ represents the event that the $j^{\textrm{th}}$
agent is actually located in bin $B[i]$ at the $k^{\textrm{th}}$
time instant (i.e., $\boldsymbol{r}_{k}^{j}[i]=1$) and $X_{k,i}^{j}=0$
otherwise (i.e., $\boldsymbol{r}_{k}^{j}[i]=0$). We get from (\ref{eq:def_prob_vector})
that $\mathbb{P}(X_{k,i}^{j}=1)=\boldsymbol{x}_{k}^{j}[i]$. Therefore
$\mathbb{E}[X_{k,i}^{j}]=\boldsymbol{x}_{k}^{j}[i]$ and $\textrm{Var}(X_{k,i}^{j})=\mathbb{E}[X_{k,i}^{j}](1-\mathbb{E}[X_{k,i}^{j}])\leq\frac{1}{4}$,
where $\mathbb{E}[\cdot]$ and $\textrm{Var}(\cdot)$ respectively
denote the expected value and the variance of the random variable.
It follows from Theorem~\ref{thm:DL1_x_Theta} that $D_{\mathcal{L}_{1}}(\boldsymbol{x}_{k}^{j},\boldsymbol{\Theta})\leq\sigma_{k}$
for all $k\geq k_{\epsilon,\mathrm{lim}}$. Therefore $|\mathbb{E}[X_{k,i}^{j}]-\boldsymbol{\Theta}[i]|\leq\sigma_{k}$
for all $k\geq k_{\epsilon,\mathrm{lim}}$. 

The swarm distribution in bin $B[i]$ at the $k^{\textrm{th}}$ time
instant is given by $\boldsymbol{\mu}_{k}^{\star}[i]=\frac{1}{m_{\kappa}}\sum_{j=1}^{m_{\kappa}}X_{k,i}^{j}$.
Therefore, $\mathbb{E}[\boldsymbol{\mu}_{k}^{\star}[i]]=\frac{1}{m_{\kappa}}\sum_{j=1}^{m_{\kappa}}\mathbb{E}[X_{k,i}^{j}]$.
The random variables $X_{k,i}^{j},\thinspace j\in\{1,\ldots,m_{\kappa}\}$
are negatively correlated because:
\begin{align*}
 & \textrm{Cov}(X_{k,i}^{j1},X_{k,i}^{j2})=\mathbb{E}[X_{k,i}^{j1}X_{k,i}^{j2}]-\mathbb{E}[X_{k,i}^{j1}]\mathbb{E}[X_{k,i}^{j2}]\\
 & =\frac{\ensuremath{\left(\begin{smallmatrix}m_{k}-2\\
n_{k,i}-2
\end{smallmatrix}\right)}}{\ensuremath{\left(\begin{smallmatrix}m_{k}\\
n_{k,i}
\end{smallmatrix}\right)}}-\left(\frac{\ensuremath{\left(\begin{smallmatrix}m_{k}-1\\
n_{k,i}-1
\end{smallmatrix}\right)}}{\ensuremath{\left(\begin{smallmatrix}m_{k}\\
n_{k,i}
\end{smallmatrix}\right)}}\right)^{2}\leq0\thinspace,
\end{align*}
where $n_{k,i}$ is the number of agents in bin $B[i]$ at the $k^{\textrm{th}}$
time instant and {\tiny{}$\left(\begin{array}{c}
\cdot\\
\cdot
\end{array}\right)$} represents the Binomial coefficient. Therefore, we get:
\begin{align*}
\textrm{Var}(\boldsymbol{\mu}_{k}^{\star}[i]) & =\frac{1}{m_{\kappa}^{2}}\left(\sum_{j=1}^{m_{\kappa}}\textrm{Var}(X_{k,i}^{j})+2\!\!\!\!\!\!\!\!\!\!\sum_{1\leq j1<j2\leq m_{\kappa}}\!\!\!\!\!\!\!\!\!\!\textrm{Cov}(X_{k,i}^{j1},X_{k,i}^{j2})\right)\\
 & \leq\frac{1}{m_{\kappa}^{2}}\sum_{j=1}^{m_{\kappa}}\textrm{Var}(X_{k,i}^{j})\leq\frac{1}{4m_{\kappa}}\thinspace.
\end{align*}
It follows from Chebychev's inequality (cf. \cite[Theorem 1.6.4, pp. 25]{Ref:Durrett05})
that for any $\varepsilon_{\mathrm{bin}}$, the pointwise error probability
for each bin is bounded by:
\[
\mathbb{P}\left(\left|\boldsymbol{\mu}_{k}^{\star}[i]-\mathbb{E}[\boldsymbol{\mu}_{k}^{\star}[i]]\right|\geq\varepsilon_{\mathrm{bin}}\right)\leq\frac{1}{4m_{\kappa}\varepsilon_{\mathrm{bin}}^{2}}\thinspace.
\]
It follows from the triangle inequality that $\left|\boldsymbol{\mu}_{k}^{\star}[i]-\mathbb{E}[\boldsymbol{\mu}_{k}^{\star}[i]]\right|\geq\left|\boldsymbol{\mu}_{k}^{\star}[i]-\boldsymbol{\Theta}[i]\right|-\sigma_{k}$,
therefore:
\begin{align*}
\mathbb{P}\left(\left|\boldsymbol{\mu}_{k}^{\star}[i]-\boldsymbol{\Theta}[i]\right|\geq\varepsilon_{\mathrm{bin}}+\sigma_{k}\right) & \leq\frac{1}{4m_{\kappa}\varepsilon_{\mathrm{bin}}^{2}}\thinspace.
\end{align*}
The bound on $\mathcal{L}_{1}$ distance is obtained using Boole's
inequality:
\begin{align*}
 & \mathbb{P}\left(D_{\mathcal{L}_{1}}(\boldsymbol{\mathcal{\mu}}_{k}^{\star},\boldsymbol{\Theta})\geq\varepsilon_{\mathrm{bin}}+\sigma_{k}\right)\\
 & \leq\sum_{i=1}^{n_{\textrm{rec}}}\mathbb{P}\left(\left|\boldsymbol{\mu}_{k}^{\star}[i]-\boldsymbol{\Theta}[i]\right|\geq\varepsilon_{\mathrm{bin}}+\sigma_{k}\right)\leq\frac{n_{\textrm{rec}}}{4m_{\kappa}\varepsilon_{\mathrm{bin}}^{2}}\thinspace.
\end{align*}
The bound on HD follows from (\ref{eq:HD_L1_compare}). 

Since $\lim_{k\rightarrow\infty}\boldsymbol{x}_{k}^{j}=\boldsymbol{\Theta}$,
we get $\lim_{k\rightarrow\infty}\mathbb{E}[X_{k,i}^{j}]=\boldsymbol{\Theta}[i]$
and $\lim_{k\rightarrow\infty}\sigma_{k}=0$. Setting $\varepsilon_{\mathrm{bin}}=2\xi_{\textrm{des}}^{2}$,
we get:
\[
\mathbb{P}\left(D_{H}\Bigl(\lim_{k\rightarrow\infty}\boldsymbol{\mathcal{\mu}}_{k}^{\star},\boldsymbol{\Theta}\Bigr)\geq\xi_{\textrm{des}}\right)\leq\frac{n_{\textrm{rec}}}{16m_{\kappa}\xi_{\textrm{des}}^{4}}\thinspace.
\]
The lower bound on the number of agents is given by $\frac{n_{\textrm{rec}}}{16m_{\kappa}\xi_{\textrm{des}}^{4}}\leq\varepsilon_{\mathrm{conv}}$.
\hfill$\blacksquare$ \end{proof}

\begin{remark} It follows from Theorem~\ref{thm:Convergence-error}
and the weak law of large numbers \cite[pp. 86]{Ref:Billingsley95}
that the final swarm distribution $\lim_{k\rightarrow\infty}\boldsymbol{\mathcal{\mu}}_{k}^{\star}$
converges in probability to the desired formation $\boldsymbol{\Theta}$
as the number of agents $m_{\kappa}$ tends to infinity.\hfill$\Box$
\end{remark}

Thus, we have proved the convergence of the PSG\textendash IMC algorithm
for shape formation. We now discuss its property of robustness and
some extensions.

\begin{remark} \label{rem:Robustness} \textit{(Robustness of the
PSG\textendash IMC Algorithm)} The PSG\textendash IMC algorithm satisfies
the Markov (memoryless) property because the action of each agent
depends only on its present bin location and the current swarm distribution.
This property ensures that all the agents re-start their guidance
trajectory from their present bin location during every time instant.
Thus, the swarm continues to converge to the desired shape even if
agents are added to or removed from the swarm, or if some agents have
not reached their target bin during the previous time instant. 

Moreover, the PSG\textendash IMC algorithm can tolerate estimation
errors $\epsilon_{\mathrm{est}}$ (\ref{eq:estimation_error}) in
the feedback of the current swarm distribution $\boldsymbol{\mathcal{\mu}}_{k}^{j}$.
The distance between the feedback gains $\xi_{k}^{j}=D_{H}(\boldsymbol{\Theta},\boldsymbol{\mathcal{\mu}}_{k}^{j})$
and $\xi_{k}^{\star}=D_{H}(\boldsymbol{\Theta},\boldsymbol{\mathcal{\mu}}_{k}^{\star})$
is bounded by \cite{Ref:Pollard00asymptopia}:
\begin{equation}
\left|\xi_{k}^{j}-\xi_{k}^{\star}\right|\leq D_{H}(\boldsymbol{\mathcal{\mu}}_{k}^{\star},\boldsymbol{\mathcal{\mu}}_{k}^{j})\leq\frac{1}{\sqrt{2}}\epsilon_{\mathrm{est}}^{\frac{1}{2}}\thinspace.
\end{equation}
Even though $\xi_{k}^{j}$ might differ from $\xi_{k}^{\star}$ substantially,
the resulting Markov matrix $\boldsymbol{M}_{k}^{j}$ still has $\boldsymbol{\Theta}$
as it stationary distribution. Therefore the agent's PMF vector $\boldsymbol{x}_{k}^{j}$
still converges to $\boldsymbol{\Theta}$, and consequently the swarm
distribution also converges to $\boldsymbol{\Theta}$.\hfill$\Box$
\end{remark}

\begin{remark} \label{rem:Collision-avoidance} \textit{(Collision
Avoidance in the PSG\textendash IMC Algorithm) }The PSG\textendash IMC
algorithm can handle inter-agent collision avoidance using line 15
in \textbf{Method~2}. We have not explicitly considered collision
avoidance with stationary obstacles because they can be easily handled
by the existing framework. If the stationary obstacles are comparable
or larger than the bin size, then the bins are designed such that
they do not overlap with these obstacles. Furthermore, the motion
constraints matrices are designed to prevent transitions that are
not allowed due to these obstacles. If the stationary obstacles are
significantly smaller than the bins, then they can be handled by the
lower-level collision avoidance algorithm as shown in Remark~\ref{rem:collision-free-motion}
in Appendix. \hfill$\Box$ \end{remark}

\begin{figure}[h]
\begin{centering}
\includegraphics[bb=0bp 10bp 640bp 205bp,clip,width=2in]{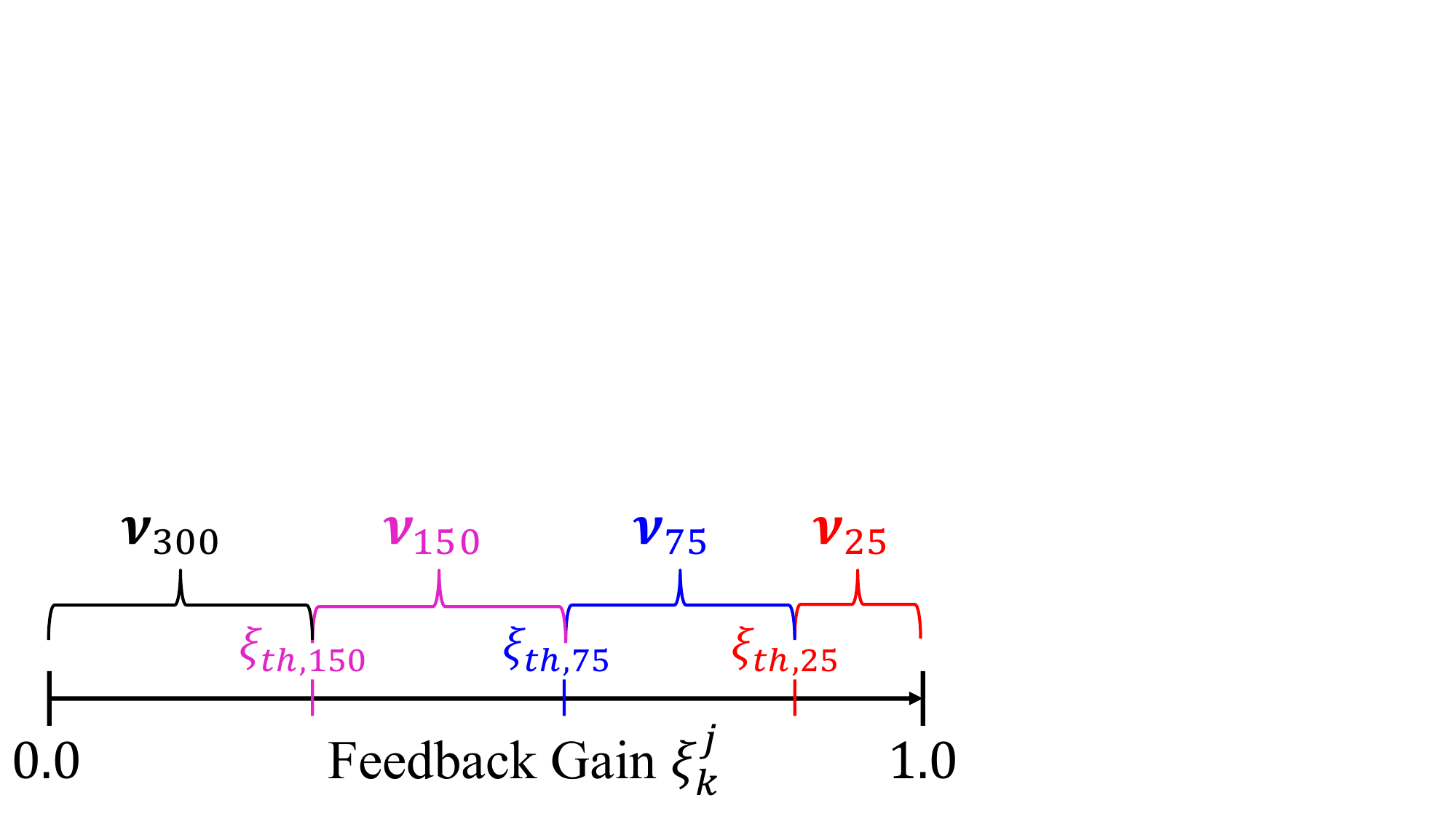}
\par\end{centering}
\caption{This image shows the different thresholds and the corresponding resolution
of the desired formation that should be used when $\xi_{k}^{j}$ is
within those thresholds. \label{fig:multiresolution-HD-thresholds}}
\end{figure}

\begin{remark} \label{rem:Multiresolution-PSG-IMC}\textit{(Multiresolution
PSG\textendash IMC Algorithm for Shape Formation) }We can take advantage
of multiresolution representation of the desired formation in our
guidance strategy (see Fig.~\ref{fig:Multi-resolution-Images-of-Eiffel-Tower}).
The key idea is that the agents use an appropriate resolution of the
desired formation depending on the feedback gain. For example, as
shown in Fig.~\ref{fig:multiresolution-HD-thresholds}, we select
thresholds $\xi_{\textrm{th},150}$, $\xi_{\textrm{th},75}$, and
$\xi_{\textrm{th},25}$ so that the agents use the appropriate resolution
of the desired formation ($\boldsymbol{\nu}_{300}$, $\boldsymbol{\nu}_{150}$,
$\boldsymbol{\nu}_{75}$ or $\boldsymbol{\nu}_{25}$ in Fig.~\ref{fig:Multi-resolution-Images-of-Eiffel-Tower})
if the feedback gain $\xi_{k}^{j}$ is within these thresholds. The
main advantage of this approach is its computational efficiency.\hfill$\Box$
\end{remark}

\begin{remark} \label{rem:Time-varying-physical-space} \textit{(Time-Varying
Physical Space of the Swarm)} The compact physical space over which
the swarm is distributed need not be time-invariant in the global
reference frame. The local reference frame of the swarm can follow
a predefined trajectory in the global reference frame (e.g., an orbit
in space or a trajectory in the sea) and the time-varying position
of each bin can be computed from this known trajectory. Consequently,
all the algorithms discussed in this paper are also applicable in
this scenario. \hfill$\Box$ \end{remark}

\section{PSG\textendash IMC Algorithm for Area Exploration \label{sec:Goal-Searching}}

In this section, we present an extension of the PSG\textendash IMC
algorithm for area exploration in which a swarm of distributed agents
are driven to match the unknown target distribution of some physical
or artificial phenomena (e.g., oil spill). This problem is commonly
called \textit{goal searching} \cite{Ref:Hereford10}.

\begin{definition} \label{def:Target-distribution} \textit{(Unknown
Target Distribution $\boldsymbol{\Omega}$)} The unknown target distribution
$\boldsymbol{\Omega}$ is a probability (row) vector in $\mathbb{R}^{n_{\textrm{bin}}}$,
where each element $\boldsymbol{\Omega}[i]$ represents the target
distribution in the corresponding bin $B[i]$. Each agent can measure
the target distribution in its present bin. \hfill$\Box$ \end{definition}

Each agent independently determines its bin-to-bin trajectory using
the PSG\textendash IMC algorithm for area exploration so that the
overall swarm converges to this unknown target distribution $\boldsymbol{\Omega}$.
The key idea of this algorithm is that the waiting time in a bin is
directly proportional to the target distribution in that bin. 

\begin{algorithm}[h]
\begin{tabular}{cl}
\hline 
\multicolumn{2}{l}{\textbf{Method 3}: PSG\textendash IMC Algorithm for Area Exploration}\tabularnewline
\hline 
{\small{}1:} & {\small{}Lines 1\textendash 2 in}\textbf{\small{} Method~2}{\small{}
and given $\tau_{c}$, $\xi^{j}$}\tabularnewline
{\small{}2:} & {\small{}Measure target distribution in present bin $\boldsymbol{\Omega}[i]$ }\tabularnewline
{\small{}3:} & \textbf{\small{}if}{\small{} $k-k_{0}<\tau_{c}\boldsymbol{\Omega}[i]$,
}\textbf{\small{}then}{\small{} }\tabularnewline
{\small{}4:} & {\small{}$\qquad$Wait in bin $B[i]$ }\tabularnewline
{\small{}5:} & \textbf{\small{}else }{\small{}Set $\boldsymbol{\Theta}=\frac{\boldsymbol{1}^{T}}{n_{\mathrm{bin}}}$
and $\xi_{k}^{j}=\xi^{j}$ }\tabularnewline
{\small{}6:} & \textbf{\small{}$\qquad$}{\small{}Compute the term $\eta_{k,i}^{j}$
using (\ref{eq:eta_ki-goal}) }\tabularnewline
{\small{}7:} & \textbf{\small{}$\qquad$}{\small{}Lines 9, 11\textendash 13, 15 in}\textbf{\small{}
Method~2}{\small{}, Set $k_{0}=k$ }\tabularnewline
{\small{}8:} & \textbf{\small{}end if}{\small{} }\tabularnewline
\hline 
\end{tabular}
\end{algorithm}

The pseudo-code for this PSG\textendash IMC algorithm for area exploration
is given in \textbf{Method~3}. The $j^{\textrm{th}}$ agent first
measures the target distribution in its present bin $\boldsymbol{\Omega}[i]$
(line 2). The waiting time in bin $B[i]$ is greater than or equal
to $\tau_{c}\boldsymbol{\Omega}[i]$, where $\tau_{c}$ is the constant
of proportionality.

\begin{assumption} \label{assump:tau_c} \textit{(Waiting Time Constant
$\tau_{c}$)} Let $\Delta$ denote the time step of the algorithm.
Each agent has a-priori information about the common waiting time
constant $\tau_{c}$, which is selected such that $\tau_{c}>\frac{\Delta}{\min^{+}\boldsymbol{\Omega}}$.
The actual value of $\tau_{c}$ plays a crucial role in the convergence
analysis. \hfill$\Box$ \end{assumption}

The agent checks if it has spent enough time instants in bin $B[i]$
(line 3), where $k_{0}$ is is set in line 7. When the algorithm starts
($k=1$), we set $k_{0}=1$. If the agent has not spent enough time
instants in bin $B[i]$, then it continues to wait in bin $B[i]$
(line 4). 

If the agent has spent enough time instants in bin $B[i]$, then it
sets $\boldsymbol{\Theta}=\frac{\boldsymbol{1}^{T}}{n_{\mathrm{bin}}}$
because it wants to uniformly explore all the bins (line 5). The agent
sets the feedback gain $\xi_{k}^{j}$ to some positive known constant
$\xi^{j}\in(0,1)$ because the it does not know the target distribution
$\boldsymbol{\Omega}$ (line 5). In order to suppress undesirable
transitions from bins that are deficient in agents, the agent computes
the term $\eta_{k,i}^{j}$ as follows (line 6):
\begin{align}
\eta_{k,i}^{j} & =\exp(-\tau^{j}k)\frac{\exp\left(\beta^{j}(\boldsymbol{\Omega}[i]-\boldsymbol{\mathcal{\mu}}_{k}^{j}[i])\right)}{\exp\left(\beta^{j}|\boldsymbol{\Omega}[i]-\boldsymbol{\mathcal{\mu}}_{k}^{j}[i]|\right)}\thinspace.\label{eq:eta_ki-goal}
\end{align}
Then, the agent computes the transition probabilities $\boldsymbol{P}_{k}^{j}[i,1\!:\!n_{\textrm{bin}}]$
using lines 9 and 11 from \textbf{Method~2}, selects the next bin
using lines 12\textendash 13 from \textbf{Method~2}, and goes to
the selected bin using line 15 from \textbf{Method~2} (line 7). Finally,
the agent sets $k_{0}$ equal to the current time instant $k$ (line
7). We now discuss the convergence analysis of this algorithm.

\begin{theorem} \label{thm:Convergence-goal-searching} \textit{(Convergence
of IMC for Area Exploration)} Let $\hat{\boldsymbol{\Omega}}\in\mathbb{R}^{n_{\textrm{bin}}}$
be a probability (row) vector. Let $n_{\Omega}$ denote the number
of bins with nonzero elements in $\boldsymbol{\Omega}$. According
to \textbf{Method~3}, the PMF vector $\boldsymbol{x}_{k}^{j}$ converges
pointwise to the distribution $\hat{\boldsymbol{\Omega}}$ irrespective
of the initial condition, where
\[
D_{\mathcal{L}_{1}}(\hat{\boldsymbol{\Omega}},\boldsymbol{\Omega})\leq\frac{n_{\Omega}\frac{\Delta}{\tau_{c}}}{(n_{\mathrm{bin}}-n_{\Omega})\frac{\Delta}{\tau_{c}}+1}\thinspace.
\]
If $\tau_{c}\gg\Delta$, then $\lim_{k\rightarrow\infty}\boldsymbol{x}_{k}^{j}=\boldsymbol{\Omega}$
pointwise for all agents. \end{theorem} 

\begin{proof} Here, all bins are recurrent bins because $\boldsymbol{\Theta}=\frac{\boldsymbol{1}^{T}}{n_{\mathrm{bin}}}$.
It follows from Theorem~\ref{thm:Convergence-of-imhomo-MC} that
as $k\rightarrow\infty$, an agent is equally likely to transition
to any bin $B[i]$. But the waiting time in each bin, under Assumption~\ref{assump:tau_c},
is given by:
\begin{align*}
\begin{array}{c}
\textrm{Waiting time}\\
\textrm{ in bin }B[i]
\end{array} & =\begin{cases}
\Delta\left\lceil \frac{\tau_{c}\boldsymbol{\Omega}[i]}{\Delta}\right\rceil  & \textrm{ if }i\in\mathcal{V}\\
\Delta & \textrm{ if }i\in\bar{\mathcal{V}}
\end{cases}\thinspace,
\end{align*}
where $\mathcal{V}$ represents the set of all bins $B[i]$ where
$\boldsymbol{\Omega}[i]>0$, and $\bar{\mathcal{V}}$ represents the
set of all bins $B[i]$ where $\boldsymbol{\Omega}[i]=0$. Note that
$|\mathcal{V}|=n_{\Omega}$. Therefore, $\lim_{k\rightarrow\infty}\boldsymbol{x}_{k}^{j}=\hat{\boldsymbol{\Omega}}$
pointwise for all agents, where:
\begin{align*}
\hat{\boldsymbol{\Omega}}[i] & =\begin{cases}
\frac{\Delta\left\lceil \frac{\tau_{c}\boldsymbol{\Omega}[i]}{\Delta}\right\rceil }{\sum_{\ell\in\bar{\mathcal{V}}}\Delta+\sum_{\ell\in\mathcal{V}}\Delta\left\lceil \frac{\tau_{c}\boldsymbol{\Omega}[\ell]}{\Delta}\right\rceil } & \textrm{ if }i\in\mathcal{V}\\
\frac{\Delta}{\sum_{\ell\in\bar{\mathcal{V}}}\Delta+\sum_{\ell\in\mathcal{V}}\Delta\left\lceil \frac{\tau_{c}\boldsymbol{\Omega}[\ell]}{\Delta}\right\rceil } & \textrm{ if }i\in\bar{\mathcal{V}}
\end{cases}\thinspace.
\end{align*}
The $\mathcal{L}_{1}$ distance between $\hat{\boldsymbol{\Omega}}$
and $\boldsymbol{\Omega}$ is given by:
\begin{align*}
 & D_{\mathcal{L}_{1}}(\hat{\boldsymbol{\Omega}},\boldsymbol{\Omega})=\sum_{i\in\bar{\mathcal{V}}}\frac{\Delta}{\sum_{\ell\in\bar{\mathcal{V}}}\Delta+\sum_{\ell\in\mathcal{V}}\Delta\left\lceil \frac{\tau_{c}\boldsymbol{\Omega}[\ell]}{\Delta}\right\rceil }\\
 & \quad+\sum_{i\in\mathcal{V}}\left|\frac{\Delta\left\lceil \frac{\tau_{c}\boldsymbol{\Omega}[i]}{\Delta}\right\rceil }{\sum_{\ell\in\bar{\mathcal{V}}}\Delta+\sum_{\ell\in\mathcal{V}}\Delta\left\lceil \frac{\tau_{c}\boldsymbol{\Omega}[\ell]}{\Delta}\right\rceil }-\boldsymbol{\Omega}[i]\right|\thinspace,\\
 & \leq\frac{\sum_{i\in\bar{\mathcal{V}}}\Delta+\sum_{i\in\mathcal{V}}\left(\Delta-\boldsymbol{\Omega}[i](n_{\mathrm{bin}}-n_{\Omega})\Delta\right)}{(n_{\mathrm{bin}}-n_{\Omega})\Delta+\tau_{c}}\thinspace,\\
 & \leq\frac{n_{\Omega}\Delta}{(n_{\mathrm{bin}}-n_{\Omega})\Delta+\tau_{c}}\thinspace.
\end{align*}
If $\tau_{c}\gg\Delta$, then $D_{\mathcal{L}_{1}}(\hat{\boldsymbol{\Omega}},\boldsymbol{\Omega})=0$
and $\lim_{k\rightarrow\infty}\boldsymbol{x}_{k}^{j}=\boldsymbol{\Omega}$
pointwise for all agents. \hfill$\blacksquare$\end{proof} 

The remaining convergence analysis straightforwardly follows that
of the previous algorithm given in Section \ref{sec:Convergence-Analysis}.

\begin{figure*}[t]
\begin{centering}
\begin{tabular}{cc}
\includegraphics[height=2.5in]{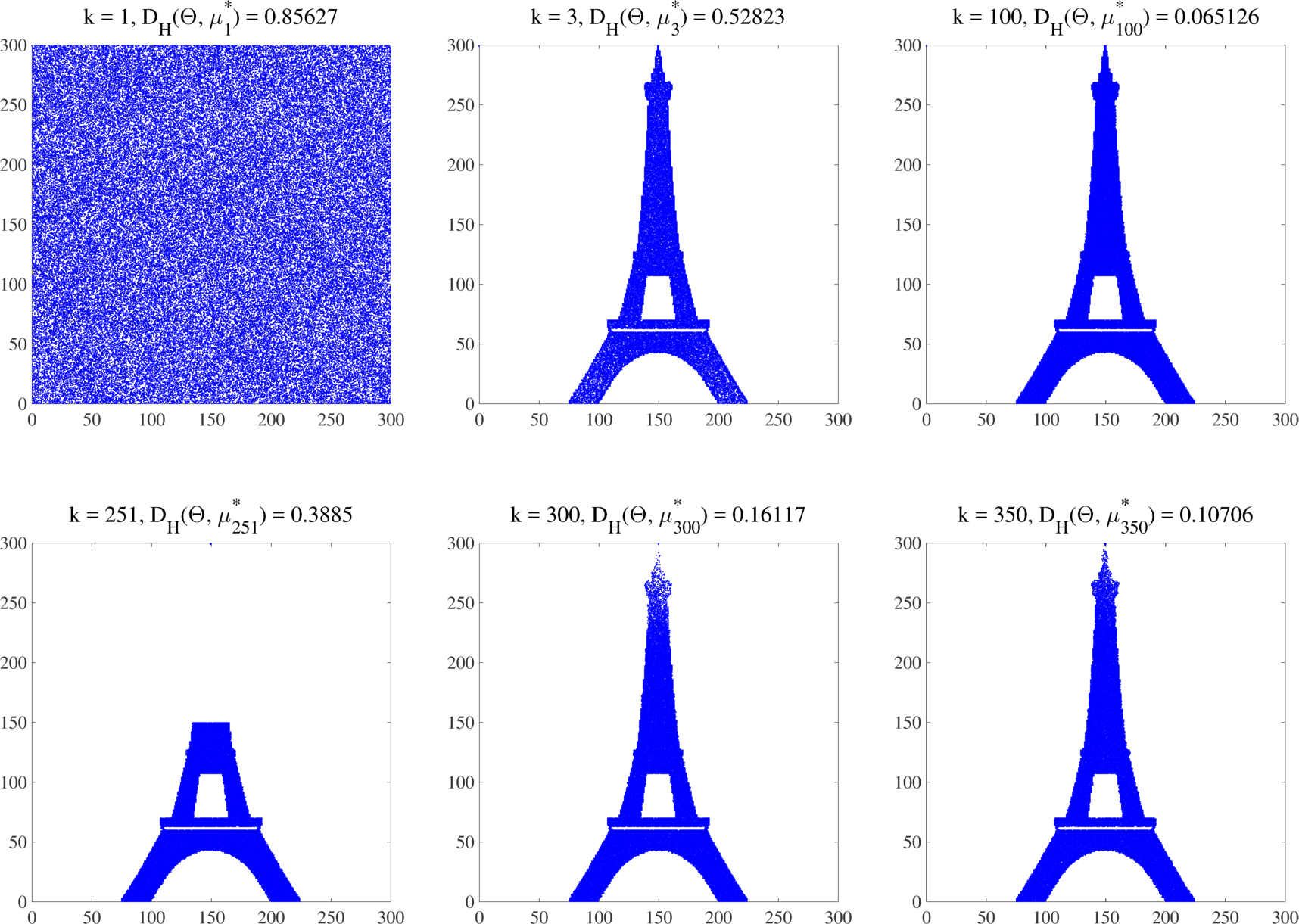} & \includegraphics[height=2.5in]{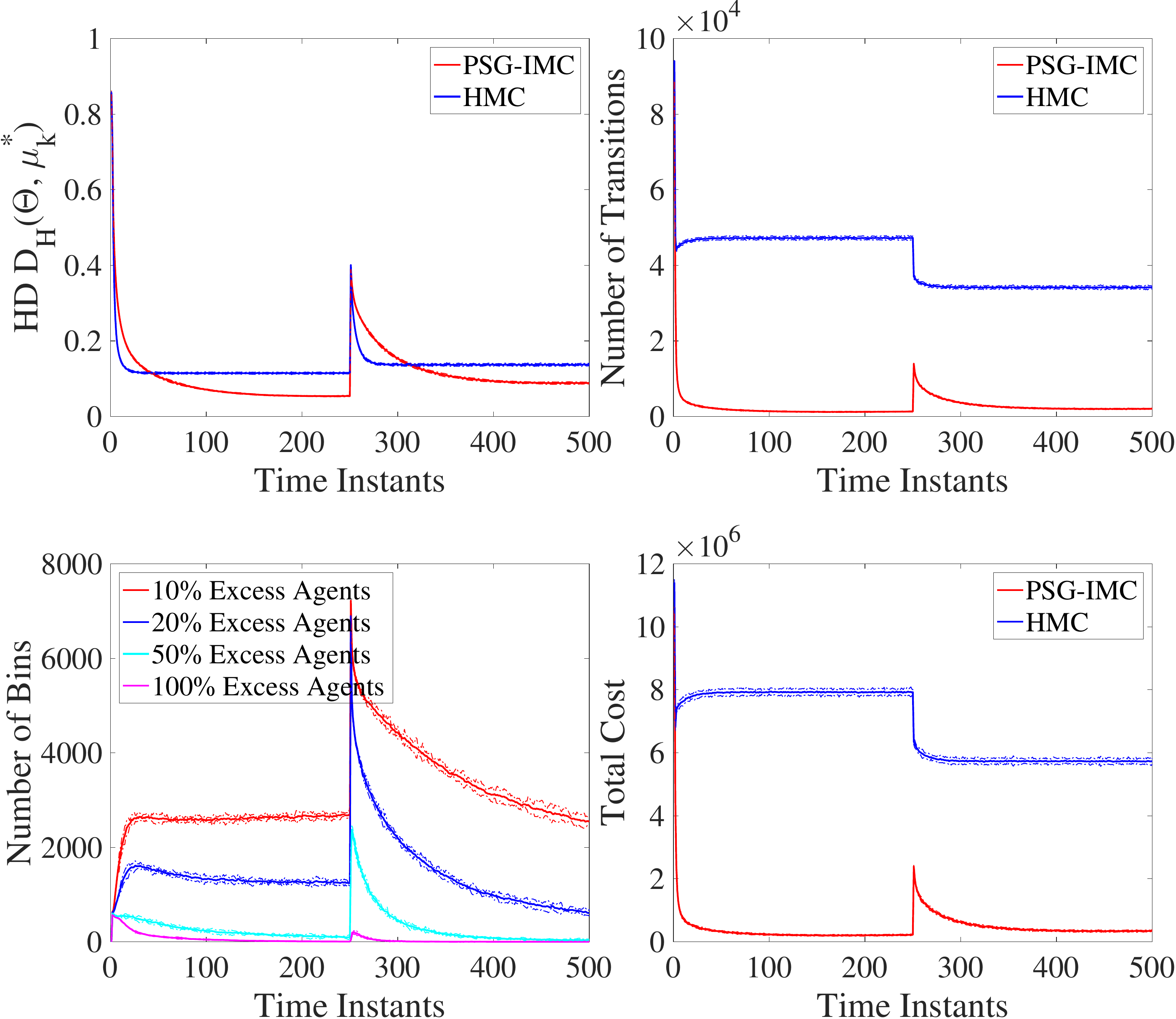}\tabularnewline
(a) & (b)\tabularnewline
\end{tabular}
\par\end{centering}
\caption{\textbf{(a)} This plot shows the swarm distribution at different time
instants, in a sample run of the Monte Carlo simulation. Starting
from a uniform distribution, the swarm converges to the desired formation
of the Eiffel Tower. At the $251^{\textrm{st}}$ time instant, the
agents in the top half of the formation are removed and the remaining
agents reconfigure into the desired formation. See the supplementary
video (SV2). \textbf{(b)} The cumulative results of $10$ Monte Carlo
simulations are shown. The discontinuity at the $251^{\textrm{st}}$
time instant is because of the removal of agents from the top half
of the formation. \label{fig:Eiffel-Tower-pattern-formation-MC}}
\end{figure*}

\section{Numerical Simulations and Experiments \label{sec:Numerical-Simulation-pattern-formation}}

In this section, we present results of numerical simulation and experimentation
for the shape formation algorithm in Sections~\ref{subsec:Numerical-Simulations-fine-resolution}\textendash \ref{subsec:Experimental-Validation}
and results of numerical simulations for the area exploration algorithm
in Section~\ref{subsec:Numerical-Simulation-goal searching}.

\subsection{Numerical Simulations for Shape Formation with Fine Spatial Resolution
\label{subsec:Numerical-Simulations-fine-resolution}}

In this subsection, we show that the PSG\textendash IMC algorithm
for shape formation can be used to accomplish multiple complex formation
shapes with fine spatial resolutions. At the start of each simulation,
a swarm of $10^{5}$ agents are uniformly distributed across the physical
space. During each time instant, each agent gets the error-free feedback
of the current swarm distribution $\boldsymbol{\mu}_{k}^{\star}$.
The cost of transition is equal to the $\ell_{1}$ distance between
bins, therefore it is symmetric. We use the following constants $\varepsilon_{M}=1$,
$\varepsilon_{C}=0.1$, $\tau^{j}=10^{-3}$, and $\beta^{j}=1.8\times10^{5}$.

In the first example, the desired formation $\boldsymbol{\Theta}$
is that of the Eiffel Tower ($\boldsymbol{\nu}_{300}$ in Fig.~\ref{fig:Multi-resolution-Images-of-Eiffel-Tower},
with $300\times300$ bins). Each agent is allowed to transition to
only those bins that are at most $50$ steps away. Starting from a
uniform distribution, the agents attain the desired formation in $100$
time instants (see Fig.~\ref{fig:Eiffel-Tower-pattern-formation-MC}(a)).
At the $251^{\textrm{st}}$ time instant, approximately $3\times10^{4}$
agents are removed from the top half of the formation and the remaining
agents reconfigure into the desired formation. 

For the HMC-based algorithm, we generate the homogeneous Markov matrix
using (\ref{eq:scalable_Markov1})-(\ref{eq:scalable_Markov2}) and
setting $\xi_{k}^{j}=1$. The cumulative results of $10$ Monte Carlo
simulations for the PSG\textendash IMC and HMC-based algorithms are
shown in Fig.~\ref{fig:Eiffel-Tower-pattern-formation-MC}(b). Compared
to the HMC-based algorithm, the PSG\textendash IMC algorithm provides
approximately $2$ times improvement in HD, $16$ times reduction
in the cumulative number of transitions in $500$ time instants, and
$16$ times reduction in the total cost incurred by all the agents
in $500$ time instants. The key reasons behind the superior performance
of the PSG\textendash IMC algorithm are as follows: \\
(i) In Fig.~\ref{fig:Eiffel-Tower-pattern-formation-MC}(b), the
HD of the HMC algorithm reaches an equilibrium at $0.115$ after approximately
$40$ time instants. The HMC algorithm allows undesirable transitions
(i.e., transitions from bins with fewer agents to bins with surplus
agents) which increases the HD. Therefore, the HD for the HMC algorithm
reaches an equilibrium because these undesirable transitions reach
an equilibrium with the other favorable transitions. Such undesirable
transitions are largely avoided in the PSG\textendash IMC algorithm
(due to lines 10\textendash 11 in \textbf{Method~2}), hence the resulting
HD after $250$ time instants is $0.055$ (i.e., $\approx2$ times
improvement compared to HMC). The final HD can be further reduced
by tuning $\tau^{j}$ and $\beta^{j}$. But such undesirable transitions
prevent both these Markovian approaches from truly achieving zero
convergence error.\\
(ii) In the HMC algorithm, there are $1.9\times10^{6}$ transitions
in the first $40$ time instants. This is significantly more than
that of the PSG\textendash IMC algorithm (i.e., $5.6\times10^{5}$
transitions in $250$ time instant). In the PSG\textendash IMC algorithm,
the number of transitions at each time instant is proportional to
the HD. This helps in achieving faster convergence (when HD is large)
while avoiding unnecessary transitions (when HD is small). This also
ensures that the agents settle down after the desired formation is
achieved. Note that the total number of transitions in the HMC algorithm
in $250$ time instant is extremely large (i.e., $1.2\times10^{7}$
transitions). \\
(iii) There are $7$\textendash $9$ agents in each recurrent bin.
For the PSG\textendash IMC algorithm, the number of bins with $1$\textendash $2$
excess agents (i.e., $10$\textendash $20\%$) is shown in Fig.~\ref{fig:Eiffel-Tower-pattern-formation-MC}(b).
The number of bins with a large number of excess agents (i.e., $50$\textendash $100\%$)
is a small fraction of the total number of bins. Hence this algorithm
also avoids traffic jams or bottlenecks. \\
Consequently, the PSG\textendash IMC algorithm achieves a smaller
convergence error than the HMC-based algorithm and significantly reduces
the number of transitions for achieving and maintaining the desired
formation. Moreover, these two key reasons depend on the feedback
and, therefore, do not hold true for HMC-based algorithms.

\begin{figure}[h]
\begin{centering}
\includegraphics[width=3in]{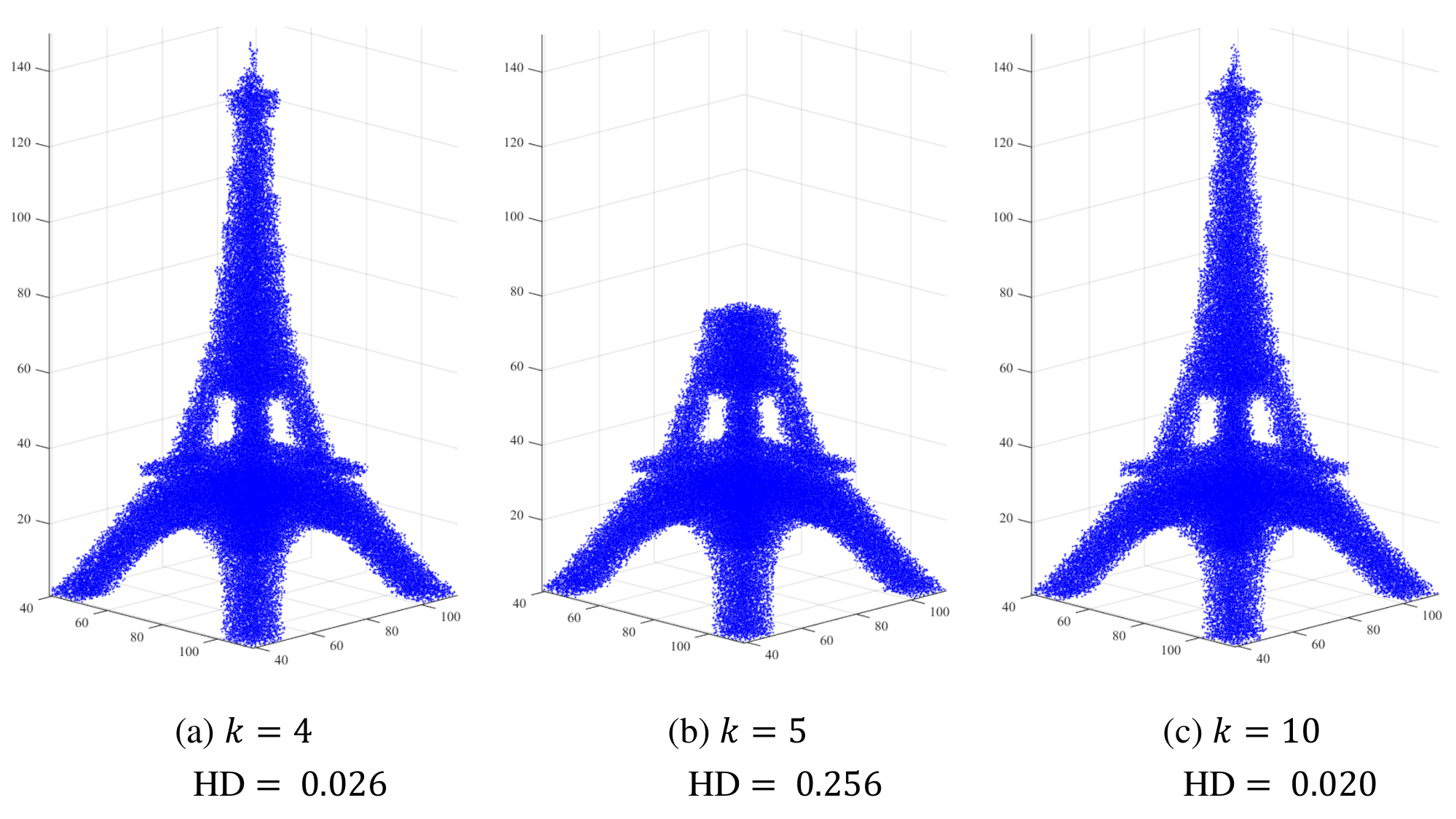}
\par\end{centering}
\caption{The swarm attains the 3-D shape of the Eiffel Tower. When agents are
removed from the top half of the formation, the remaining agents reconfigure
to the desired formation. \label{fig:Eiffel_3D}}
\end{figure}

In the next example, the desired formation $\boldsymbol{\Theta}$
is that of the 3-D Eiffel Tower (see Fig.~\ref{fig:Eiffel_3D}, with
$150\times150\times150$ bins). Starting from a uniform distribution
and no motion constraints, a swarm of $10^{5}$ agents achieve the
desired formation in a few time instants. When $1.25\times10^{4}$
agents are removed from the top half of the formation, the remaining
agents reconfigure to the desired formation in a few more time instants. 

\begin{figure}[h]
\begin{centering}
\includegraphics[bb=0bp 0bp 830bp 540bp,clip,width=3in]{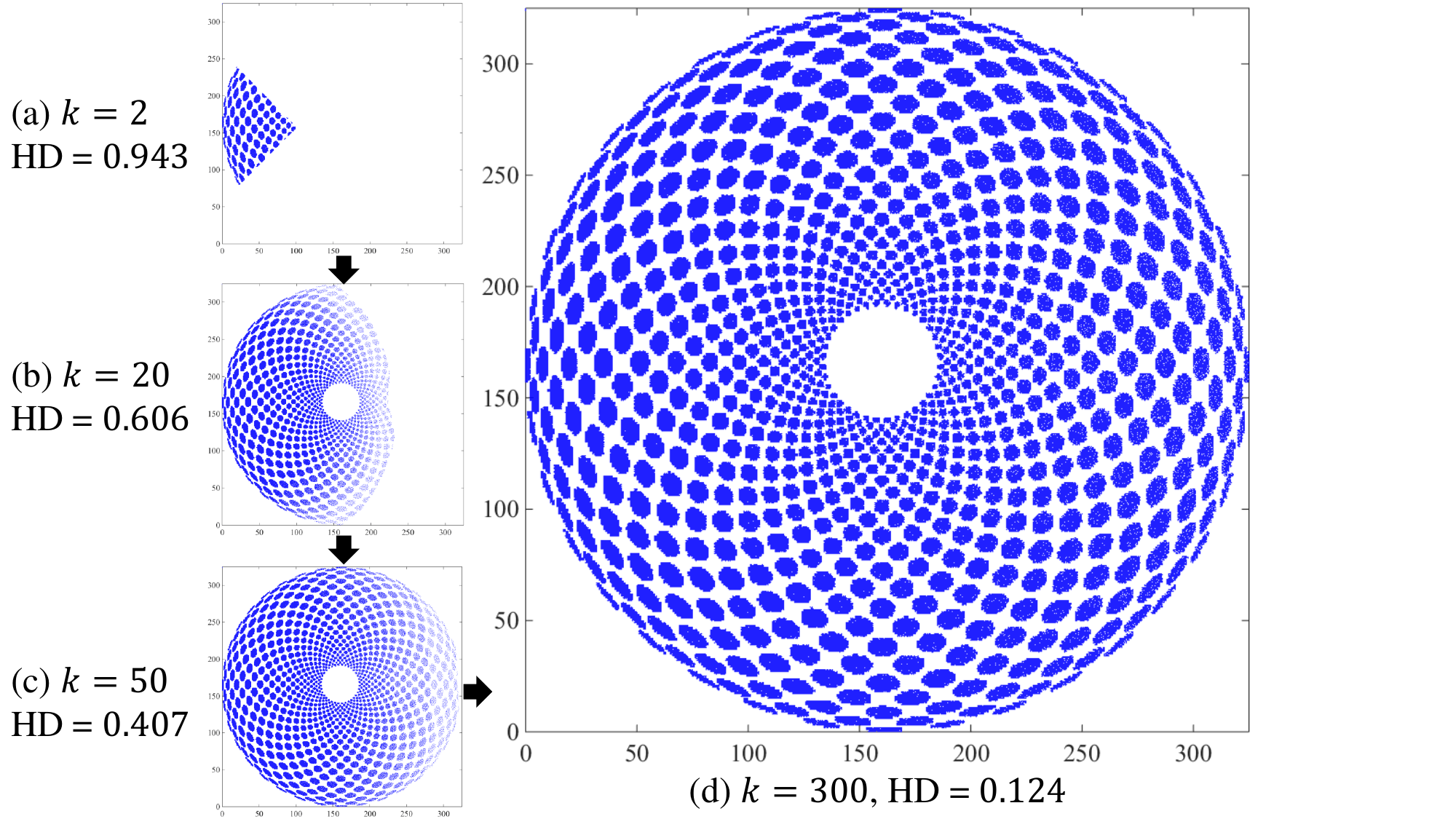}
\par\end{centering}
\caption{This plot shows the swarm distribution at different time instants,
where the swarm attains the desired formation shape with multiple
disconnected parts. See the supplementary video (SV3). \label{fig:BH_sim}}
\end{figure}

In the next example, the desired formation $\boldsymbol{\Theta}$
in Fig.~\ref{fig:BH_sim}(d), with $325\times325$ bins, has multiple
disconnected parts. Each agent is allowed to transition to only those
bins that are at most $50$ steps away. In this case, the recurrent
bins are not contiguous, but they satisfy property (3) in Definition~\ref{def:Motion-Constraints}.
A swarm of $10^{6}$ agents starts from the left-most bin (located
at $(1,163)$) and attains the desired formation in $300$ time instants
(as shown in Fig.~\ref{fig:BH_sim}(a)\textendash (d)). 

\begin{table}[h]
\caption{Computation Times for the PSG\textendash IMC and HMC-based algorithms
\label{tab:Computation-Times}}

\centering{}%
\begin{tabular}{|c|c|c|c|}
\hline 
\begin{tabular}{c}
Desired\tabularnewline
Distribution\tabularnewline
\end{tabular} & %
\begin{tabular}{c}
Number of\tabularnewline
Agents\tabularnewline
\end{tabular} & PSG\textendash IMC & HMC\tabularnewline
\hline 
\hline 
 & $5000$ & $15$ seconds & $23$ seconds\tabularnewline
\cline{2-4} 
$\boldsymbol{\nu}_{25}$ & $10^{4}$ & $32$ seconds & $43$ seconds\tabularnewline
\cline{2-4} 
($625$ bins) & $10^{5}$ & $5$ minutes & $6$ minutes\tabularnewline
\cline{2-4} 
 & $10^{6}$ & $45$ minutes & $54$ minutes\tabularnewline
\hline 
\hline 
$\boldsymbol{\nu}_{75}$ & $10^{4}$ & $3$ minutes & $3.5$ minutes\tabularnewline
\cline{2-4} 
($5625$ bins) & $10^{5}$ & $30$ minutes & $36$ minutes\tabularnewline
\cline{2-4} 
 & $10^{6}$ & $5$ hours & $5.9$ hours\tabularnewline
\hline 
\hline 
$\boldsymbol{\nu}_{150}$ & $10^{5}$ & $2.0$ hours & $2.4$ hours\tabularnewline
\cline{2-4} 
($2.25\times10^{4}$ bins) & $10^{6}$ & $23$ hours & $1$ day\tabularnewline
\hline 
\hline 
$\boldsymbol{\nu}_{300}$ & $10^{5}$ & $5.3$ hours & $8$ hours\tabularnewline
\cline{2-4} 
($9\times10^{4}$ bins) & $10^{6}$ & $2.5$ days & $3.3$ days\tabularnewline
\hline 
\end{tabular}
\end{table}

In Table~\ref{tab:Computation-Times}, the computation times using
the PSG\textendash IMC and HMC-based algorithms on a desktop computer
are shown. The simulation setup for all these runs is exactly the
same as shown in Fig.~\ref{fig:Eiffel-Tower-pattern-formation-MC}.
Although both the algorithms scale well with the spatial resolution
of the desired distribution and the number of agents in the swarm,
the PSG\textendash IMC algorithm performs better because of the smaller
number of transitions. Thus, the robustness and scalability properties
of the PSG\textendash IMC algorithm for shape formation are evident
in these simulation results.

\subsection{Numerical Simulations for Shape Formation with Coarse Spatial Resolution
and Estimation Errors \label{subsec:Numerical-Simulations-Estimation-error}}

The objective of this subsection is to study the effect of estimation
errors on the three Eulerian algorithms, namely PSG\textendash IMC,
HMC, and the Probabilistic Swarm Guidance using Optimal Transport
(PSG\textendash OT) algorithm \cite{Ref:Bandyopadhyay14MSC} (see
Remark~\ref{rem:PSG-OT-algo} in Appendix). The desired formation
$\boldsymbol{\Theta}$ is given by the coarse Eiffel Tower ($\boldsymbol{\nu}_{25}$
in Fig.~\ref{fig:Multi-resolution-Images-of-Eiffel-Tower}, with
$25\times25$ bins) because the PSG\textendash OT's LP (\ref{eq:OTP_LP})
cannot be solved for finer resolutions. The simulation setup is similar
to that in Section~\ref{subsec:Numerical-Simulations-fine-resolution}.
A swarm of 5000 agents is used and each agent is allowed to transition
to only those bins which are at most $9$ steps away. 

\begin{figure}[h]
\begin{centering}
\begin{tabular}{cc}
\includegraphics[height=1.8in]{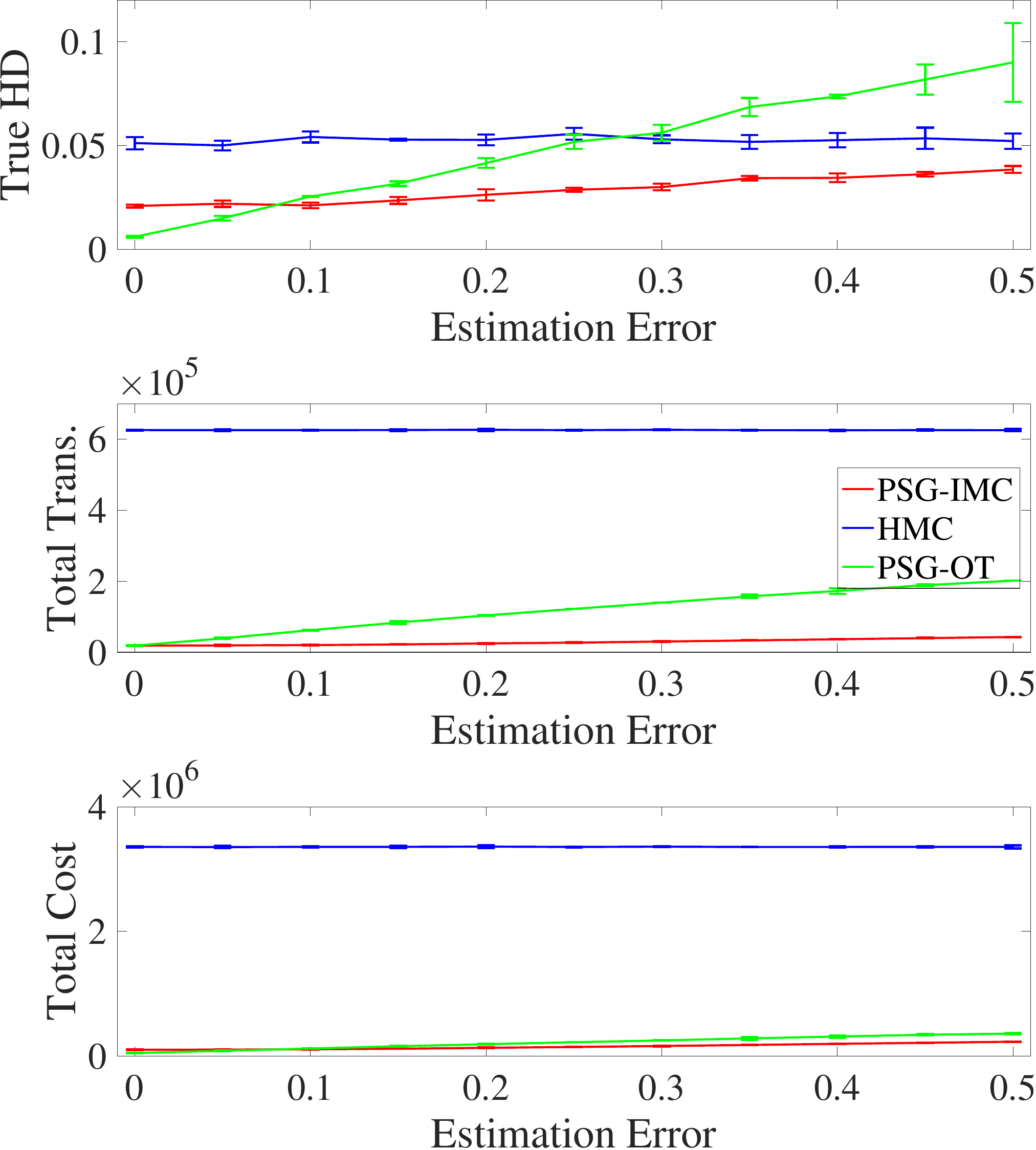} & \includegraphics[height=1.8in]{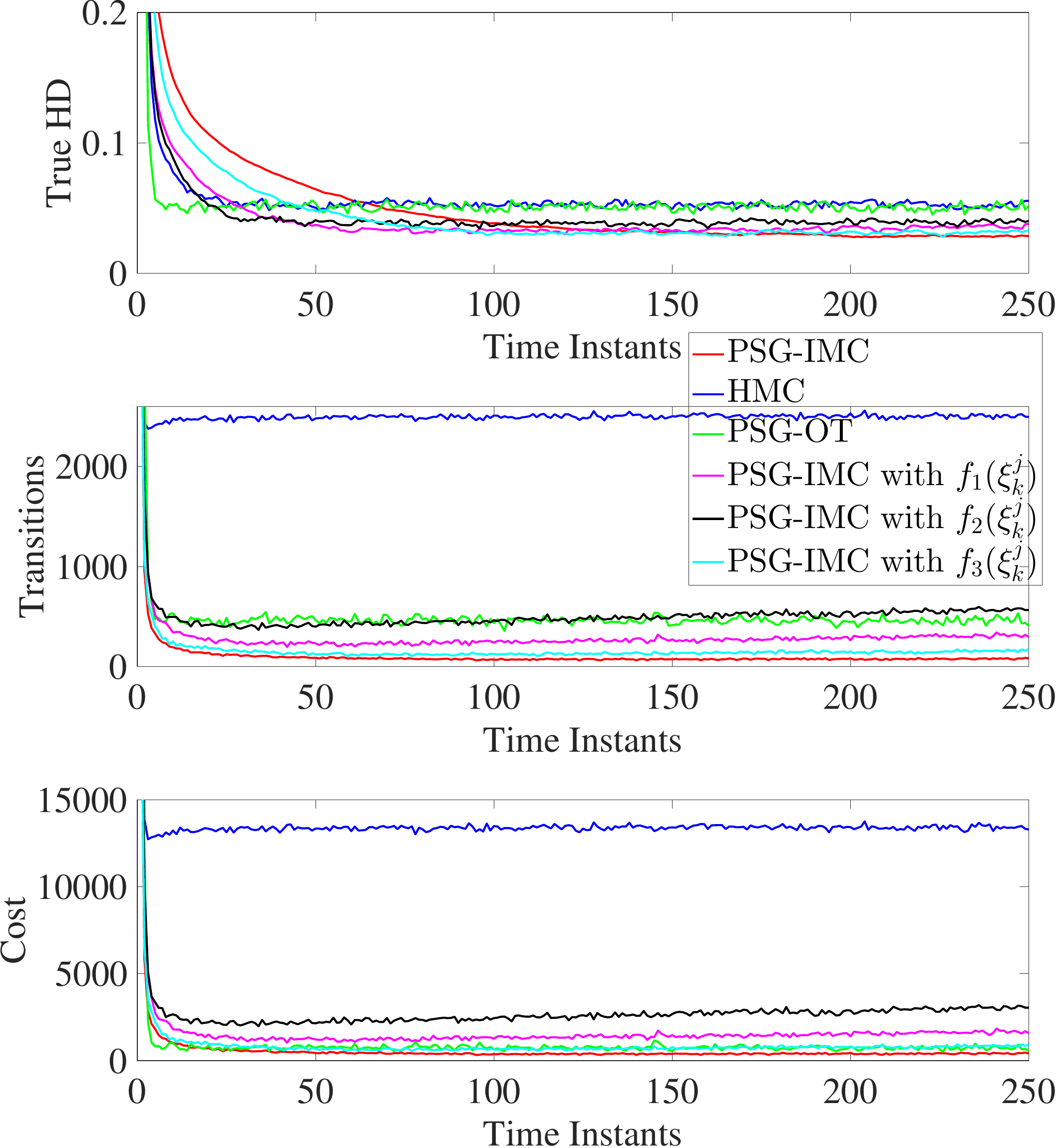}\tabularnewline
(a) & (b)\tabularnewline
\end{tabular}
\par\end{centering}
\caption{\textbf{(a)} The estimation error is varied from $0.0$ to $0.5$.
The performance of the three algorithms, along with $1\sigma$ error-bars,
are shown for the true HD $D_{H}(\boldsymbol{\Theta},\boldsymbol{\mathcal{\mu}}_{250}^{\star})$
between the actual swarm distribution after $250$ time instants and
the desired formation, the cumulative number of transitions in $250$
time instants, and the total cost incurred by all the agents in $250$
time instants. \textbf{(b)} The cumulative results of the three algorithms
and alternative functions for $\xi_{k}^{j}$ are shown, with an estimation
error $\epsilon_{\mathrm{est}}=0.25$. \label{fig:estimation-error-plot-compare}}
\end{figure}

During each time instant, each agent receives the feedback of the
current swarm distribution $\boldsymbol{\mathcal{\mu}}_{k}^{j}$ with
an estimation error $\epsilon_{\mathrm{est}}$. The cumulative results
of Monte Carlo simulations are shown in Fig.~\ref{fig:estimation-error-plot-compare}(a).
The PSG\textendash OT algorithm performs slightly better than the
PSG\textendash IMC algorithm in the absence of an estimation error
($\epsilon_{\mathrm{est}}=0.0$), but such a situation does not arise
in practical scenarios. The PSG\textendash OT algorithm's true HD
between the swarm distribution $\boldsymbol{\mathcal{\mu}}_{250}^{\star}$
and the desired formation $\boldsymbol{\Theta}$ drops precipitously
in the presence of an estimation error and it performs worse than
the open-loop HMC-based algorithm when $\epsilon_{\mathrm{est}}\geq0.25$.
On the other hand, the PSG\textendash IMC algorithm works reliably
well for all estimation errors and performs much better than the other
two algorithms. Thus, the PSG\textendash IMC algorithm can tolerate
large estimation errors in the feedback of the current swarm distribution.

The cumulative results for the three algorithms are shown in Fig.~\ref{fig:estimation-error-plot-compare}(b),
where the estimation error $\epsilon_{\mathrm{est}}$ is equal to
$0.25$. Compared to the HMC and PSG\textendash OT algorithms, the
PSG\textendash IMC algorithm achieves a smaller convergence error
with fewer transitions. The results of a few alternative functions
for $\xi_{k}^{j}$ are also shown in Fig.~\ref{fig:estimation-error-plot-compare}(b)
(see Remark~\ref{rem:Alternative_func}). The two functions $f_{1}(\xi_{k}^{j})=\tanh(\pi\xi_{k}^{j})$
and $f_{2}(\xi_{k}^{j})=\sin\left(\cos^{-1}(1-\xi_{k}^{j})\right)$
are always larger than $\xi_{k}^{j}$. The sigmoid function $f_{3}(\xi_{k}^{j})=\left(\xi_{k}^{j}+0.1\sin(2\pi\xi_{k}^{j})\right)$
is larger than $\xi_{k}^{j}$ when $\xi_{k}^{j}<0.5$. Fig.~\ref{fig:estimation-error-plot-compare}(b)
shows that the rate of convergence increases with these functions,
but there is also a corresponding increase in the number of transitions.
The collision-free motion of the agents, where the estimation error
$\epsilon_{\mathrm{est}}=0.25$ and the agents use the Voronoi-based
lower-level algorithm in Remark~\ref{rem:collision-free-motion},
is shown in the supplementary video (SV4).

\begin{figure}[h]
\begin{centering}
\includegraphics[bb=130bp 0bp 870bp 540bp,clip,width=3in]{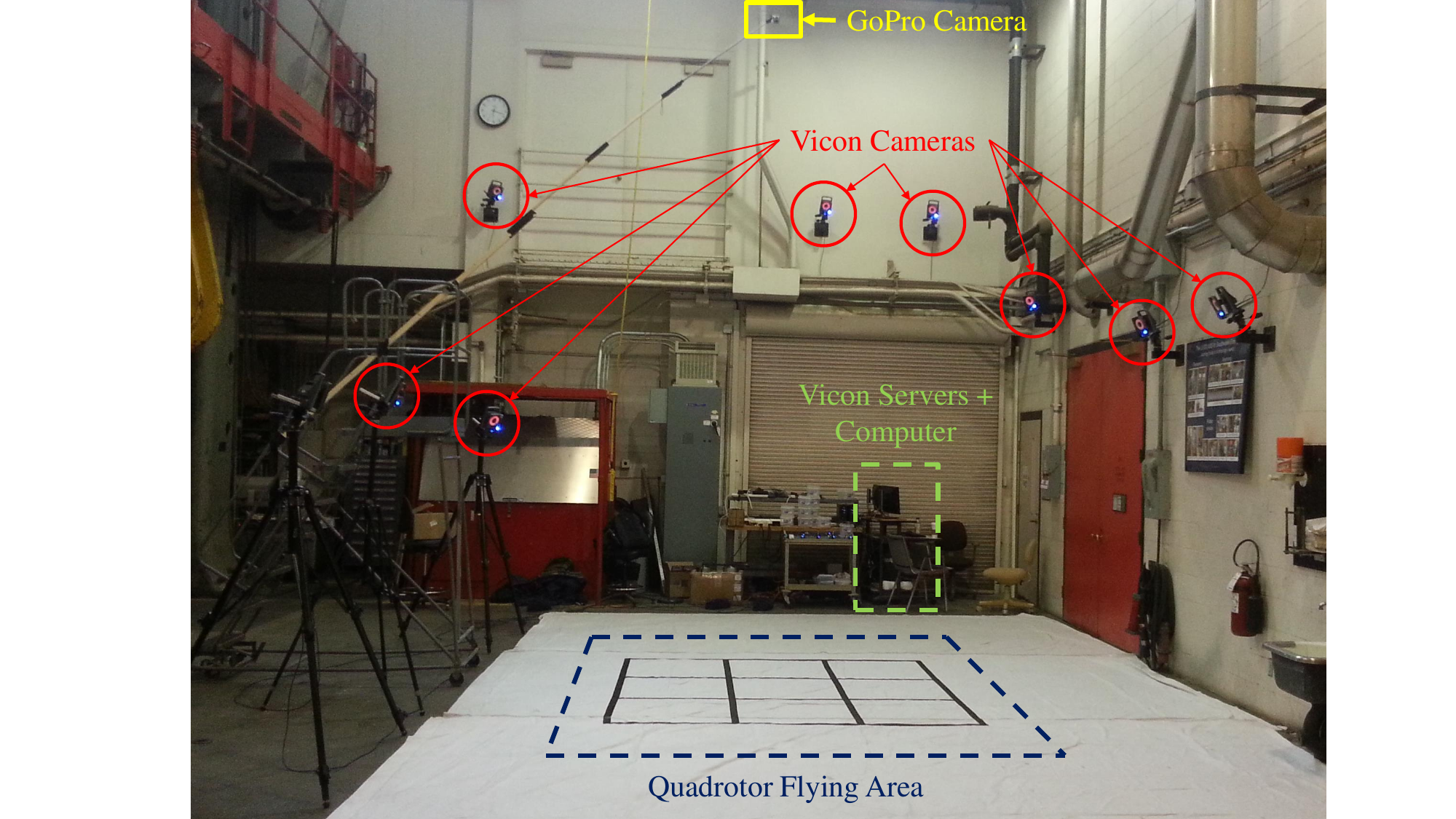}
\par\end{centering}
\caption{The experimental setup is shown. Notice the $3\times3$ grid in the
quadrotor flying area. \label{fig:experimental-setup}}
\end{figure}

\subsection{Experimental Results for Shape Formation to Demonstrate Real-Time
Execution \label{subsec:Experimental-Validation}}

In this subsection, we show that the PSG\textendash IMC algorithm
and the lower-level algorithm in Remark~\ref{rem:collision-free-motion}
can be executed in real-time using quadrotors. The experimental setup
is shown in Fig.~\ref{fig:experimental-setup} and described in \cite{Ref:Morgan15_SATO,Ref:Giri15}.
A $3\times3$ grid is placed on the ground where the quadrotor experiments
are performed. The quadrotors are tracked using the Vicon motion capture
cameras and the Vicon server fuses the data from these cameras to
estimate the position of these quadrotors to $1-2$ cm accuracy. A
desktop computer receives the position information from the Vicon
server and executes the PSG\textendash IMC algorithm for each agent
in a virtually distributed manner, i.e., each quadrotor's computations
are performed by an independent thread on the desktop so that the
experimental setup mimics a distributed system. The trajectories computed
by each quadrotor's thread are then communicated to that quadrotor.
Finally, each quadrotor tracks its computed trajectory using an onboard
nonlinear control law \cite{Ref:Morgan15_SATO,Ref:Giri15}.

\begin{figure}[t]
\begin{centering}
\begin{tabular}{c}
\includegraphics[width=3in]{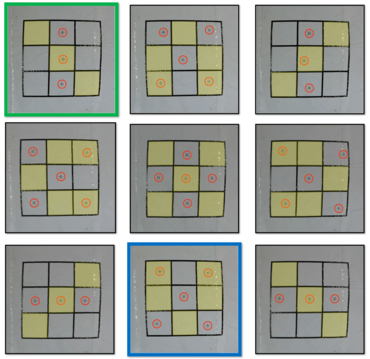}\tabularnewline
(a) Initial position\tabularnewline
\includegraphics[width=3in]{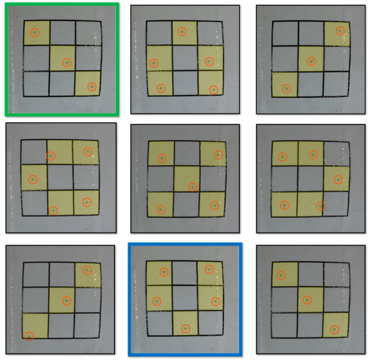}\tabularnewline
(b) Final position \tabularnewline
\end{tabular}
\par\end{centering}
\caption{Nine different experiments are shows where three or five quadrotors
execute the PSG\textendash IMC algorithm is real-time to achieve the
desired formations shown in yellow. The quadrotors are encircled in
red. See the supplementary video (SV5). \label{fig:Quadrotors-PSGIMC}}
\end{figure}

We first present nine different experiments using three or five quadrotors.
The desired formation shape for these experiments are shown in Fig.~\ref{fig:Quadrotors-PSGIMC}.
In these experiments, the time step of the PSG\textendash IMC algorithm
is $9$ seconds and the time step for the lower-level guidance algorithm
in Remark~\ref{rem:collision-free-motion} is $3$ seconds. As shown
in the supplementary video (SV5), the Crazyflie quadrotors first take
off from the ground and climb to $1$ m altitude. Thereafter, the
PSG\textendash IMC algorithm is switched on, and the quadrotors achieve
the desired formation shape within a few time instants. The quadrotor
then land inside their selected bins. Note that these exists some
parallax error in the video (SV5) because the grid is marked on the
ground, the quadrotors are flying at $1$ m altitude, and the camera
is located directly above the central square at $5$ m height. This
parallax error vanishes when the quadrotors land and the desired formation
shape is clearly visible in the end of video (SV5). The quadrotors
experience measurement errors, actuator errors, and inter-quadrotor
aerodynamic coupling due to downwash. In addition, the quadrotors
experience intermittent loss of visibility due to environmental noises
and intermittent communication losses. These experiments show that
the PSG\textendash IMC algorithm can be implemented in real-time to
achieve a variety of desired formation shapes and robustly adapt to
real-world disturbance sources.

\begin{figure}[t]
\begin{centering}
\begin{tabular}{c}
\includegraphics[clip,width=3in]{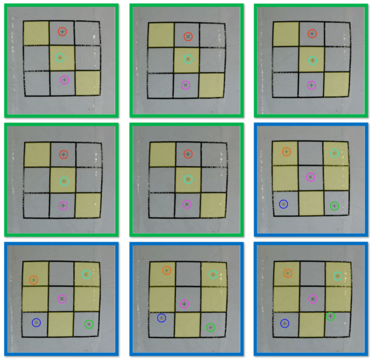}\tabularnewline
(a) Initial position\tabularnewline
\includegraphics[clip,width=3in]{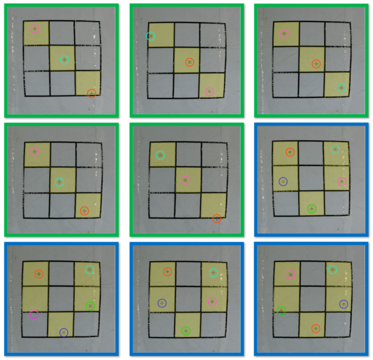}\tabularnewline
(b) Final position \tabularnewline
\end{tabular}
\par\end{centering}
\caption{Repetitions each of the same two experiments are shown. The desired
formations are shown in yellow and the quadrotors are encircled in
different colors. The actual trajectories of the quadrotors in the
these repetitions are different. See the supplementary video (SV6).
\label{fig:Quadrotors-PSGIMC-repeat}}
\end{figure}

A key feature of the PSG\textendash IMC algorithm is that each agent
probabilistically selects the bin that it transitions to. We demonstrate
this property using two sets of experiments in Fig.~\ref{fig:Quadrotors-PSGIMC-repeat}.
In the first experiment, five quadrotors start from the same initial
condition and reach the same desired formation highlighted in blue
in Fig.~\ref{fig:Quadrotors-PSGIMC}. But in each of the four experimental
runs, the actual trajectory of each quadrotor is significantly different
as shown in Fig.~\ref{fig:Quadrotors-PSGIMC-repeat}. Similarly,
in the second experiment, five experimental runs are shown where three
quadrotors reach the desired formation highlighted in green in Fig.~\ref{fig:Quadrotors-PSGIMC}.
These repetitive experiments show that the quadrotors select different
bins during different runs due to the probabilistic nature of the
PSG\textendash IMC algorithm. 

\subsection{Numerical Simulations for Area Exploration \label{subsec:Numerical-Simulation-goal searching}}

\begin{figure}[h]
\begin{centering}
\begin{tabular}{cc}
\includegraphics[height=1.7in]{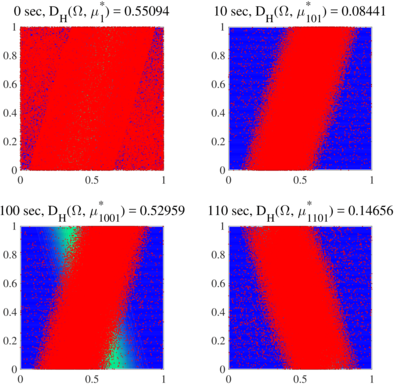} & \includegraphics[height=1.7in]{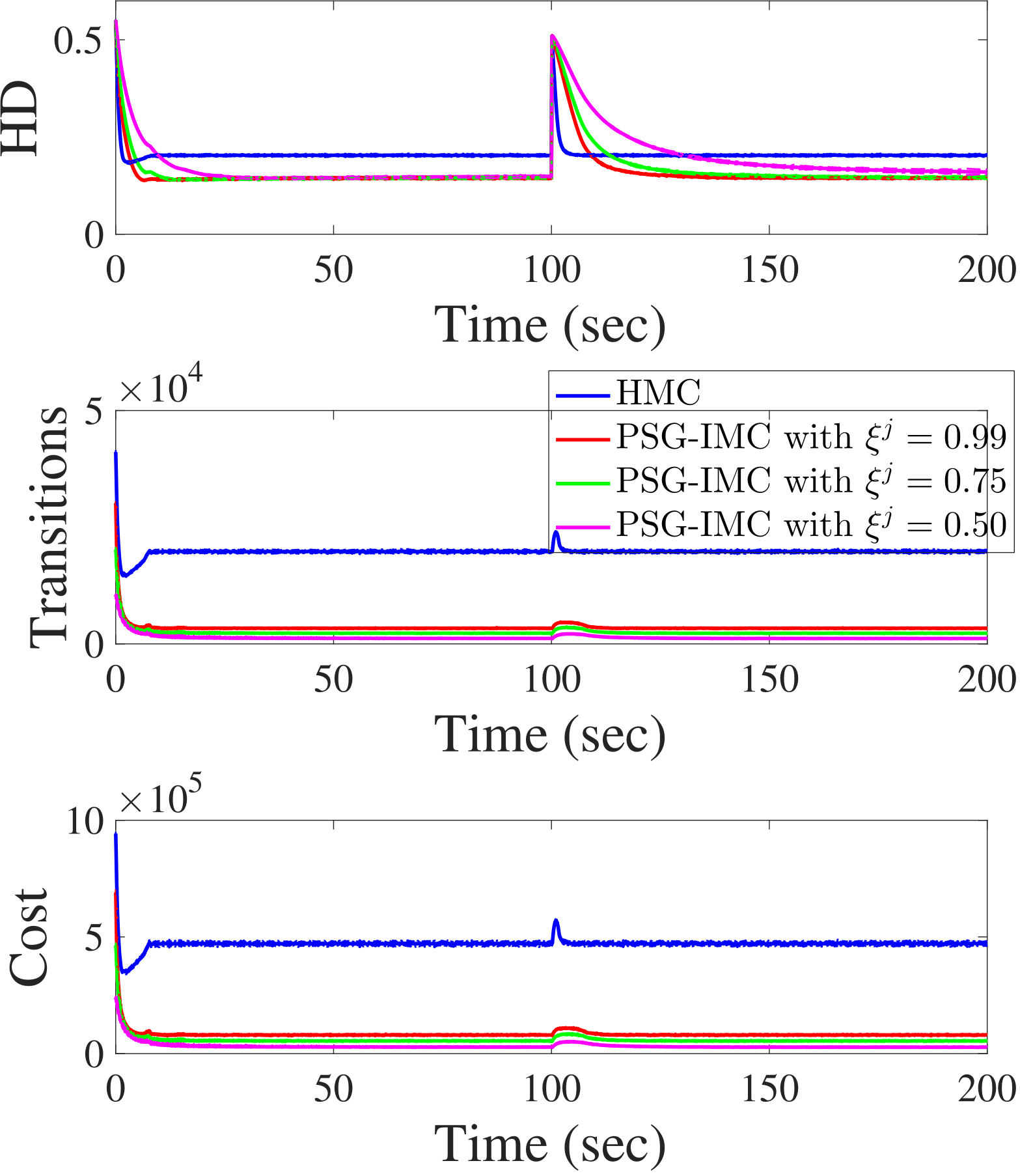}\tabularnewline
(a) & (b)\tabularnewline
\end{tabular}
\par\end{centering}
\caption{\textbf{(a)} These plots show the swarm distribution of $10^{5}$
agents (in red) and the unknown target distribution (background contour
plot), in a sample run of the Monte Carlo simulation. Starting from
a uniform distribution, the swarm converges to the unknown target
distribution. After $100$ sec, the unknown target distribution is
suddenly changed and the agents reconfigure to this new target distribution.
See the supplementary video (SV6). \textbf{(b)} The cumulative results
of $10$ Monte Carlo simulations are shown. The discontinuity or change
after $100$ sec is because of the sudden change in the unknown target
distribution. \label{fig:Swarm-goal-searching-MC}}
\end{figure}

In this numerical example, a swarm of $10^{5}$ agents use the PSG\textendash IMC
algorithm for area exploration to attain the unknown target distribution.
The physical space $[0,1]\times[0,1]$ is partitioned into $100\times100$
bins and the time step $\Delta=0.1$ sec. The unknown target distribution
$\boldsymbol{\Omega}_{1}$ for the first $100$ sec is given by the
PMF representation of the multivariate normal distribution $\boldsymbol{\mathcal{N}}\left(\left[\begin{smallmatrix}0.5 & 0.5\end{smallmatrix}\right],\left[\begin{smallmatrix}0.1 & 0.3\\
0.3 & 1.0
\end{smallmatrix}\right]\right)$, as shown in the background contour plots in Fig.~\ref{fig:Swarm-goal-searching-MC}(a).
Similarly, the unknown target distribution $\boldsymbol{\Omega}_{2}$
for the next $100$ sec is given by $\boldsymbol{\mathcal{N}}\left(\left[\begin{smallmatrix}0.5 & 0.5\end{smallmatrix}\right],\left[\begin{smallmatrix}0.1 & -0.3\\
-0.3 & 1.0
\end{smallmatrix}\right]\right)$. Here we use the constants $\tau^{j}=2.5\times10^{-3}$, and $\beta^{j}=200$
in (\ref{eq:eta_ki-goal}). 

The cumulative results of $10$ Monte Carlo simulations for different
values of $\xi^{j}$ are shown in Fig.~\ref{fig:Swarm-goal-searching-MC}(b).
A results for the HMC-based area exploration algorithm are also shown
in Fig.~\ref{fig:Swarm-goal-searching-MC}(b). Compared to the HMC-based
algorithm, the PSG\textendash IMC algorithm provides approximately
$1.5$ times improvement in HD, $6$ times reduction in the cumulative
number of transitions and the total cost incurred by the agents in
$200$ sec.

\section{Conclusion \label{sec:Conclusion}}

In this paper, we have presented a new distributed control algorithm
for large-scale swarms by using the Eulerian framework. Our highly
scalable, robust, and versatile PSG\textendash IMC algorithm ensures
that a swarm converges to the desired formation or the unknown target
distribution. 

In our PSG\textendash IMC algorithm, time-inhomogeneous Markov matrices
with a desired stationary distribution are systematically constructed
using the Hellinger-distance-based feedback gain, which incorporates
feedback from the current swarm distribution. These Markov matrices
satisfy suitable motion constraints, minimize the expected cost of
transitions at each time instant, and circumvent transitions from
bins that are deficient in the number of agents. We have also presented
a rigorous convergence analysis of the PSG\textendash IMC algorithm.
To our knowledge, this PSG\textendash IMC algorithm we have presented
in this paper is the first path planning strategy that leverages the
idea of constructing IMC, in real-time, based on feedback of the current
swarm distribution. 

We have carried out numerical simulations which show that the PSG\textendash IMC
algorithm achieves $6-16$ times reduction in total cost of transitions
and $1.5-2$ times reduction in HD, as compared to existing HMC-based
algorithms for shape formation and area exploration applications.
This is because the PSG\textendash IMC algorithm avoids undesirable
transitions, and the number of transitions at each time instant is
proportional to the HD. In the presence of estimation errors, our
PSG\textendash IMC algorithm also outperforms the optimal transport-based
algorithm, because the PMF of the predicted position of each agent
converges to the desired formation regardless of estimation errors.
We have demonstrated the robustness and computational benefits of
the PSG\textendash IMC algorithm using hardware experiments with multiple
quadrotors, where a Voronoi-based lower-level guidance and control
algorithm is used. This also provides an avenue for future research
in tightly integrating a lower-level guidance and control algorithm
with our PSG\textendash IMC algorithm. 

The PSG\textendash IMC algorithm presented in this paper can also
solve other cooperative control tasks, such as surveillance, task
allocation, and coverage, since such problems can also be cast as
shape formation or area exploration problems. We envisage that the
proposed algorithm will facilitate the development of autonomous swarm
robotic systems that are capable of performing a variety of complex
tasks, by providing a versatile, robust, and scalable path planning
strategy.

\section*{Acknowledgment}

The authors would like to thank D. Morgan, G. Subramanian, and J.
Yu for their valuable inputs. This research was carried out in part
at the Jet Propulsion Laboratory, California Institute of Technology,
under a contract with the National Aeronautics and Space Administration.
$\copyright$ 2016 California Institute of Technology. 

\section*{Appendix}

In this section, we first state some definitions and results used
in the proofs of Theorems~\ref{thm:Convergence-of-imhomo-MC} and
\ref{thm:DL1_x_Theta}, and then present the PSG\textendash OT algorithm.

\begin{definition} \label{def:primitive-matrix} \cite[pp. 3]{Ref:Seneta06}\textit{
(Primitive Matrix) }A square non-negative matrix $\boldsymbol{T}$
is said to be primitive if there exists a positive integer $k$ such
that $\boldsymbol{T}^{k}>0$. \hfill $\Box$ \end{definition}

\begin{definition} \label{def:Asymp-homo} \cite[pp. 92, 149]{Ref:Seneta06}\textit{
(Asymptotic Homogeneity) }A sequence of stochastic matrices $\boldsymbol{P}_{k},\thinspace k\geq1$
is said to be asymptotically homogeneous (with respect to $\boldsymbol{d}$)
if there exists a probability (row) vector $\boldsymbol{d}$ such
that $\lim_{k\rightarrow\infty}\boldsymbol{d}\boldsymbol{P}_{k}=\boldsymbol{d}$.
\hfill $\Box$ \end{definition}

\begin{definition} \label{def:Strong-ergodic} \cite[pp. 92, 149]{Ref:Seneta06}
\textit{(Strong Ergodicity)} The forward matrix product $\boldsymbol{U}_{T,r}:=\boldsymbol{P}_{T}\boldsymbol{P}_{T+1}\cdots\boldsymbol{P}_{r-1}$,
formed from a sequence of stochastic matrices $\boldsymbol{P}_{k},\thinspace k\geq1$,
is said to be strongly ergodic if for each $i,\thinspace\ell,\thinspace T$,
we get $\lim_{r\rightarrow\infty}\boldsymbol{U}_{T,r}[i,\ell]=\boldsymbol{v}[\ell]$,
where $\boldsymbol{v}$ is a probability vector and the element $\boldsymbol{v}[\ell]$
is independent of $i$. Therefore, $\boldsymbol{v}$ is the unique
limit vector and $\lim_{r\rightarrow\infty}\boldsymbol{U}_{T,r}=\boldsymbol{1v}$.
\hfill $\Box$ \end{definition}

\begin{theorem} \label{thm:asymp-homo-ns-strong-ergo}\cite[pp. 150]{Ref:Seneta06}
If (i) the forward matrix product $\boldsymbol{U}_{T,r}$ is primitive
and (ii) there exists $\gamma$ (independent of $k$) such that:

\begin{equation}
0<\gamma\leq\min_{i,\thinspace\ell}{}^{+}\boldsymbol{P}_{k}[i,\ell]\thinspace,\label{eq:condition_C}
\end{equation}
where $\min{}^{+}$ refers to the minimum of the positive elements,
then the asymptotic homogeneity of $\boldsymbol{P}_{k}$ is necessary
and sufficient for strong ergodicity of $\boldsymbol{U}_{T,r}$. \end{theorem}

\begin{theorem} \label{thm:limit-of-strong-ergo}\cite[pp. 149]{Ref:Seneta06}
Let $\boldsymbol{e}_{k}$ be the unique stationary distribution vector
of the matrix $\boldsymbol{P}_{k}$ (i.e., $\boldsymbol{e}_{k}\boldsymbol{P}_{k}=\boldsymbol{e}_{k}$).
If (i) all $\boldsymbol{P}_{k},\thinspace k\geq1$ are irreducible
and (ii) there exists $\gamma$ (independent of $k$) such that (\ref{eq:condition_C})
is satisfied, then asymptotic homogeneity of $\boldsymbol{P}_{k}$
(with respect to $\boldsymbol{d}$) is equivalent to $\lim_{k\rightarrow\infty}\boldsymbol{e}_{k}=\boldsymbol{e}$,
where $\boldsymbol{e}$ is a limit vector. Moreover, $\boldsymbol{d}=\boldsymbol{e}$.
\end{theorem}

\begin{corollary} \label{cor:limit-of-strong-ergo}\cite[pp. 150]{Ref:Seneta06}
Under the prior conditions (i) and (ii) of Theorem~\ref{thm:limit-of-strong-ergo}
and if (iii) $\boldsymbol{U}_{T,r}$ is strongly ergodic with unique
limit vector $\boldsymbol{v}$, then $\boldsymbol{v}=\boldsymbol{e}$.
\end{corollary}

\begin{theorem} \label{thm:positive-matrix} The matrix $\boldsymbol{U}_{k,k+n_{\textrm{rec}}-1}^{j}$
is a positive matrix. \end{theorem}

\begin{proof} Proof by Contradiction. Assume that $\boldsymbol{U}_{k,k+n_{\textrm{rec}}-1}^{j}[i,\ell]=0$
for some $i\not=\ell$. Then $\boldsymbol{U}_{k,r}^{j}[i,\ell]=0$
for all $r\leq k+n_{\textrm{rec}}-2$ because $\boldsymbol{U}_{k,k+n_{\textrm{rec}}-1}^{j}[i,\ell]\geq\boldsymbol{U}_{k,r}^{j}[i,\ell]\Bigl(\prod_{q=r}^{k+n_{\textrm{rec}}-2}\boldsymbol{P}_{q,\textrm{sub}}^{j}[\ell,\ell]\Bigr)$.
Since the matrix $\boldsymbol{P}_{k+n_{\textrm{rec}}-2,\textrm{sub}}^{j}$
is irreducible, there exists $s_{1}\in\{1,\ldots,n_{\textrm{rec}}\}\backslash\{i,\ell\}$
such that $\boldsymbol{P}_{k+n_{\textrm{rec}}-2,\textrm{sub}}^{j}[s_{1},\ell]>0$.
Since $\boldsymbol{U}_{k,k+n_{\textrm{rec}}-1}^{j}[i,\ell]\geq\boldsymbol{U}_{k,k+n_{\textrm{rec}}-2}^{j}[i,s_{1}]\boldsymbol{P}_{k+n_{\textrm{rec}}-2,\textrm{sub}}^{j}[s_{1},\ell]$,
consequently $\boldsymbol{U}_{k,r}^{j}[i,s_{1}]=0$ for all $r\leq k+n_{\textrm{rec}}-2$.
Similarly, there exists $s_{2}\in\{1,\ldots,n_{\textrm{rec}}\}\backslash\{i,\ell,s_{1}\}$
such that either $\boldsymbol{P}_{k+n_{\textrm{rec}}-3,\textrm{sub}}^{j}[s_{2},\ell]>0$
or $\boldsymbol{P}_{k+n_{\textrm{rec}}-3,\textrm{sub}}^{j}[s_{2},s_{1}]>0$.
Therefore $\boldsymbol{U}_{k,r}^{j}[i,s_{2}]=0$ for all $r\leq k+n_{\textrm{rec}}-3$.
Continuing this argument till the $k^{\textrm{th}}$ time instant,
we see that if $\boldsymbol{U}_{k,k+n_{\textrm{rec}}-1}^{j}[i,\ell]=0$,
then $\boldsymbol{P}_{k,\textrm{sub}}^{j}[i,s]=0$ for all $s\in\{1,\ldots,n_{\textrm{rec}}\}\backslash\{i\}$.
But this is a contradiction since $\boldsymbol{P}_{k,\textrm{sub}}^{j}$
is irreducible. \hfill$\blacksquare$ \end{proof}

\begin{remark} \label{rem:PSG-OT-algo} \cite{Ref:Bandyopadhyay14MSC}
\textit{(PSG\textendash OT Algorithm)} The cost function $\boldsymbol{C}_{k}$
from Definition~\ref{def:cost-function} is first modified to capture
motion constraints, i.e., $\tilde{\boldsymbol{C}}_{k}[i,\ell]=\boldsymbol{C}_{k}[i,\ell]$
if $\boldsymbol{A}_{k}^{j}[i,\ell]=1$ and $\tilde{\boldsymbol{C}}_{k}[i,\ell]=C_{\max}$
otherwise, where $C_{\max}\gg\boldsymbol{C}_{k}[i,\ell]$ for all
$i,\ell$. The optimum transference plan $\boldsymbol{\Gamma}_{k}^{j}\in\mathbb{R}^{n_{\textrm{bin}}\times n_{\textrm{bin}}}$
is found using the following LP:
\begin{align}
 & \min_{\boldsymbol{\Gamma}_{k}^{j}[i,\ell]\geq0,\forall i,\ell}\sum_{i=1}^{n_{\mathrm{bin}}}\sum_{\ell=1}^{n_{\mathrm{bin}}}\tilde{\boldsymbol{C}}_{k}[i,\ell]\boldsymbol{\Gamma}_{k}^{j}[i,\ell]\thinspace,\label{eq:OTP_LP}\\
\textrm{sub. to } & \mathrm{(i)}\sum_{\ell=1}^{n_{\mathrm{bin}}}\boldsymbol{\Gamma}_{k}^{j}[i,\ell]=\boldsymbol{\mu}_{k}^{j}[i],\forall i,\thinspace\mathrm{(ii)}\sum_{i=1}^{n_{\mathrm{bin}}}\boldsymbol{\Gamma}_{k}^{j}[i,\ell]=\boldsymbol{\Theta}[\ell],\forall\ell.\nonumber 
\end{align}
Note that the feedback of the current swarm distribution $\boldsymbol{\mu}_{k}^{j}$
directly enters (\ref{eq:OTP_LP}), hence $\boldsymbol{\Gamma}_{k}^{j}$
is sensitive to estimation errors. The matrix $\boldsymbol{\Gamma}_{k}^{j}$
is not a Markov matrix because it is not row stochastic. \hfill$\Box$
\end{remark}

\begin{remark} \label{rem:consensus-distributed-estimation} \cite{Ref:Bandyopadhyay14_ACC_BCF}
\textit{(Distributed Estimation of $\boldsymbol{\mu}_{k}^{\star}$)}
Let the probability vector $\hat{\boldsymbol{\mu}}_{k,\nu}^{j}\in\mathbb{R}^{n_{\textrm{bin}}}$
represent the $j^{\textrm{th}}$ agent's estimate of the current swarm
distribution during the $\nu^{\textrm{th}}$ consensus loop at the
$k^{\textrm{th}}$ time step. During each consensus loop, the agents
recursively combine their local estimates with their neighboring agents
as follows:
\begin{align}
\hat{\boldsymbol{\mu}}_{k,\nu+1}^{j} & =\sum_{\ell\in\mathcal{J}_{k}^{j}}a_{k}^{\ell j}\hat{\boldsymbol{\mu}}_{k,\nu}^{j},\qquad\forall\nu\in\mathbb{N},\label{eq:agreement_equation}
\end{align}
where $\mathcal{J}_{k}^{j}$ is the set of inclusive neighbors of
the $j^{\textrm{th}}$ agent, $\sum_{\ell\in\mathcal{J}_{k}^{j}}a_{k}^{\ell j}=1$,
and $\hat{\boldsymbol{\mu}}_{k,1}^{j}=\boldsymbol{r}_{k}^{j}$ is
the initial local estimate of each agent. 

Under Assumption \ref{assump:Strongly-connected-communication}, the
matrix $A_{k}$, with entries $A_{k}[\ell,j]=a_{k}^{\ell j}$, is
irreducible. Distributed algorithms in \cite{Ref:Boyd04,Ref:Cortes12}
can ensure that the matrix $A_{k}$ is balanced. Then each agents
local estimate $\hat{\boldsymbol{\mu}}_{k,\nu}^{j}$ globally exponentially
converges to $\boldsymbol{\mu}_{k}^{\star}$ pointwise with a rate
faster or equal to the second-largest singular value of $A_{k}$ (i.e.,
$\sigma_{m_{k}-1}(A_{k})$). For some $\varepsilon_{\mathrm{cons}}>0$,
if the number of consensus loops within each consensus stage $n_{\mathrm{loop}}\geq\left\lceil \frac{\ln\left(\varepsilon_{\mathrm{cons}}/2m\right)}{\ln\sigma_{m_{k}-1}(A_{k})}\right\rceil $;
then the convergence error is bounded by $\sum_{j=1}^{m_{k}}D_{\mathcal{L}_{1}}(\boldsymbol{\mathcal{\mu}}_{k}^{\star},\hat{\boldsymbol{\mu}}_{k,n_{\textrm{loop}}}^{j})\leq\varepsilon_{\mathrm{cons}}$.
The $j^{\textrm{th}}$ agent's estimate of the current swarm distribution
at the $k^{\textrm{th}}$ time instant is given by $\boldsymbol{\mathcal{\mu}}_{k}^{j}=\hat{\boldsymbol{\mu}}_{k,n_{\textrm{loop}}}^{j}$.
\hfill$\Box$ \end{remark}

\begin{remark}\label{rem:Inverse-Transform-Sampling} \cite{Ref:Devroye86}
\textit{(Inverse Transform Sampling)} This a well-known sampling technique
for generating samples at random from a given PMF over the set of
bins. The key steps are: (1) Generate a uniform random number $z$
in the interval $[0,\thinspace1]$. (2) Represent the PMF as a cumulative
distribution function (CDF). (3) Find bin $B[i]$ such that CDF of
bins up to (but not including) bin $B[i]$ is less than $z$; and
the CDF of bin $B[i]$ is greater than or equal to $z$. Then, the
bin $B[i]$ is selected as the sample. \hfill$\Box$ \end{remark}

\begin{figure}[h]
\begin{centering}
\begin{tabular}{ccc}
\includegraphics[bb=10bp 45bp 323bp 210bp,clip,width=1in]{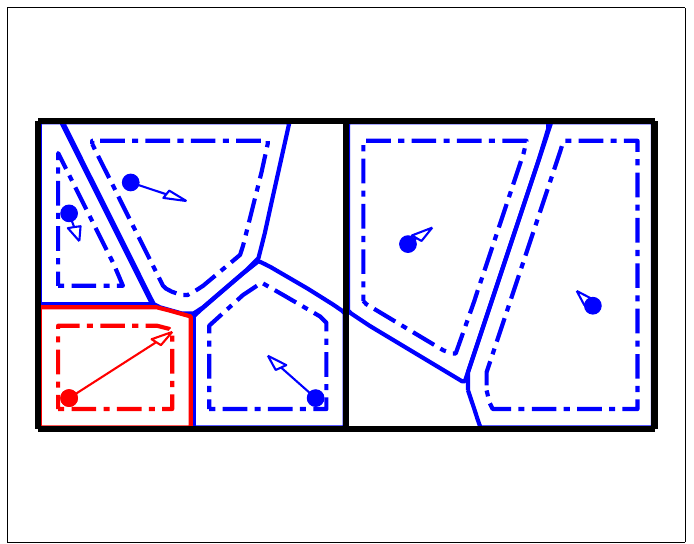} & \includegraphics[bb=10bp 45bp 323bp 210bp,clip,width=1in]{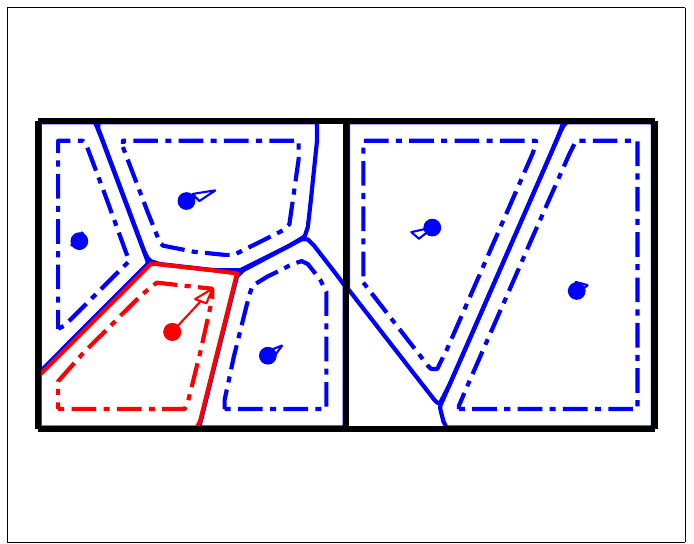} & \includegraphics[bb=10bp 45bp 323bp 210bp,clip,width=1in]{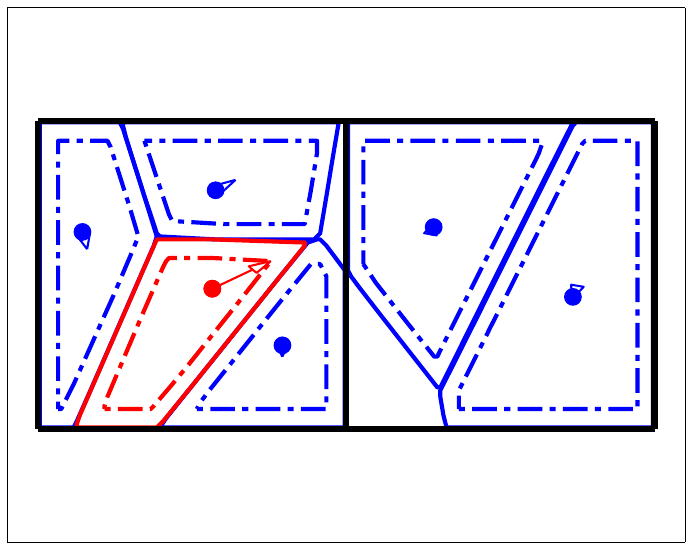}\tabularnewline
(a) & (b) & (c)\tabularnewline
\includegraphics[bb=10bp 45bp 323bp 210bp,clip,width=1in]{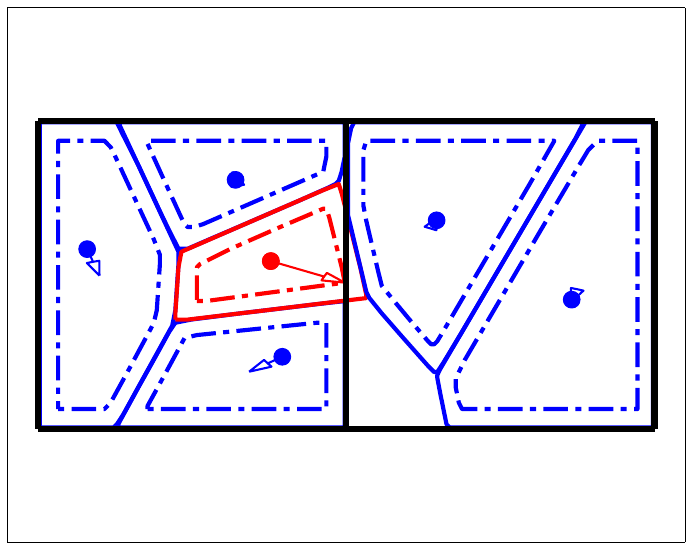} & \includegraphics[bb=10bp 45bp 323bp 210bp,clip,width=1in]{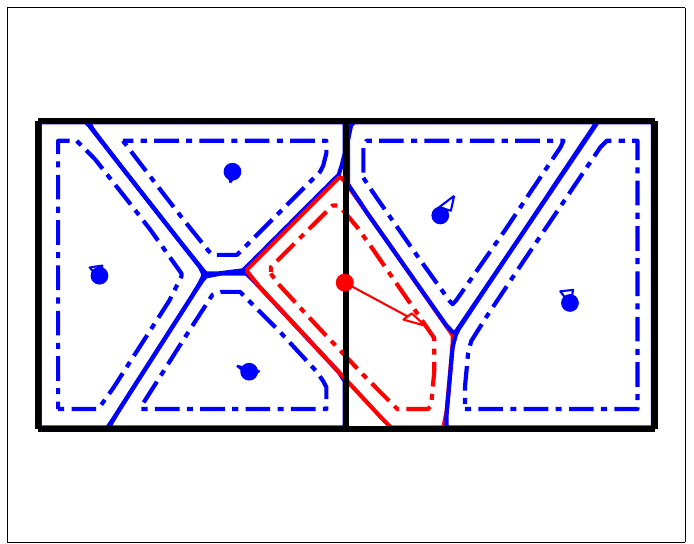} & \includegraphics[bb=10bp 45bp 323bp 210bp,clip,width=1in]{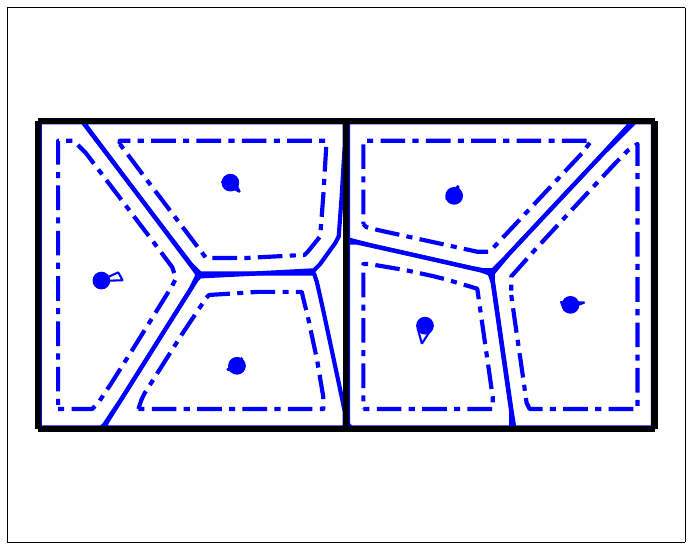}\tabularnewline
(d) & (e) & (f)\tabularnewline
\end{tabular}
\par\end{centering}
\caption{The motion of six agents in two bins (left and right bin) using the
Voronoi-based collision-free trajectory generation algorithm are shown.
The red agent goes from the left bin to the right bin. The remaining
blue agents stay in their present bin. The Voronoi sets of all the
agents along with their trajectories (denoted using arrows) are also
shown.\label{fig:Voronoi-motion}}
\end{figure}

\begin{remark}\label{rem:collision-free-motion} \cite{Ref:Bandyopadhyay14MSC}
\textit{(Voronoi-based Collision-free Motion to Selected Bin)} As
shown in Fig.~\ref{fig:Voronoi-motion}, the agents first generate
their Voronoi partitions by communicating with their neighboring agents
and considering nearby stationary obstacles. The agents that need
to transition to another bin move to the location in their Voronoi
partition that is closest to their target bin, while leaving buffer
distance for collision avoidance. The agents that want to remain in
their present bin move to the centroid of the region where their Voronoi
partition overlaps with their present bin. This results in a collision-free
trajectory for each agent. \hfill$\Box$ \end{remark}

\bibliographystyle{IEEEtran}
\bibliography{SapBib}

\end{document}